\documentclass[10pt]{article}
\usepackage{amssymb}
\usepackage[english]{babel}

\newcommand{\dem}{\noindent{\bf Proof. }}

\newtheorem{theor}{Theorem}[section]

\newtheorem{prop}[theor]{Proposition}
\newtheorem{coro}[theor]{Corollary}
\newtheorem{fcoro}[theor]{Corollaire}
\newtheorem{lema}[theor]{Lemma}
\newtheorem{rema}{Remark}

\newtheorem{defi}{Definition}[section]

\begin{Large}                \end{Large}

\let \x = \backslash

\let \n = \noindent
\author{  Amir Baklouti and Sa\"{\i}d Benayadi  \\\\ Universit\'e Paul Verlaine-Metz, LMAM, CNRS UMR 7122, \\Ile
du Saulcy, F-57045 Metz cedex, France}
\title{Pseudo-euclidean Jordan algebras}
\date{ }
\begin{document}
\maketitle
\begin{abstract}
A  pseudo-euclidean Jordan algebra  is a Jordan algebra  $\frak J$  with an associative non-degenerate symmetric bilinear form
$B$. $B$ is said an associative scalar product on $\frak J$. We  study the structure of the pseudo-euclidean Jordan $\mathbb
K$-algebras (where $\mathbb K$ is a field of null characteristic) and we obtain an inductive description of these algebras in terms of double extensions and   generalized double extensions. Next, we study the symplectic pseudo-euclidean Jordan $\mathbb
K$-algebras and we give some informations on  a particular class of these algebras namely the class of symplectic Jordan-Manin Algebras. Finally, we obtain some characterizations of semi-simple Jordan algebras among the pseudo-euclidean Jordan algebras.

\vspace*{0.2cm}
\n {\it Keywords:} Jordan Algebras,Tits-Kantor-Koecher   construction, Pseudo-euclidean symplectic Jordan  algebras, symplectic quadratic Lie  algebras, representations of Jordan algebras, Jordan bialgebras, Jordan-Manin algebras, Jordan Yang Baxter equation, r-matrices, $T^*-$extension, double extensions.

\vspace*{0.1cm}

\n {\it MSC:} 17C10, 17C20, 17C50, 17C55, 17B05, 17B62.

\vspace*{0.1cm}

\n {\it Corresponding Author:} Sa\"{\i}d Benayadi,
{\it e-mail:} benayadi@univ-metz.fr
\end{abstract}

\section{Introduction}
In this paper, we consider finite dimensional algebras over a
commutative field $K$ with characteristic zero.
 Let $\frak J$ be a Jordan algebra. A bilinear form $B$ on $\frak J$ is said associative if $B$ satisfies
$B(xy,z)=B(x,yz),\,\,\,\forall x,y,z\in \frak J$. Moreover, if  $B$
is nondegenerate and symmetric, we say that $B$ is an associatif
scalar product of $\frak J$. In this case, $({\frak J}, B)$ is
called a pseudo-euclidean Jordan algebra. It is well Known that,
the semi-simple Jordan algebras are pseudo-euclidean (see
$\cite{far},\cite{Jacob}, \cite{scha}$). But there are many nilpotent Jordan algebras which
are pseudo-euclidean (see the second section  of this paper and $\cite{Bord}$).

It is well known that, to any pseudo-euclidean unital Jordan algebra, one
can apply the Tits Kantor Koecher   construction $(TKK$
construction) to obtain  a quadratic Lie algebra
$(\cite{koe1},\,\cite{koe2},\,\cite{Tits})$. In section $\ref{gen}$, we slightly
modify this construction in order to obtain  a $3$-graded
quadratic Lie algebras ${\cal G}={\cal G}_{-1}\oplus{\cal
G}_0\oplus{\cal G}_1,$ starting from the pseudo-euclidean Jordan
algebras $({\cal G},B)$ which are not necessarily unital. We also
call this modification the TKK construction. Recall that  the  Lie algebra $\cal
G$ of Lie group $\frak G$ which admits a bi-invariant
pseudo-Riemannian structure is quadratic (i.e. $\cal G$ is endowed
with a symmetric non degenerate invariant (or associative) bilinear form $B$ . Conversely, any connected Lie
group whose Lie algebra ${\cal G}$ is quadratic is endowed with
bi-invariant pseudo-Riemannian structure $\cite{nei}$. In $\cite{med},$ A. Medina
and Ph. Revoy have introduced the notion of double extension and
used this notion to give an inductive classification
 of quadratic Lie algebras.

 The principal objective    of this paper  is to  study    structures  of
pseudo-euclidean Jordan algebras and to give  inductive descriptions of   these algebras. In $\cite {Bord} $, it is
 shown  that any pseudo-euclidean Jordan algebra is the
orthogonal direct sum of irreducible ideals. So  the study of  pseudo-euclidean Jordan algebras is reduced to study the irreducible ones. In order to obtain informations and inductive descriptions of the irreducible pseudo-euclidean Jordan algebras, we introduce in section $\ref {ext} $ some notions of double extensions namely the double extension of pseudo-euclidean Jordan algebras and  the generalized double extension of pseudo-euclidean Jordan algebras by the one dimensional Jordan algebra with zero product. Next, in  the section $\ref {desc} $, we show that any pseudo-euclidean Jordan algebra can be obtained  from a finite number of elements of $\cal U$ (where $\cal U$ is the set constituted by the  algebra $ \{0 \} $, the one dimensional
 Jordan algebra with zero product and the simple Jordan algebras)   by a finite sequence of orthogonal
 direct sums of pseudo-euclidean Jordan algebras and / or generalized double extensions by one dimensional Jordan algebra with zero product and / or double extensions by simple Jordan algebras. recall that the classification of simple Jordan algebra is we
  description is well known (see. $\cite{bert}$, $\cite{far},\,\,\cite{Jacob},\,\,\cite{mac}$. In particular, we can construct all
nilpotent pseudo-euclidean Jordan algebras from the null algebra  $\{0\}$ and the one dimensional Jordan algebra with zero product using the generalized double extension by the one dimensional Jordan algebra with zero product.
We shall construct, in section $\ref{5}$, all nilpotent pseudo-euclidean Jordan algebras of dimension less then or equal
to five. The list of these algebras shows that all nilpotent pseudo-euclidean Jordan algebras of dimension  less then or equal
to four are associative and that the smallest dimension of non associative nilpotent  pseudo-euclidean Jordan algebra is five.

Sections  $\ref {big} $ and $\ref {dajps} $ are devoted to study
symplectic pseudo-euclidean Jordan algebras and  their connections
with the symplectic qudratic Lie algebras. A  pseudo-euclidean Jordan
algebra $ ({\frak J}, B) $ is said symplectic if it is endowed
with a skew symmetric nondegenerate bilinear form $\omega:{\frak
J}\times{\frak J}\longrightarrow\mathbb K$ which satisfies
\begin{eqnarray}\label{si}
\omega(xy,z)+\omega(yz,x)+\omega(zx,y)=0,\,\,\,\forall
x,y,z\in{\frak J}.
\end{eqnarray}
Such form, is said a symplectic form of ${\frak J}$. Recall that
in $\cite{zhel}$ and $\cite{zhel1}$, a Jordan algebra $\frak J$ is
symplectic if it is endowed with an antisymmetric bilinear form
not necessarily nondegenerate  $\omega:{\frak J}\times{\frak
J}\longrightarrow\mathbb K$ which satisfies $(\ref{si})$. In the
case where $\omega$ is nondegenerate, V.N Zhelyabin
$(\cite{zhel},$ $\cite{zhel1})$ has obtained some results which
similar to known results in the symplectic Lie algebras. If $({\frak J},B,\omega)$ is a symplectic pseudo-euclidean Jordan
algebra, we show that  there exist an invertible derivation $D$ of $ {\frak
J}$ satisfying $\omega(x,y)=B(D(x),y),\,\,\,\forall x,y\in {\frak
J}.$ Next we prove that ${\frak J}$ is nilpotent and  ${\cal U}=D^{-1}$  is a solution of 
the Jordan  Yang Baxter equation. We show in section  $\ref{big},$ that the symplectic form
$\omega$ may be extended to the quadratic Lie algebra obtained by
the $TKK$ construction starting from $({\frak J}, B)$, if and only
if the derivation $D$ checks the following condition:
\begin{eqnarray}
&& \hskip0.5cm\left[D,\left[ R({\frak J}),R({\frak J})\right]\right]
=\left[ R({\frak J}),R({\frak J})\right],\hskip1.5cm\label{2}
\end{eqnarray}
where $R({\frak J})=\{R_x,\,\,\, x\in{\frak J}\}$ and $R_x:{\frak
J}\longrightarrow{\frak J},\,\,\,y\longmapsto yx$. If $D$
satisfies the condition $(\ref{2})$ above, we obtain a symplectic quadratic Lie algebra from a symplectic pseudo-euclidean Jordan algebra. Moreover this symplectic quadratic Lie algebra is a $\mathbb Z_2-$graded
  symplectic quadratic Lie algebra (see Remark $\ref{sympz2}$). In section $\ref{big}$, we give examples of non
  associative symplectic pseudo-euclidean Jordan algebras such that  the quadratic Lie algebra obtained from these examples by the $TKK$ construction are symplectic and pseudo-euclidian. Let us recall that    the  Lie algebra of Lie group which admits a
bi-invariant pseudo-Riemannian metric and also a left-invariant
symplectic form is a symplectic quadratic Lie algebra. These Lie
groups are nilpotent and their geometry (and, consequently, that
of their assiciated homogeneous spaces) is very rich. In
particular, they carry two left-invariant affine structures: one
defined by the symplectic form and another which is compatible
with a left-invariant pseudo-Riemannian metric. Moreover, if the
symplectic form is viewed as a solution $r$ of the classical Yang
Baxter equation of Lie algebras (i.e. $r$ is an r-matrix), then
the Poisson-Lie tensor $\pi=r^+-r^-$ and the geometry of double
Lie groups $D(r)$ can be nicely described in $\cite{med2}.$ In addition, the symplectic quadratic Lie algebras were described by the methods of double extensions in $\cite{aub}$ and $\cite{ben2}$. Further, in $\cite{ben2}$, it is proved that every  symplectic
quadratic Lie algebra $({\cal G},B,\omega)$, over an algebraically closed fields $\mathbb K$, may be
 constructed by $T^*-$extension of nilpotent Lie algebra which admits an invertible derivation.

In order to give an inductive description of the  symplectic
pseudo-euclidean Jordan algebras, we  introduce  in section
$\ref{dajps}$  the notion of symplectic pseudo-euclidean double
extension. More precisely, we prove that every symplectic
pseudo-euclidean Jordan algebra may be constructed from the algebra $\{0\}$ by a finite number of 
symplectic pseudo-euclidean double extension of symplectic
pseudo-euclidean Jordan algebra 

In the section $\ref{mani}$, we introduce Jordan-Manin algebra $({\frak J},B)$ (i.e. $({\frak J}={\cal
U}\oplus{\cal V},B)$ where ${\cal U}$ and ${\cal V}$ are two completely isotropic subalgebras of ${\frak J}$) and the notion of
double extension of Jordan-Manin algebras in order to describe the nilpotent Jordan-Manin algebras. We will show,
 in Proposition $\ref{spe},$ that every symplectic pseudo-euclidean Jordan algebra over an algebraically  closed fields $\mathbb K$ is a symplectic Jordan-Manin algebra (i.e.
(${\frak J}={\cal U}\oplus{\cal V},B,\omega)$ where ${\cal U}$ and
${\cal V}$ two subalgebras   of ${\frak J}$ completely isotropic with respect to $B$ and $\omega$). Finally, we give an inductive description of symplectic Jordan-Manin algebras over an algebraically  closed field by using the notion of symplectic double extension of symplectic Jordan-Manin algebras.

By using some results obtained in the previous  sections, we give  in the section $9$  some new characterizations of semisimple Jordan algebras among the pseudo-euclidean Jordan algebras. These charcterizations are based on representations of Jordan algebras, operators of Casimir type and the index of a pseudo-euclidean Jordan algebra.

\section{Definitions and
preliminary results}\label{gen}
 Jordan algebra $ \frak J$ is a commutative non necessary associative algebra which satisfy:
\begin{eqnarray}\label{eq-1}
x(yx^2)=(xy)x^2 ,\,\,\,\forall x,y\in{\frak J},
\end{eqnarray}

where $R_x$ is the endomorphism of ${\frak J}$ defined by: $R_x(y):= xy=yx, \forall y \in    {\frak J}.$

This equality is equivalent to $\left[ R_x,R_{x^2}\right]=0
,\,\,\,\forall x,y\in\frak J$. In Jordan algebra, the following
equality are satisfied:
\begin{eqnarray}
&&2(x,y,zx)+(z,y,x^2)=0,\,\,\,\forall x,y,z\in\frak J \mbox{ where
}
(x,y,z)=(xy)z-x(yz),\,\,\,\forall x,y,z\in\frak J.\label{eqq}\\
&&[R_{wz},R_{x}]+[R_{zx},R_{w}]+[R_{xw},R_{z}]=0 \,\,\,\,\,\,\,\, \forall x,z,w\in
{\frak J}.\label{eq-4}\\
&&R_{xy}R_{z}-R_{x}R_{y}R_{z}+R_{zx}R_{y}-R_{y(zx)}+R_{zy}R_{x}-R_{z}R_{y}R_{x}=0\hskip 0,2cm x,y,z \in
\frak J.\label{eq-5}\\
&&\left[ R_x,\left[ R_y,R_z\right] \right] =R_{(y,x,z)}=R_{\left[ R_z,R_y\right] (x)},\,\,\,\forall x,y,z
\in\frak J.\label{eg-2}
\end{eqnarray}
 For more informations about Jordan algebra see   $\cite{far},\cite{scha}, \cite{Jacob}$ and $\cite{mac}$ .\\
\begin{defi} \label{def1}
Let $\frak{J}$ be a Jordan algebra.
\begin{enumerate}
\item[(i)]A bilinear form $B$ on $\frak{J}$ is called associative if
$  B(xy,z)=B(x,yz),\,\,\ \forall x,y,z \in \frak{J}.$
\item[(ii)] If $B$ is a nondegenerate symmetric and associative bilinear form on $\frak{J}$, we say that (${\frak J},B)$ is a pseudo-euclidean Jordan algebra and $B$ is an associatif scalar product on $\frak{J}$.
\end{enumerate}

\end{defi}
 In the following proposition, we are going to give a characterization of pseudo-euclidean Jordan algebra.
\begin{prop} 
Let  $\frak{J}$ be a Jordan  algebra. Then  $\frak{J}$  is pseudo-euclidean  if and only if  its  adjoint  representation and  co-adjoint  representation are equivalent.
\end{prop}
\dem
Let  $\frak{J}$ be a Jordan algebra, $R$ (resp.$\rho$) be the adjoint (resp.coadjoint) representation of  $\frak J$ (see the definition or a representation of Jordan algebra in Proposition $\ref{rep}$  of section 3). Recall that:
\begin{eqnarray*}
R&:&\frak{J}\longrightarrow End(\frak{J})\hskip 0,5 cm and \hskip 0,5 cm\rho:\frak{J}\longrightarrow End(\frak{J}^*) 
 \end{eqnarray*}
are defined by  $R(x)y:= R_x(y)= xy \mbox{ et } \Bigl(\rho(x)f\Bigl):=f\circ R_x,\,\, \forall x,y\in {\frak J},\forall f\in {\frak J}^*.$

Assume that $\frak{J}$ is  pseudo-euclidean Jordan algebra, then  there exists a bilinear form  $B:\frak{J}\times \frak{J}\rightarrow \mathbb{K}$ which is   symmetric, non-degenerate  and  associative. Consequently, the map $\varphi: \frak{J}\rightarrow \frak{J}^*$ defined by:  $\varphi(x):=B(x,.), \forall x \in \frak{J},$ is an isomorphism of vector spaces which verifies:
\begin{eqnarray*}
\varphi\Bigl(R_x(y)\Bigl)(z)=B(xy,z)=B(y,xz)=\Bigl(\rho(x)\varphi(y)\Bigl)(z),\,\,\forall x,y,z \in {\frak J}.  
\end{eqnarray*}
Which proves that  $\varphi\circ R_{x} =\rho(x)\circ \varphi,\,\,\,\forall x \in {\frak J} $. It follows that  the  representations $R$ and  $\rho$ are  equivalent.  

Conversely, suppose that the  representations $R$ and  $\rho$  are  equivalent, then   there exists  an isomorphism  of vector spaces $\phi:\frak{J}\rightarrow \frak{J}^*$ such that $\phi\circ R_{x} =\rho(x)\circ \phi, \forall x \in {\frak J}$. Now, we consider the map  $T:\frak{J}\times \frak{J}\rightarrow \mathbb{K}$ defined by: $T(x,y):=\phi(x)(y), \forall x, y \in \frak{J}.$ Since  $\phi$ is invertible, then  $T$ is  non-degenerate. Moreover, $T$ is associative, in fact  $T(xy,z)= \phi(xy)(z)= \phi(R(y)x)(z)= (\rho(y)(\phi(x))(z)= \phi(x)(zy)= T(x,zy)= T(x,yz), \forall x, y, z \in {\frak J}.$  It is clear that $T$ is not necessarily  symmetric. We consider the symmetric    (resp.  anti-symmetric) part $T_{s}$ (resp. $T_{a}$) of  $T$ defined by: $T_{s}(x,y)= \frac{1}{2}(T(x,y)+T(y,x)), \,\, (resp. T_{a}(x,y)= \frac{1}{2}(T(x,y)-T(y,x)), \,\, \forall  x,y \in {\frak J}.$ It is clear that $T$ is  associative if and only if  $T_{s}$ and $T_{a}$ are associative. 

Let us consider  $\frak{J}_{s}=\{x\in \frak{J} \,/\, T_{s}(x,y)=0, \,\forall y\in \frak{J}\}$ and $\frak{J}_{a}=\{x\in \frak{J}\,/\,T_{a}(x,y)=0,\,\, \,\forall y\in \frak{J}\}.$ The fact that  $T$ est non-degenerate implies that  $\frak{J}_{s}\cap \frak{J}_{a} =\{0\}.$  Let $x,y,z \in \frak{J}$, since $T_{a}(xy,z)=T_{a}(x,yz)=-T_{a}(z,xy)=-T_{a}(x,zy)$ then $T_{a}(xy,z)=0.$ It follows that  $ \frak{J}^2=\frak{J}\frak{J} $ is contained in   $\frak{J}_{a}.$    Moreover $\frak{J}_{s}$ and $\frak{J}_{a}$ are ideals of  $\frak{J}$ because $T_{s}$ and  $T_{a}$ are associative. Consequently,    $\frak{J}_{s}^2\subset \frak{J}_{s}\cap \frak{J}_{a} =\{0\}.$ Now, let   $\cal W$ be a sub-vector space of  $\frak{J}$ such that  $\frak{J}_{a} \subset {\cal W}$ and  $\frak{J}={\cal W}\oplus \frak{J}_{s}$. Consider   $F: {\frak{J}_{s}}\times {\frak{J}_{s}}\rightarrow \mathbb{K}$ be a non-degenerate  symmetric bilinear on $\frak{J}_{s}$. Since $\frak{J}_{s}^2=\{0\}, $ then $F$ est associative. Therefore the bilinear form $L: \frak{J}\times \frak{J} \rightarrow \mathbb{K}$ defined by:  $L_{\vert_{W\times W}}={T_{s}}_{\vert_{W\times W}},\,\, L_{\vert_{\frak{J}_{s}\times \frak{J}_{s}}}=F,\,\, L(W,{\frak{J}_{s}})= L({\frak{J}_{s}},W)= \{0\},$ is symmetric  non-degenerate and  associative. Then  $(\frak{J},L)$ is a pseudo-euclidean   Jordan algebra. $\Box$

\begin{defi}
\begin{enumerate}
 \item Let $(\frak J,B)$ be a pseudo-euclidean Jordan algebra. An ideal ${\cal I}$ is called nondegenerate  (resp. degenerate), if the restriction of $B$ on ${\cal I}\times {\cal I}$ is a nondegenerate  (resp. degenerate)bilinear form.
 \item A pseudo-euclidean Jordan algebra is called  $B-$irreducible, if $\frak J$
 contains  no non-trivial nondegenerate ideals.
\end{enumerate}
\end{defi}

One has the following lemma whose proof is straightforward.

\begin{lema}\label{lem1}
 Let $({\frak J},B)$ be a pseudo-euclidean Jordan algebra and ${\cal I}$ be an ideal of ${\frak J}.$ Then,
\begin{enumerate}
\item[(i)] ${\cal I}^\perp $ is an ideal of ${\frak J}$ and ${\cal I}{\cal I}^{\perp}=\{0\};$ 
 \item[(ii)]If ${\cal I}$ is nondegenrate, then ${\frak J}={\cal I}\oplus {\cal I}^\perp$ and ${\cal I}^\perp$  is a nondegenerate ideal of ${\frak J}$;
\item[(iii)]If ${\cal I}$ is semisimple, then ${\cal I}$ is nondegenerate.
\end{enumerate}
\end{lema}
\begin{prop}\label{deco}
Let $(\frak J,B) $ be a pseudo-euclidean Jordan algebra.
Then, ${\frak J}=\oplus^{r}_{i=1}{\frak J}_i,$ where for all
$1\leq i\leq r,$
\begin{enumerate}
\item[(i)] ${\frak J}_i$ is a non degenerate ideal;
\item[(ii)] ${\frak J}_i$ contains no nondegenerate ideal of ${\frak J}_i$.
\item[(iii)] For all $i\neq j$,  ${\frak J}_i$ and ${\frak J}_j$ are orthogonal.
\end{enumerate}
\end{prop}
For more details about this decomposition, see $\cite{Bord}$.
 \begin{defi}
Let ${\frak J}$ be a   algebra. \\
\begin{enumerate}
\item[(i)] $({\frak
J},{\frak J},{\frak J}):=Vect\{ (x,y,z):=(xy)z-x(yz);,\,\,\,x,y,z\in
{\frak J}\}$ is a vector sub-space of ${\frak J}$ called the associator of ${\frak J}.$
\item[(ii)]  The vector subspace  $Ann({\frak J}):=\{x\in {\frak J};\,\,\,\, xy= yx= 0,\forall y\in {\frak J}\}$ of ${\frak J}$ is called the annulator of ${\frak J}.$ 
\item[(iii)]  The vector subspace $N({\frak J}):=\{x\in {\frak J};\,\,\,\,
(x,y,z)=(y,x,z)=(y,z,x)=0, \,\, \forall y,z\in{\frak J} \}$ is called the nucleus of ${\frak J}.$
\end{enumerate}
\end{defi}

\begin{rema}

\begin{enumerate}
 \item   In the case of Jordan algebra ${\frak J}$,  the nucleus $N({\frak J})$ of ${\frak J}$    coincide with the center    

$Z({\frak J}):=\{x\in N({\frak J});\,\,\,\,xy=yx, \,\,\forall y \in{\frak J} \}$ of ${\frak J}.$
\item  Let  $(\frak J,B) $ be  a pseudo-euclidean Jordan algebra. Since $B((x,y,z),t)= B((y,x,t),z)= B((z,t,x),y)\,\,, \forall x,y,z\in{\frak J}, $ then $N({\frak J})=Z({\frak J}):=\{x\in {\frak J};\,\,\,\, (x,y,z)= 0, \,\, \forall y,z\in{\frak J} \}.$
\end{enumerate}
\end{rema}

\begin{prop}\label{cent}
Let $({\frak J},B)$ be a pseudo-euclidean Jordan algebra. Then,
 $$1.\,\,\,\Bigl(Ann({\frak J})\Bigl)^{\perp}={\frak J}^2.\hskip 2cm
2.\,\,\,\Bigl(Z({\frak J})\Bigl)^{\perp}=({\frak J},{\frak J},{\frak J})$$
\end{prop}
\dem $1.$ The fact that $B$ is associative implies that  ${\frak J}^2 \subset \Bigl(Ann({\frak J})\Bigl)^{\perp}.$ Conversely, let $x\in ({\frak J}^2)^\perp$ and $y\in  {\frak J} $. Then  $B(xy,z)=B(x,yz)=0, \forall z \in {\frak J}.$ Since $B$ is non-degenerate, then $xy=0.$ Thus,  $x\in Ann({\frak J})$. Hence, $({\frak
J}^2)^\perp\subset Ann({\frak J}).$ We conclude that $\Bigl(Ann({\frak J})\Bigl)^{\perp}= {\frak J}^2.$

$2.$ Let $x,y,z\in{\frak J}$ and  let $u\in Z({\frak J}), \,$ $B(X,u)=B((xy)z-x(yz),u)=B(x,y(zu))-B(x,(yz)u)=B(x,y(zu)-(yz)u)=0.$ Which proves that  $({\frak J},{\frak J},{\frak J})\subset \Bigl(Z({\frak J})\Bigl)^{\perp}$. Conversely, let $x\in ({\frak
J},{\frak J},{\frak J})^{\perp}$.
$B((x,y,z),u)=B(x,(u,z,y))=0 \, \mbox{and}\, B((y,x,z),u)=B(x,(z,u,y))=0, \,  \forall y,z, u \in{\frak J}$. Since $B$ is nondegenerate, then $(x,y,z)= (y,x,z)= 0$. Hence, $x\in Z({\frak J})$. Therefore, $ ({\frak J},{\frak J},{\frak J})^{\perp}\subset \Bigl(Z({\frak J})\Bigl)$. Consequently, $\Bigl(Z({\frak J})\Bigl)^{\perp}=({\frak J},{\frak J},{\frak J}).$ $\Box$

\begin{coro}
Let  $({\frak J},B)$ be  a  Jordan algebra. If ${\frak J}\neq\{0\},$ then $({\frak J},{\frak J},{\frak J})\neq {\frak J}.$
\end{coro}
%%%%%%%%%%%%%%%%%%%%%%%%%%%%%%%%%%%
%%%%%%%%%%%%%%%%%%%%%%%%%%%%%%ùits Albert form   ${\frak A}
\dem Suppose that $({\frak J},{\frak J},{\frak J})= {\frak J}.$ Let $Rad({\frak J})$ be the radical of ${\frak J}$, $S$
be a semi-simple subalgebra of ${\frak J}$ such that ${\frak
J}= S \oplus Rad({\frak J}).$  The fact that   $({\frak J},{\frak J},{\frak J})= {\frak J}$ implies that  $(S,S,S)=S$. Since, by Proposition \ref{cent}, $(S,{\frak A})$ is a pseudo-euclidean Jordan algebra (where ${\frak A}$ is the Albert form of $S$), then  $Z(S)=(S,S,S)^{\perp}= \{0\}.$ Therefore, By Theorem 4.7 of \cite{scha}, $S=\{0\}.$ It follows   that ${\frak J}= Rad({\frak J}),$ consequently ${\frak J}= \{0\}$ because ${\frak J}$ is nilpotent and  $({\frak J},{\frak J},{\frak J})= {\frak J}.$  $\Box$

 \begin{coro}
Let  $({\frak J},B)$ be  a pseudo-euclidean Jordan algebra. If ${\frak J}\neq\{0\},$ then  $Z({\frak J})\neq \{0\}.$
\end{coro}
 
Now, we shall recall the notion of the T$^*$-extension of the
Jordan algebras. This notion was introduced by M.Bordemann in
$\cite{Bord}$ in order to study algebras endowed with associative
non-degenerate symmetric bilinear forms. Let
 $\frak{J}$ be a Jordan algebra,
 ${\frak J}^*$ be the dual space of $\frak{J}$ and
$\theta :\frak{J}\times \frak{J}\longrightarrow {\frak J}^* $
be a bilinear map. On the vector space
  ${\tilde {\frak J}}= {\frak J}\oplus  {\frak J}^*,$ we define the following product:
\begin{eqnarray}\label{eq-2}
(x+f)(y+h)=xy+h\circ R_x+f\circ R_y+\theta(x,y), \,\,  \forall x,y \in \frak{J}, f,g \in \frak{J}^*.
\end{eqnarray}
 
\begin{prop}
${\tilde {\frak J}}$ endowed with the product above, is a Jordan
algebra if and only if  $\theta$ is symmetric and satisfies the
following identity:
\begin{eqnarray}\label{id1}
\theta(xy,x^2)+ \theta(x,x)\circ R_{xy} +\theta(x,y)\circ R_{x^2}=\theta(x,yx^2)+\theta(y,x^2)\circ R_x
+ \theta(x,x)\circ R_y R_x,\, \forall x,y \in \frak{J}.&&
\end{eqnarray}
In this case, the bilinear form $B:\, \tilde{\frak{J}}\times\tilde{\frak{J}}\rightarrow \mathbb{K}$ defined by:
$$B(x+f,y+h)=f(y)+h(x),\,\,  \forall x,y \in {\frak J}, f,g \in {\frak J}^*, $$ is an associatif scalar product on $\tilde {\frak J}$ if and only if $\theta$ satisfies  \begin{eqnarray}\label{id2}
\theta(x,y)(z)=\theta(z,x)(y),\,\, \, \forall x,y,z\in \frak{J}.
\end{eqnarray}
 \end{prop}
\begin{defi}
 If $\frak{J}$ is a Jordan algebra and $\theta :{\frak J}\times {\frak J}\longrightarrow {\frak J}^* $ is a bilinear map which satisfies
  $(\ref{id1})$ and $(\ref{id2}),$ then the pseudo-euclidean Jordan algebra $\tilde {\frak J}$ is called the $T^*-extension$ of $\frak{J}$ by mean of $\theta$ and denoted $T^*_{\theta}{\frak J}.$
\end{defi}
\begin{rema} \label{etoile}
 If $\frak{J}$ is a nilpotent Jordan algebra, then
$T^*_0\frak{J}$ is a nilpotent pseudo-euclidean Jordan algebra.
But if $\frak{J}$ is a semi-simple Jordan algebra then
$T^*_0\frak{J}$ is a pseudo-euclidean Jordan algebra which is not
nilpotent nor semi-simple, because $\frak{J}^*$ is a nilpotent
ideal of $T^*_0\frak{J}$ and $(T^*_0\frak{J})^2=T^*_0\frak{J}.$
Which show that the class of pseudo-euclidean Jordan algebra is
very rich.
\end{rema}
Now, we will give an inductive description of pseudo-euclidean
Jordan algebras by using the notion of double extension. Let us
start by  introducing some  Jordan algebras extensions.
\section{Some extensions of Jordan algebras}\label{ext}
\subsection{Central extension of Jordan algebras }

Let $\frak{J}_{1}$ be a Jordan algebra, $\cal{V}$ be a vector space and $\varphi:\frak{J}_{1}\times \frak{J}_{1}\longrightarrow \cal{V}$ be a bilinear map. On the vector space $\frak{J}=\frak{J}_{1}\oplus \cal{V}$ we  define the following product:
\begin{eqnarray}\label{EC1}
(x+v)(y+w)=xy+\varphi(x,y), \, \,  \forall x,y \in \frak{J}_{1}, v,w \in \cal{V}.
\end{eqnarray}

It is easy to cheks that $\frak{J}$ endowed with the product
(\ref{EC1}) is a Jordan algebra if and only if $\varphi$ is
symmetric and satisfies:
$$\varphi(xy,x^2)=\varphi(x,yx^2) \, \, \forall x,y \in \frak{J}_{1}.$$
We say that $\frak{J}$ is the central extension of $\frak{J}_{1}$
by $\cal{V}$ (by means of $\varphi$). In this case, $\cal{V}$ is
contained in the annulator of $\frak{J}$.
\subsection{Representations and  semi-direct products of Jordan  algebras}
Let $\frak{J}_{1}$ be a Jordan algebra, $\cal{V}$ be a vector space and $\pi:\frak{J}_{1}\longrightarrow End(\cal{V})$ be a linear map. On $\frak{J}=\frak{J}_{1}\oplus \cal{V}$ we define the following product:
\begin{eqnarray}\label{SD1}
(x+v)(y+w)=xy+\pi(x)w+\pi(y)v, \,  \, \, \forall x,y \in \frak{J}_{1}, v,w \in \cal{V}.
\end{eqnarray}
\begin{prop}\label{rep}
 $\frak{J}$ endowed with the product (\ref{SD1}) is a Jordan algebra if and only if $\pi$ satisfies the following conditions:\\
\begin{equation} \label{eq0}
\renewcommand{\arraystretch}{1.1}
\begin{array}{rll}
&&(i)\pi(x^2)\pi(x)-\pi(x)\pi(x^2)=0\\
&&(ii)2\pi(xy)\pi(x)+\pi(x^2)\pi(y)-2\pi(x)\pi(y)\pi(x)-\pi(x^2y)=0,\, \forall x,y \in \frak{J}_{1}.
\end{array}
\end{equation}
In this case $\pi$ is called a representation of the Jordan algebra $\frak{J}_1$ in $\cal{V}$.
\end{prop}
\begin{rema}
 Recall that ${T^*}_0{\frak J}= {\frak J}\oplus {{\frak J}^*} $
endowed with the product
$$(x+f)(y+h)=xy+h\circ R_x+f\circ R_y,\,\, \forall x,y \in {\frak J}
, f,g \in {\frak J}^*,$$ is a Jordan algebra. Thus $\rho:{\frak J}
\longrightarrow End({\frak J}^*) $ defined by $\rho(x)f = f\circ R_x,\,
\forall x \in {\frak J},$ is a representation of $\frak{J}$.
 This representation is called the co-adjoint representation of $\frak{J}.$
\end{rema}
\begin{coro} \label{cor2}
 $\pi:\frak{J}_{1}\longrightarrow End(\cal{V})$ is a representation of $\frak{J}_{1}$ if and only if
\begin{eqnarray*}
&\pi(xy)\pi(z)+\pi(yz)\pi(x)+\pi(xz)\pi(y)&\\
&=\pi(y)\pi(xz)+\pi(x)\pi(yz)+\pi(z)\pi(xy)&\forall x,y,z \in \frak{J}_1.\\
&=\pi((xy)z)+\pi(x)\pi(z)\pi(y)+\pi(y)\pi(z)\pi(x),&
\end{eqnarray*}
\end{coro}
\dem We will proceed by linearization to show this Corollary. Let
$x,y,z \in {\frak{J}_1}$ and $\lambda \in \mathbb{K}.$ Replace $x$
by $x+\lambda y$ in $(i)$ of (\ref{eq0}). The term of $\lambda$ is
$0$. Thus
\begin{eqnarray}\label{eq2}
2\pi(xy)\pi(x)+\pi(x^2)\pi(y)-\pi(y)\pi(x^2)-2\pi(x)\pi(xy)=0
\end{eqnarray}
Replace $x$ in (\ref{eq2}) by $x+\lambda z.$ The fact that the term of $\lambda$ is $0$ implies that
\begin{equation} \label{eq3}
\renewcommand{\arraystretch}{1.1}
\begin{array}{rll}
&\pi(xy)\pi(z)+\pi(yz)\pi(x)+\pi(xz)\pi(y)&\\
&=\pi(y)\pi(xz)+\pi(x)\pi(yz)+\pi(z)\pi(xy).&
\end{array}
\end{equation}
Conversely, if we replace  $y $ and $z$  in (\ref{eq3}) by $x$ we obtain $(i)$ of  (\ref{eq0}).\\
By the same argument, we show that  $(ii)$ of (\ref{eq0}) is equivalent to
\begin{equation} \label{eq4}
\renewcommand{\arraystretch}{1.1}
\begin{array}{rll}
&\pi(xy)\pi(z)+\pi(yz)\pi(x)+\pi(xz)\pi(y)&\\
&=\pi((xy)z)+\pi(x)\pi(z)\pi(y)+\pi(y)\pi(z)\pi(x).\Box&
\end{array}
\end{equation}

\begin{coro}\label{lie}
Let $\frak{J}$ be a Jordan algebra, $\cal{V}$ be a vector space
and $\pi:\frak{J} \longrightarrow End(\cal{V})$ be a
representation of $\frak{J} $.  Then,
$$\pi(x,y,z)=[\pi(y),[\pi(x),\pi(z)] ,\,\, \forall x,y,z\in \frak{J}.$$
\end{coro}
\dem Let $x,y,z\in \frak{J}.$ By Proposition \ref{cor2}, we have
the following equalities:
\begin{eqnarray*}
&\pi((xy)z)=\pi(xy)\pi(z)+\pi(yz)\pi(x)+\pi(xz)\pi(y)
-\pi(x)\pi(z)\pi(y)-\pi(y)\pi(z)\pi(x)&\\
 &\pi(x(yz))=\pi((zy)x)=\pi(zy)\pi(x)+\pi(xy)\pi(z)+\pi(xz)\pi(y)
-\pi(z)\pi(x)\pi(y)-\pi(y)\pi(x)\pi(z).&
\end{eqnarray*}
Consequently, $\pi(x,y,z)=[\pi(y),[\pi(x),\pi(z)].$  $\Box$

%%%%%%%%%%%%%%%%%%%%%%%%%%%%%%%%%
%%%%%%%%%%%%%%%%%%%%%%%%%%%%%%%%%%%
 \begin{rema} \label{srep}
\begin{enumerate}
\item [(i)] Let  $\frak{J}$ be a Jordan algebra and $\cal{V},W$ be two vector spaces. If
$\pi :\frak{J} \longrightarrow End(\cal{V})$ and
 $\rho :\frak{J}\longrightarrow End(W)$ are two representations of Jordan algebras, then
 $\tilde{\pi}:=\pi\oplus \rho :\frak{J} \longrightarrow
 End(\cal{V}\oplus W)$ defined by
$$(\pi\oplus \rho)(x)(v+w)=\pi(x)v+\rho(x)w,\, \forall x \in {\frak J}, v\in {\cal V}, w\in W,$$
is a representation of $\frak{J}$ called the direct sum of the representations $\pi$ and $\rho .$
\item [(ii)] Let $\frak{J}$ be a Jordan algebra. Replace $z$ by $x$ in the identity (\ref{eq-5}), we obtain:
$$ 2R_{xy} R_{x}+R_{x^2} R_{y}-2R_{x} R_{y} R_{x}-R_{x^2y}=0,\,\forall x,y \in \frak{J}.$$

Further, The fact that $\frak{J}$ is a Jordan algebra implies that
$R_xR_{x^2}=R_{x^2}R_x, \forall x \in \frak{J}.$ Which proves that
the linear map $R:\frak{J} \longrightarrow End(\frak{J})$ defined
by $ R(x)= R_x,\,\,\forall x \in \frak{J}$ is a representation of
$\frak{J}$. This representation is called the  adjoint or the
regular representation of $\frak{J}.$
\end{enumerate}
\end{rema}

\begin{prop} \label{pr3}
Let  $\frak{J}_{1},\frak{J}_{2}$ be two Jordan algebras and
$\pi:\frak{J}_{1}\longrightarrow End(\frak{J}_{2})$ be a linear
map. We define on the vector space $\frak{J}=\frak{J}_{1}\oplus
\frak{J}_{2}$ the following product:
\begin{eqnarray*}
(x+y)\star(x'+y')=xx'+\pi(x)y'+\pi(x')y+yy',\,\forall x,x' \in \frak{J}_{1},y,y' \in \frak{J}_{2}.
\end{eqnarray*}
Then, $(\frak{J},\star)$ is a Jordan algebra if and only if $\pi$ satisfies the following conditions:
 \begin{eqnarray*}
 &(i)& \pi(x^2)\pi(x)y'+\pi(x^2)(yy')+2\Bigl(\pi(x)y'\Bigl)\Bigl(\pi(x)y\Bigl)+\Bigl(\pi(x)y'\Bigl)y^2+
2\Bigl(yy'\Bigl)\Bigl(\pi(x)y\Bigl) \\
&&=\pi(x)\pi(x^2)y'+2\pi(x)\Bigl(y'\Bigl(\pi(x)y\Bigl)\Bigl)+\pi(x)(y'y^2)+\Bigl(\pi(x^2)y'\Bigl)y+
2\Bigl(y'\Bigl(\pi(x)y\Bigl)\Bigl)y,
\end{eqnarray*}
 \begin{eqnarray*}
&(ii)&2\pi(xx')\pi(x)y+\pi(xx')y^2+\pi(x^2)\pi(x')y+\Bigl(\pi(x')y\Bigl)y^2+2\Bigl(\pi(x')y\Bigl)\Bigl(\pi(x)y\Bigl)\\
&&=2\pi(x)\pi(x')\pi(x)y+\pi(x)\pi(x')y^2+\pi(x'x^2)y+2\Bigl(\pi(x')\pi(x)y\Bigl)y+\Bigl(\pi(x')y^2\Bigl)y,
\end{eqnarray*}
$$\forall x,x' \in \frak{J}_{1},y,y' \in \frak{J}_{2}.$$
\end{prop}
 \dem Suppose that $(\frak{J},\star)$ is a Jordan algebra.
Let $x,x' \in \frak{J}_{1}$ and $y,y'\in \frak{J}_{2}.$
\begin{eqnarray*}
 0&=&\Bigl((x+y)\star y'\Bigl)\star(x+y)^2-(x+y)\star\Bigl(y'\star(x+y)^2\Bigl)\\
&=&\pi(x^2)\pi(x)y'+\pi(x^2)(yy')+2\Bigl(\pi(x)y'\Bigl)\Bigl(\pi(x)y\Bigl)
+\Bigl(\pi(x)y'\Bigl)y^2+2\Bigl(yy'\Bigl)\Bigl(\pi(x)y\Bigl) \\
&=&\pi(x)\pi(x^2)y'+2\pi(x)\Bigl(y'\Bigl(\pi(x)y\Bigl)\Bigl)+\pi(x)(y'y^2)+\Bigl(\pi(x^2)y'\Bigl)y+
2\Bigl(y'\Bigl(\pi(x)y\Bigl)\Bigl)y,
\end{eqnarray*}
\begin{eqnarray*}
0&=& \Bigl((x+y)\star x'\Bigl)\star(x+y)^2-(x+y)\star\Bigl(x'\star(x+y)^2\Bigl)\\
&=&2\pi(xx')\pi(x)y+\pi(xx')y^2+\pi(x^2)\pi(x')y+\Bigl(\pi(x')y\Bigl)y^2
+2\Bigl(\pi(x')y\Bigl)\Bigl(\pi(x)y\Bigl)\\
 &=&2\pi(x)\pi(x')\pi(x)y+\pi(x)\pi(x')y^2+\pi(x'x^2)y+
2\Bigl(\pi(x')\pi(x)y\Bigl)y+\Bigl(\pi(x')y^2\Bigl)y.
\end{eqnarray*}
Which give the identities $(i)$ and $(ii).$

Conversely, suppose that $(i)$ and $(ii) $ are satisfied,
then for all  $x,x' \in \frak{J}_{1}$ and $y,y'\in \frak{J}_{2},$
we have:

$$  \Bigl((x+y)\star y'\Bigl)\star(x+y)^2-(x+y)\star\Bigl(y'\star(x+y)^2\Bigl)= 0,$$
$$ \Bigl((x+y)\star x'\Bigl)\star(x+y)^2-(x+y)\star\Bigl(x'\star(x+y)^2\Bigl)= 0. $$

Consequently,
$$  \Bigl((x+y)\star(x'+ y')\Bigl)\star(x+y)^2-(x+y)\star\Bigl((x'+y')\star(x+y)^2\Bigl)= 0, \forall x,x' \in \frak{J}_{1}, y,y'\in \frak{J}_{2}.$$
 Since the product $\star$ is commutatif, then $\frak{J}$ endowed with $\star$ is a Jordan algebra.
  $\Box$

\begin{coro}\label{pr4}
Let $\frak{J}_{1},\frak{J}_{2}$ be two Jordan algebras and
$\pi:\frak{J}_{1}\longrightarrow End(\frak{J}_{2})$ be a
representation of Jordan algebra. We define on the vector space
$\frak{J}=\frak{J}_{1}\oplus \frak{J}_{2}$ the following product:
\begin{eqnarray*}
(x+y)\star(x'+y')=xx'+\pi(x)y'+\pi(x')y+yy',\,\forall x,x' \in \frak{J}_{1},y,y' \in \frak{J}_{2}.
\end{eqnarray*}
Then, $(\frak{J},\star)$ is a Jordan algebra if and only if $\pi$ satisfies the following conditions:\\

\begin{enumerate}
\item[(1)]$\,\,\,\pi(x^2)(yy')+2\Bigl(\pi(x)y'\Bigl)\Bigl(\pi(x)y\Bigl)
+\Bigl(\pi(x)y'\Bigl)y^2+
2\Bigl(yy'\Bigl)\Bigl(\pi(x)y\Bigl) $\\

$=2\pi(x)\Bigl(y'\Bigl(\pi(x)y\Bigl)\Bigl)+\pi(x)(y'y^2)
+\Bigl(\pi(x^2)y'\Bigl)y+
2\Bigl(y'\Bigl(\pi(x)y\Bigl)\Bigl)y $\\

\item[(2)]$\,\,\,\Bigl(\pi(x)y\Bigl)y^2=\Bigl(\pi(x)y^2\Bigl)y$\\

\item[(3)]$\,\,\,\pi(xx')y^2+
2\Bigl(\pi(x')y\Bigl)\Bigl(\pi(x)y\Bigl)=\pi(x)\pi(x')y^2
+2\Bigl(\pi(x')\pi(x)y\Bigl)y,$
\end{enumerate}

\end{coro}
\dem Let $x,x'\in \frak{J}_{1},\,\, y,y' \in \frak{J}_{2}$. First,
since $\pi$ is a representation, then
 $\pi(x^2)\pi(x)y'=\pi(x)\pi(x^2)y'.$ consequently,  the identity $(i)$ of the Proposition \ref{pr3}
 is equivalent to the identity $(1).$ Second, we obtain the identity $(2)$  by replacing $x$ by $0$ in $(ii)$ of the
 Proposition $\ref{pr3}.$
Finaly, since $\pi$ is a representation, then
 \begin{eqnarray}\label{sup}
&&2\pi(xx')\pi(x)y+\pi(x^2)\pi(x')y=2\pi(x)\pi(x')\pi(x)y+\pi(x^2x')y.
\end{eqnarray}Hence, the identity $(2)$ implies that $(3)$ is satisfies.
Finally, remark that the sum of the identities $(2)$, $(3)$ and (\ref{sup}), give
the identity $(ii)$  of the Proposition \ref{pr3}. $\Box$
\begin{defi}
Let  $\frak{J}_{1},\frak{J}_{2}$ be two Jordan algebras and
$\pi:\frak{J}_{1}\longrightarrow End(\frak{J}_{2})$ be a
representation of Jordan algebras. $\pi$ is called admissible, if
$\pi$ satisfies the identities $(1), (2)$ and $(3)$ of
 Corollary \ref{pr4}. In this case, the Jordan algebra $(\frak{J}=\frak{J}_{1}\oplus \frak{J}_{2},\star)$ where the product
 $\star$ is defined by :
\begin{eqnarray*}
(x+y)\star(x'+y')=xx'+\pi(x)y'+\pi(x')y+yy',\,\forall x,x' \in \frak{J}_{1},y,y' \in \frak{J}_{2},
\end{eqnarray*}
is called the semi-direct product of $\frak{J}_{2}$ by $\frak{J}_1$ by means of $\pi.$
\end{defi}
\begin{coro}\label{adj}
 The adjoint representation of a Jordan algebra  $\frak{J}$   is an admissible representation of $\frak{J}$.
\end{coro}
\dem The identities $(1)$ and $(3)$ of Corollary $\ref{pr4}$ are a
consequences of the equality $(\ref{eqq}).$ Further, the fact that
$\frak J$ is a Jordan algebra implies that $(2)$ is satisfies.
$\Box$
 \subsection{Generalized semi-direct product}

Let $\frak{J}$ be a Jordan algebra, $(D,x_{0}) \in
End(\frak{J})\times \frak{J},$ and $\mathbb{K}a$ be
one-dimensional algebra with zero product. On the vector space
$\tilde{\frak{J}}=\mathbb{K}a\oplus \frak{J}$, we define the
following product:
$$\begin{array}{l}
x\star y:= xy,
x\star a= a\star x= D(x), \hskip 0.5cm
a\star a= x_{0},\,\,\,\forall x,y \in \frak{J}.
 \end{array}$$
 $\tilde{\frak{J}}$ endowed with the product above is a Jordan algebra, if and only if, for all
 $x,y$ in $\frak{J}$, the pair $(D,x_{0})$ satisfies the following conditions
 $$\begin{array}{lll} (C_1) D(x^2y)=x^2D(y)+2D(x)(xy)-2x(D(x)y),
&(C_5)D^2(x^2)=2(D(x))^2-2xD^2(x)+x_{0}x^2,
\\
(C_2):D(x)D(y)-D(D(x)y)=\frac{1}{2}(x_{0},y,x),
&(C_6):D^3(x)=\frac{3}{2}x_{0}D(x) -\frac{1}{2}xD(x_{0}) ,
\\
(C_3): D(x_{0}x)=x_{0}D(x),\hskip0.5cm(C_4):xD(x^2)=x^2D(x),
&(C_7):D^2(x_{0})=x_{0}^{2}.
\end{array}$$
In this case, $(D,x_0)$ is called an admissible pair of
$\frak{J}$ and the Jordan algebra  $\tilde{\frak{J}},$ is called the generalized
 semi-direct product of $\frak{J}$ by the one-dimensional algebra with zero product by means of the pair
$(D,x_{0})$. It is easy to see that $\frak{J}$ is an ideal of $\tilde{\frak{J}}$ and   $\mathbb{K}a$ is not
 in general a subalgebra of $\tilde{\frak{J}}.$

\subsection{Double extension of Jordan algebras}
\begin{defi}
Let $\frak{J}$ be a pseudo-euclidean Jordan algebra and $B$ be an associatif scalar product on $\frak{J}$.
 An endomorphism $f$ of $\frak{J}$ is called $B-$symmetric, (resp. $B-$antisymmetric) if
  $B(f(x),y)=B(x,f(y)),$   (resp. $B(f(x),y)=-B(x,f(y))), \forall x,y\in
  \frak{J}.$ Denoted by $End_s({\frak J})$ (resp. $End_a({\frak
  J})$) the subspace of $B-$symmetric (resp. $B-$antisymmetric)
  endomorphism of $\frak J$.
\end{defi}

Let $\frak{J}_{1}$ be a pseudo-euclidean Jordan algebra, $B_{1}$
be an associatif scalar product on $\frak{J}_{1}$, $\frak{J}_{2}$
be a Jordan algebra not necessarily pseudo-euclidean and ${\pi}:
\frak{J}_{2} \longrightarrow End_s(\frak{J}_1)$ be an admissible
representation of $\frak{J}_{2}$ in $(\frak{J}_1).$ Now,
 Let $\varphi: \frak{J}_{1}\times \frak{J}_{1} \longrightarrow End\Biggl((\frak{J}_{2})^*\Biggl)$ be 
  the bilinear map defined by:
  $$\varphi(x,y)(z)=B_{1}(\pi(z)x,y), \,\forall x,y\in \frak{J}_{1},\, z\in \frak{J}_{2}.$$

Since, by Corollary $\ref{pr4},$ for all $(x,z)\in \frak{J}_{1}\times \frak{J}_{2},$ $\pi(z)$
is $B_1-$symmetric and satisfies $(\pi(z)x)x^2=(\pi(z)x^2)x,$
 then $\varphi$ is a symmetric bilinear
mapping which satisfies: $\varphi(xy,x^2)=\varphi(x,yx^2), \,
\forall x,y\in \frak{J}_{1}.$ Consequently, we can consider the
central extension $\frak{J}=\frak{J}_{1}\oplus \frak{J}_{2}^{*}$
of $\frak{J}_{1}$ by $\frak{J}_{2}^{*}$ by means of $\varphi.$
Recall that the product on  $\frak{J}$ is defined by:
\begin{eqnarray*}
(x+f)(x'+f')=xx'+\varphi(x,x'), \, \forall x,x'\in \frak{J}_{1},\, f,f'\in \frak{J}_{2}^{*}
\end{eqnarray*}
By the Remark \ref{srep}, the linear map $\tilde{\pi}: \frak{J}_{2} \longrightarrow End(\frak{J})$ defined by
\begin{eqnarray*}
\tilde{\pi}(y)(x+f)=\pi(y)x + \rho(y)f=\pi(y)x + f\circ R_y,\,\, \forall y \in \frak{J}_{2}, x \in \frak{J}_{1}, f\in \frak{J}_{2}^{*},
\end{eqnarray*}
where $\rho$ is the coadjoint representation of $\frak{J}_{2},$  is a
representation of $\frak{J}_2$ in $\frak{J}.$ Now, we shall prove
that $\tilde{\pi}$ is admissible.

Since $\pi$ is an admissible representation, then $\tilde \pi$
satisfied the condition $(1)$ of Corollary $\ref{pr4}$ if and only
if $\Omega(x,x',y)= 0, \forall x, x' \in {\frak{J}_{1}}, \forall y \in {\frak{J}_{2}}, $ where
\begin{eqnarray*}
&\Omega(x,x',y):=\rho(y^2)\varphi(x,x')+2\varphi(\pi(y)x',\pi(y)x)+
\varphi(\pi(y)x',x^2)+2\varphi(xx',\pi(y)x)&\\&-2\rho(y)\varphi(x',\pi(y)x)
-\rho(y)\varphi(x',x^2)-\varphi(\pi(y^2)x',x)-2\varphi(x'\pi(y)x,x).\hskip-2cm&
\end{eqnarray*}

It is clear that for all $x, x' \in {\frak{J}_{1}},   y, y' \in {\frak{J}_{2}},$ we have $\Omega(x,x',y)(y')= B(x',\Gamma(x,y,y')), $  where
\begin{eqnarray*}
& \Gamma(x,y,y')=\Bigl(\pi(y'y^2)+
2\pi(y)\pi(y')\pi(y)-2\pi(yy')\pi(y)-\pi(y^2)\pi(y')\Bigl)x&\\&-\biggl(\pi(yy')x^2+
2\Bigl(\pi(y)x\Bigl)\Bigl(\pi(y')x\Bigl)-\pi(y)\pi(y')x^2-
2\Bigl(\pi(y')\pi(y)x\Bigl)x\biggl)=0.&
\end{eqnarray*}
By $(ii)$ of $(\ref{eq0})$ and the identity $(3)$ of
Corollary $\ref{pr4}$ we have $\Gamma(x,y,y')= 0.$ Consequently
$\Omega(x,x',y)=0.$ Hence $\tilde{\pi}$ satisfies the condition
$(1)$ of corollary $\ref{pr4}$. \\
we need the following lemma to show that $\tilde{\pi} $ satisfies the conditions $(2)$ et $(3)$ of Corollary \ref{pr4}.
\begin{lema}\label{res}
Let $x\in \frak{J}_1,\, y,z,u\in \frak{J}_2$, then:
\begin{enumerate}
\item[(i)] $B_1([\pi(y),\pi(z)]x^2 ,x)=0;$
 \item[(ii)]$B_1(\Bigl(\pi\bigl((y,z,u)\bigl)+2[\pi(z),\pi(u)\pi(y)]\Bigl)x,x)=0.$
\end{enumerate}
\end{lema}
\dem  Let $x\in \frak{J}_1,\, y,z,u\in \frak{J}_2.$
%%%%%%%%%%%%%%%%%%%%%%%%%%%%%%%%%%%%%%%%%%%%%%%%
$(i)$ The identity $B_1([\pi(y),\pi(z)]x^2 ,x)=0$ is equivalent to 
\begin{eqnarray*}
2B_1([\pi(y),\pi(z)]x,xx')+B_1([\pi(y),\pi(z)]x',x^2)=0.
\end{eqnarray*}

In fact we proceed by linearization, we replace $x$ by $x+\lambda
x'$ in the identity $B_1([\pi(y),\pi(z)]x^2 ,x)=0,$ where $\lambda \in
\mathbb{K}.$
%%%%%%%%%%%%%%%%%%%%%%%%%%%%%%%%%%%%%%%%%%%%%%%%%%%%%%
 
Since $[\pi(y),\pi(z)]$ is a $B_1-$antisymmetric endomorphism and
$B_1$ is associative, then
 \begin{eqnarray*}
2B_1([\pi(y),\pi(z)]x,xx')+B_1([\pi(y),\pi(z)]x',x^2)
 =B_1 (2\Bigl([\pi(y),\pi(z)]x\Bigl)x-[\pi(y),\pi(z)]x^2,x')=0. &&
\end{eqnarray*}
 Which equivalent to
\begin{eqnarray}\label{eq10}
 2\Bigl([\pi(y),\pi(z)]x\Bigl)x-[\pi(y),\pi(z)]x^2= 0.
\end{eqnarray}
 By $(3)$ of  Corollary $\ref{pr4}$ we have:
\begin{eqnarray} \label{eq11}
 \pi(yz)x^2+2\Bigl(\pi(z)x\Bigl)\Bigl(\pi(y)x\Bigl)= \pi(y)\pi(z)x^2+2\Bigl(\pi(z)\pi(y)x\Bigl)x
\end{eqnarray}
and
\begin{eqnarray} \label{eq12}
\pi(zy)x^2+2\Bigl(\pi(y)x\Bigl)\Bigl(\pi(z)x\Bigl)=
\pi(z)\pi(y)x^2+2\Bigl(\pi(y)\pi(z)x\Bigl)x.
 \end{eqnarray}
 ${\frak J}_1$ and ${\frak J}_2$ are commutative, then
\begin{eqnarray*}
\pi(y)\pi(z)x^2+2(\pi(z)\pi(y)x)x=\pi(z)\pi(y)x^2+2(\pi(y)\pi(z)x)x.
\end{eqnarray*}
Which give the identity (\ref{eq10}) or the identity $(i)$ of
Lemma.\\
 $(ii)$ One pose
$F(y,z,u)=\pi\bigl((y,z,u)\bigl)+2[\pi(z),\pi(u)\pi(y)],\,\,\,\forall
y,z,u\in\frak{J}_1 $. Then $F(y,z,u)$ is a $B_1-$antisymmetric
endomorphism. In fact,
\begin{eqnarray*}
&&B_1(F(y,z,u)x,x')+B_1(x,F(y,z,u)x')\\&&=
B_1\Bigl(x,\pi\bigl((y,z,u)\bigl)(x')\Bigl)+2B_1\Bigl(x,[\pi(y)\pi(u),\pi(z)]x'\Bigl)
+B_1\Bigl(x,F(y,z,u)x'\Bigl)\\&&= 2B_1\Bigl(x,
(\pi\bigl((y,z,u)\bigl)-[\pi(z),[\pi(y),\pi(u)]])(x')\Bigl)=0,\,\,\,\,\,(by Corollary
\ref{lie}).
\end{eqnarray*}
Which prove that, $F(y,z,u)$ is a $B_1-$antisymmetric
endomorphism. Consequently, $B_1(F(y,z,u)x,x)=-B_1(x,F(y,z,u)x)$.
The fact that $B_1$ is symmetric implies that
$B_1(F(y,z,u)x,x)=0.$ $\Box$

{\bf 2.} Since $\pi$ satisfies the equality $(2)$ of Corollary $\ref{lie}$ and $\varphi$ is symmetric, then $\tilde{\pi}$
satisfies the identity $(2)$ of Corollary \ref{pr4} if and only if
$$\varphi(\pi(y)x,x^2)-\varphi(\pi(y)x^2,x)=0,\,\,\,\forall x\in {\frak J}_1,\,\,\,y\in{\frak J}_2.$$
By $(i)$ of Lemma \ref{res}, we have
$\Bigl(\varphi(\pi(y)x,x^2)-\varphi(\pi(y)x^2,x)
\Bigl)(z)=B_{1}([\pi(y),\pi(z)]x^2,x)=0.$ Hence,
$$\Bigl(\tilde\pi(y)(x+f)\Bigl)(x+f)^2=
\Bigl(\tilde\pi(y)(x+f)^2\Bigl)(x+f).$$

{\bf3.} Since $\pi$ is admissible, then $\tilde{\pi}$ satisfies the identity $(3)$ of  Corollairy \ref{pr4} if and
only if
$$\rho(yy')\varphi(x,x)+2\varphi(\pi(y')x,\pi(y)x) -\rho(y)\rho(y')\varphi(x,x)
-2\varphi(\pi(y')\pi(y)x,x)=0.$$ By $(i)$ and $(ii)$ of
lemma \ref{res}, we have:
\begin{eqnarray*}&&\Biggl(\rho(yy')\varphi(x,x)+2\varphi\Bigl(\pi(y')x,\pi(y)x\Bigl)
-\rho(y)\rho(y')\varphi(x,x)
-2\varphi\Bigl(\pi(y')\pi(y)x,x\Bigl)\Biggl)(z)\\ \\
&&=B_1\Bigl((\pi((y',y,z)+2[\pi(y),\pi(z)\pi(y')])x,x\Bigl)
-B_1\Bigl([\pi(y'),\pi(z)]x^2,x\Bigl)=0.
\end{eqnarray*}
Consequently, $\tilde{\pi}$ satisfies the identity $(3)$ of
Corollary \ref{pr4}.
We conclude that $\tilde \pi$ is an admissible representation of $\frak{J}_2$ in $\frak{J}_1\oplus{\frak{J}_2}^*.$\\
%%%%%%%%%%%%%%%%%%%%%%%%%%%%%%%%%%%%%%%%%%%%%%%%%%%%%
%%%%%%%%%%%%%%%%%%%%%%%%%%%%%%%%%%%%%%%%%%%%%%%%%%%%%%%%%
Hence, we can consider,
 the
semi-direct product, $\tilde{{\frak{J}}}={\frak{J}}_2\oplus
\frak{J}$ of $\frak{J}$ by ${\frak{J}}_2$ by means of
$\tilde{\pi}$. The product in  $\tilde{{\frak{J}}}$ is given by:
$$
(x+y+f)(x'+y'+f')=xx'+yy'+\pi(x)y'+\pi(x')y+\rho(x)f'+\rho(x')f+\varphi(y,y'),$$
$\forall x,x'\in {\frak J}_2,y,y'\in {\frak J}_1,f,f'\in
{\frak{J}}_2^{*}.$
 Further, let $\gamma$ be an associative symmetric bilinear form  not necessarily nondegenerate
on ${\frak{J}}_{2}\times {\frak{J}}_{2}$. Then, the bilinear form
$B$ defined on $\tilde{\frak{J}}\times \tilde{\frak{J}}$ by:
\begin{eqnarray*}
B:({\frak{J}}_{2}\oplus {\frak{J}}_{1} \oplus {\frak{J}}_{2}^{*})\times ({\frak{J}}_{2}\oplus {\frak{J}}_{1} \oplus {\frak{J}}_{2}^{*})&
\longrightarrow& \mathbb K\\
(x+y+f,x'+y'+f')\,\,\,\,\,\, &\longmapsto &\gamma(x,x')+B_{1}(y,y')+f(x')+f'(x)
\end{eqnarray*}
 is an associatif scalar product on $\tilde{\frak J}$.\\

Then we have proved the following Theorem.

\begin{theor} \label{th1}
Let $({\frak{J}}_{1},B_{1})$ be a pseudo-euclidean Jordan algebra,
${\frak{J}}_{2}$ be a Jordan algebra and
$\pi:{\frak{J}}_{2}\longrightarrow \mbox{End}_s({\frak{J}}_{1})$
be an admissible representation of Jordan algebras. Let us
consider the symmetric bilinear map $\varphi:{\frak{J}}_{1}\times
{\frak{J}}_{1}\longrightarrow {\frak{J}}_{2}^{*}$
 defined by: $\varphi(y,y')(x)= B_{1}(\pi(x)y,y'),\,\,
\forall x \in {\frak{J}}_{2},y,y'\in {\frak{J}}_{1}.$ Then, the vector space $\tilde{\frak{J}}={\frak{J}}_{2}\oplus {\frak{J}}_{1}
\oplus {\frak{J}}_{2}^{*}$ endowed with the product
\begin{eqnarray*}
&&(x+y+f)(x'+y'+f')=xx'+yy'+\pi(x)y'+\pi(x')y+f'\circ R_{x}+f\circ R_{x'}+\varphi(y,y'),\\
&&\forall x,x'\in{\frak{J}} _2,y,y'\in{\frak{J}} _1,f,f'\in {\frak{J}}_2^{*},
\end{eqnarray*}
is a Jordan algebra. Moreover, if  $\gamma$ is an associative bilinear form on ${\frak{J}}_{2}\times {\frak{J}}_{2},$
then the bilinear form $B_\gamma$ defined on
$\tilde{\frak{J}}\times \tilde{\frak{J}}$ by:
\begin{eqnarray*}
B_\gamma:({\frak{J}}_{2}\oplus {\frak{J}}_{1} \oplus {\frak{J}}_{2}^{*})\times ({\frak{J}}_{2}\oplus {\frak{J}}_{1} \oplus {\frak{J}}_{2}^{*})&
\longrightarrow& \mathbb K\\
(x+y+f,x'+y'+f')\,\,\,\,\,\, &\longmapsto &\gamma(x,x')+B_{1}(y,y')+f(x')+f'(x)
\end{eqnarray*}
is an associatif scalar product on $\tilde{\frak{J}}.$

The Jordan algebra $(\tilde {\frak J},B_0)\,$ (or ${\frak J}$)  is called the double extension of $({\frak{J}}_1,B_1)$ by ${\frak{J}}_2$ by means of $\pi.$
\end{theor}
\begin{rema}\label{rqq1}
If ${\frak{J}}_{2}$ admits an  associative bilinear form
$\gamma_1\not= 0$, then $\tilde{\frak{J}}$ has at least two linearly
independent associatif scalar products.
\end{rema}
\subsection{Generalized double extension of pseudo-euclidean Jordan algebras by the one-dimensional algebra with zero product}
Let $(\frak{J}_1,B_1)$ be a pseudo-euclidean Jordan algebra,
$\mathbb{K}b$ be the one-dimensional algebra with zero product and
$(D,x_0)\in End_s({\frak{J}}_1,B_1)\times \frak{J}_1$ be an
admissible pair. Let $\varphi:\frak{J}_1\times
\frak{J}_1\longrightarrow \mathbb K$ be the symmetric bilinear
form defined by:
\begin{eqnarray*}
\varphi(x,y)=B_1(Dx,y), \,\, \forall x,y\in \frak{J}_1.
\end{eqnarray*}
Then, by the identity $(C_4),$
$\varphi(x^2,xy)=\varphi(x^2y,x),\,\,\,\forall x,y\in\frak{J}_1$.
Hence the vector space $\frak{J}=\frak{J}_1\oplus \mathbb{K} b $
endowed with the product: $ (x+\lambda
b)(y+\lambda'b):=xy+\varphi(x,y)b, \,\forall x,y \in
\frak{J}_1,\lambda,\lambda'\in \mathbb K $ is the central
extension of $\frak{J}_1$ by $\mathbb{K}b$ by means of $\varphi$.

Let $(\tilde{D},w_0)\in End(\frak{J})\times\frak{J}$ defined by:
\begin{eqnarray*}
&&\tilde{D}(\alpha b+x):=D(x)+B_1(x_{0},x)b,\, \,\,\forall x\in
\frak{J}_1,\alpha\in \mathbb{K},\,\,\,w_{0}=x_{0}+kb,\mbox{ o\`u
}k\in \mathbb{K}.
\end{eqnarray*}
Then the pair $(\tilde{D},w_{0})$ is admissible.
Hence, we can consider
 the generalized semi-direct product $\tilde{{\frak{J}}}=\mathbb K a\oplus \frak{J}$ of
$\frak{J}$ by the one-dimensional algebra with zero product $\mathbb
K a$ by means of $(\tilde{D},w_{0})$.

Moreover, the symmetric bilinear form $B:\tilde{\frak J}\times\tilde{\frak J}\longrightarrow\mathbb K$ defined by:
\begin{eqnarray*}
B_{\vert_{{{\frak J}_1\times{\frak J}_1}}}=B_1,\,\, B_1(a,b)=1,\,\,B_1(a,{\frak J}_1)=B_1(b,{\frak J}_1)=\{0\}\mbox{ et }B_1(a,a)=B_1(b,b)=0
\end{eqnarray*}
is an associatif scalar product on $\tilde{\frak J}$. Thus, we
have the following Theorem:
 \begin{theor}
Let $(\frak{J}_1,B_1)$ be a pseudo-euclidean Jordan algebra,
$\mathbb{K}b$ be the one-dimensional algebra with zero product and
$(D,x_0)\in End_s({\frak{J}}_1,B_1)\times \frak{J}_1$ be an
admissible pair. Let $\varphi:\frak{J}_1\times
\frak{J}_1\longrightarrow \mathbb K$ be the symmetric bilinear
form defined by:
\begin{eqnarray*}
\varphi(x,y)=B_1(Dx,y), \,\, \forall x,y\in \frak{J}_1.
\end{eqnarray*}
Then, the vector space $\tilde\frak{J}=\mathbb
Ka\oplus\frak{J}_1\oplus\mathbb Kb$, (where $\mathbb Ka$ is a
one-dimensional vector space), endowed with the following product:
$$\begin{array}{lllll}
\tilde{\frak{J}}\star b= b\star \tilde{\frak{J}}=\{0\},&
a\star a= w_{0}=x_{0}+kb,&
x\star y= xy+B_1(D(x),y)b,
&a \star x= x \star a= D(x)+B_1(x_{0},x)b,
\end{array}$$
$\forall x,y \in \frak{J}_1$ and with the symmetric bilinear form
$B$ defined by:
\begin{eqnarray*}
B_{\vert_{{{\frak J}_1\times{\frak J}_1}}}=B_1,\,\, B_1(a,b)=1,\,\,B_1(a,{\frak J}_1)=B_1(b,{\frak J}_1)=\{0\}\mbox{ et }B_1(a,a)=B_1(b,b)=0
\end{eqnarray*}
is a pseudo-euclidean Jordan algebra.
\end{theor}
\begin{defi}
The pseudo-euclidean  Jordan algebra $(\tilde{\frak J},B)$ is
called the generalized double extension of the pseudo-euclidean
Jordan algebra $({\frak J_1},B_1)$ by the one-dimensional Jordan
algebra with zero product $\mathbb{K}b$ by means of $(D,x_0,k)\in
\mbox{End}_s({\frak J}_1,B_1)\times {\frak J}_1\times\mathbb K $
(or by means of $(D,x_0)\in \mbox{End}_s({\frak J}_1,B_1)\times
{\frak J}_1$).
\end{defi}
\section{ Inductive description of pseudo-euclidean Jordan algebras}\label{desc}
\begin{theor}
Let $(\frak J,B)$ be an irreducible pseudo-euclidean Jordan algebra.
 If $\frak J=\cal I\oplus \cal V ,$ where $\cal I$ is a maximal ideal of $\frak J$ and $\cal V$ is a subalgebra of
  $\frak J$, then $\frak J$ is the double extension of the pseudo-euclidean Jordan algebra $(\cal W=\cal I/\cal
  I^\perp,$ $\tilde{B}$) by $\cal V$ by means of the representation
$\pi:\cal V\longrightarrow \mbox{End$_s$}({\cal W},\tilde{B})$
defined by $\pi(v)(s(i)):=s(R_v(i))=s(vi),\forall v\in\cal V,
\mbox{i}\in \mathcal I,$ where $\cal I^\perp$ is the orthogonal of
$\cal I$, $s$ is the canonical surjection of
$\cal I$ onto $\cal I/\cal I^\perp$ and $\tilde{B}$ is defined by:
$\tilde B(s(x),s(y)):=B(x,y),\,\,\forall x,y\in \cal \cal I.$
\end{theor}
\dem Since $\cal I$ is a maximal ideal of $\frak J$, Then, $\cal
I^\perp$ is a minimal one. Further, $\frak J$ is irreducible, then
${\cal I}\cap {\cal I}^\perp\neq\{0\}.$ Consequently, $\cal
I^\perp\subset \cal I.$ Consider $\cal A= \cal I^\perp \oplus\cal
V$ and $\cal A^\perp$ the orthogonal of $\cal A$. It is easy to
check that $B_{\vert_{{\cal A} \times {\cal A}}}$ is
nondegenerate. Thus, $\cal I=\cal I^\perp \oplus \cal A^\perp$. It
follows that, ${\frak J}=\cal A^\perp\oplus \cal I^\perp
\oplus\cal  V $. Now, let $a,b\in\cal A^\perp$, then
$ab=\alpha(a,b)+\beta(a,b)
 \mbox{ where }
  \alpha(a,b)\in {\cal I}^\perp \mbox{ and }\beta(a,b) \in {\cal A}^\perp.$ It is clear that
$\cal A^\perp$ endowed with $\beta$ is a Jordan algebra. Moreover,
the bilinear form $Q:=B_{\vert_{{\cal A}^{\perp}
 \times {\cal A}^{\perp}}}$
 is an associatif scalar product on $\cal A^\perp.$ Which implies that
$({\cal A}^\perp,Q)$ is a pseudo-euclidean Jordan algebra. consider 
$\theta:=s_{\vert_{{\cal A}^{\perp}}} :\cal
A^{\perp}\longrightarrow \cal I/\cal I^\perp,$  $\theta$ is an
  isomorphism of Jordan algebras. Let $v,w\in \cal{V} $ and $x\in \cal
A^{\perp}$. Then, $vx=i+a\in \cal I=\cal I^{\perp}\oplus \cal
A^{\perp}$ where $ i\in \cal{I}^{\perp}$ and $ a\in \cal
A^{\perp}.$ Thus,
 $$B(i,w)=B(vx-a,w)=B(vx,w)=B(x,vw)=0,$$
Thus $i=0$. Which implies that  $\cal{V}\cal A^{\perp}\subset \cal
A^{\perp}$. Consequently the map $\pi:\cal{\cal V} \longrightarrow
\mbox{End}(\cal A^{\perp}) $ defined by: $\pi(v)a=R_v(a),\,\,
\forall v\in \cal{V},$ $a\in \cal A^{\perp}$ is well defined and
it is an admissible representation of $\cal{V}$ because $R$ is an admissible
representation of $\frak J$. Hence, we can consider the double
extension $\cal{V}\oplus \cal A^{\perp}\oplus \cal{V}^*$ of $\cal
A^{\perp}$ by $\cal{V}$ by means of $\pi$. Now, one consider
$\nu:\cal{I}^{\perp} \longrightarrow  \cal{V}^*$ (resp.
$\delta:\cal{V}\longrightarrow (\cal{I}^{\perp})^{*} $) defined by
$\nu(i):=B(i,.),\,\,\,\forall i\in\cal{I}^{\perp}$ (resp.
$\delta(v)=B(v,.),\,\,\,\forall v\in\cal{V}$). Since $B$ is
nondegenerate, then $\nu$ (resp. $\delta$) is injective. Hence 
$\mbox{dim} \cal{I}^{\perp} =\mbox{dim} \cal{V}^*$ and $\nu$ is an
 isomorphism of  vector spaces. By the Corollary $\ref{adj},$ $R:\frak
J\longrightarrow End(\frak J);$ $ x\longmapsto R(x):=R_x$ is an
admissible representation of $\frak J$. Further, $\cal V$ is a
subalgebra of $\frak J$ and $\cal I^\perp$ is an ideal of $\frak
J$, then $\tilde R:{\cal V}\longrightarrow \mbox{End}({\cal
I}^\perp);\,\, v\longmapsto {\tilde{R}}(v):= R(v)/_{{\cal
I}^\perp},$ is a representation of $\cal V$ on $\cal I^\perp$.
Recall that, $\rho:{\cal V}\longrightarrow \mbox{End}({\cal
V}^*);\,\,x\longmapsto \rho(x)$ where $ \rho(x)(f)(y)=f(xy),\,\,
\forall x,y\in\cal  V$ is the coadjoint representation of $\cal
V$. $\nu$ is a representations isomorphism $(i.e. \,\,\,\nu\circ
\tilde R(v)=\rho(v)\circ \nu,\,\,\,\forall v\in \cal V)$ .

It is clear that the map $ \nabla:\cal{I}^{\perp}\oplus\cal
A^{\perp}\oplus \cal{V}\longrightarrow \cal{V}\oplus \cal
A^{\perp}\oplus \cal{V}^* $;
 $(i+a+v)\mapsto v+a+\nu(i),$
is an   isomorphism of Jordan algebras. In fact, recall that
${\cal{I}}\cap{\cal{I}^{\perp}}=\{ 0 \}$. Let $X=i+a+v , Y=j+b+w\in\frak{J}=
\cal{I}^{\perp}\oplus \cal A^{\perp}\oplus \cal{V}$ where $i,j\in \cal{I}^{\perp},$ $a,b\in
\cal A^{\perp},$ $v,w\in \cal{V}.$
$$XY=iw+jv+\alpha(a,b)+\beta(a,b)+\pi(w)a+\pi(v)b+vw.$$
\begin{eqnarray*}
 \nabla(XY)&=&vw+\beta(a,b)+\pi(w)a
+\pi(v)b+\nu(iw)+\nu(jv)+\nu(\alpha(a,b)).
\end{eqnarray*}
Moreover,
$$\nu(iw)(w')=B(iw,w')=B(i,ww')=\nu(i)(ww')=\Bigl(\nu(i)\circ\pi(w)\Bigl)(w').$$
and,
$$\nu(\alpha(a,b))(w')=B(\alpha(a,b),w')=B(\alpha(a,b)+\beta(a,b),w')
=B(ab,w')=\varphi(a,b),$$ where $\varphi(a,b)\in {\cal V}^{*}$
defined by $\varphi(a,b)(v):=Q(va,b)$. It follows that,
\begin{eqnarray*}
 \nabla(XY)=vw+\beta(a,b)+\pi(w)a
+\pi(v)b+\nu(i)\circ \pi(w)+\nu(j)\circ \pi(v)+\varphi(a,b)
=\nabla(X)\nabla(Y).&&
\end{eqnarray*}
Hence $\nabla$ is an algebras isomorphism. Further, if we consider on $\frak J$ the symmetric bilinear forms
$B':{\frak J}\times {\frak J}\longrightarrow \mathbb K,$ and
$B'':\frak J\times \frak J\longrightarrow \mathbb K$  defined by:
$$B'_{\vert_{{\cal I}\times {\frak J}}}=B/_{\vert_{{\cal I}\times {\frak J}}},\,\,\, B'_{\vert_{{\cal V}\times{\cal  V}}}=0\,\,\,\,\mbox{ and
}\,\,\,\, B''_{\vert_{{\cal V}\times {\frak J}}}=0,\,\,\,
B''_{\vert_{{\cal I}\times {\cal I}}}=0 , \,\,\,B''_{\vert{_{\cal
V}\times {\cal V}}}=B_{\vert_{{\cal V}\times {\cal V}}}, $$
then $B'$ is nondegenerate and associative, $B''$ is associative and $B=B'+B''.$\\
Now, let us consider on $\cal V\oplus \cal A^{\perp}\oplus \cal
V^*$ the associatif scalar products $T$ and $T'$  defined by:
\begin{eqnarray*}
&&T'(v+a+f,w+b+h)=Q(a,b)+f(b)+h(a),\, \,\forall v,w\in{\cal V},a,b\in {\cal A}^{\perp},f,h\in {\cal V}^{*} ,\\
&&T=T'+\gamma \hskip0.2cm
\mbox{ where }\hskip0.2cm\gamma:(\cal V\oplus \cal A^{\perp}\oplus \cal V^*)\times( \cal V\oplus \cal A^{\perp}\oplus
\cal V^*)\longrightarrow \mathbb{K}\\
&&\,\,\hskip 4.8cm(v+a+f,w+b+h)\,\,\,\longmapsto \,\,\,B(v,w).
\end{eqnarray*}
It follows that  $\nabla$ is an isometry of pseudo-euclidean
Jordan algebras from  $(\frak J,B) $
  (resp.$(\frak J,B')$) to  $({\cal V}\oplus{\cal  A}^{\perp}\oplus {\cal V}^*,{T})$ (resp.(${\cal V}\oplus
  {\cal A}^{\perp}\oplus {\cal V}^*,{T'})).$ \\It is clear that
$\phi:\cal V\oplus \cal A^{\perp}\oplus \cal V^*\longrightarrow
\cal V\oplus (\cal I/ \cal I^{\perp}) \oplus \cal V^*$;
$(v+a+f)\mapsto v+\theta(a)+f$ is an isomorphism  of Jordan algebras,
where $\cal V\oplus\cal A^{\perp}\oplus \cal V^*$ is the double
extension of $\cal A^{\perp}$ by $\cal V$ by means of $\pi $ and
$\cal V\oplus \cal I/ \cal I^{\perp}\oplus \cal V^*$ is the double
extension of $\cal I/\cal I^{\perp}$ by $\cal V$ by means of
$\tilde \pi: {\cal I}/{\cal I^{\perp}}\longrightarrow {\mbox
End}_s(V)$ defined by $\tilde \pi(s(i))(v):=\pi(vi),\forall i\in
I,v\in \cal V$. Now, if one considers the associatif scalar
products, $\Gamma$ and $\Gamma'$ on $\cal V\oplus \cal I/ \cal
I^{\perp}\oplus \cal V^*$ defined by:
\begin{eqnarray*}
&&\Gamma'(v+s(i)+f,w+s(j)+h):=\tilde B(s(i),s(j))+f(w)+h(v)=B(i,j)+f(w)+h(v)\\
&&\forall v,w\in {\cal V},\,\,i,j\in{\cal I},\,\, f,h\in {\cal V}^{*},
\end{eqnarray*}
\begin{eqnarray*}
\Gamma=\Gamma'+\gamma'\mbox{ where }
\gamma':(\cal V\oplus \cal I/\cal  I^{\perp}\oplus \cal V^* )\times (\cal V\oplus \cal I/ \cal I^{\perp}\oplus \cal V^*)\longrightarrow \mathbb{K}\hskip 2.2cm&&\\
\, (v+s(i)+f,w+s(j)+h)\,\,\,\longmapsto \,\,\,B(v,w), \,\,\hskip 1cm&&
\end{eqnarray*}
then $\phi$ is an isometry of pseudo-euclidean Jordan algebras from  $({\cal V}\oplus {\cal A}^{\perp}\oplus {\cal
V}^*,B)$ (resp. $({\cal V}\oplus {\cal A}^{\perp}\oplus {\cal
V}^*,B'))$ to  $({\cal V}\oplus ({\cal I}/ {\cal
I}^{\perp})\oplus{\cal  V}^*,\Gamma)$ (resp. $({\cal V}\oplus ({\cal
I}/ {\cal I}^{\perp})\oplus{\cal  V}^*,\Gamma').$ $\Box$
\begin{coro}
Let $(\frak{J},B) $ be an irreducible pseudo-euclidean Jordan algebra which is not simple. If $\frak J $ is not nilpotent,
 then $\frak J$ is a double extension of a pseudo-euclidean Jordan algebra $({\cal W},T)$ by a simple Jordan algebra.
\end{coro}
\dem $\frak J$ is not nilpotent, then ${\frak J}=\frak S\oplus
Ra(\frak J)$ where $\frak S$ is a semi-simple subalgebra of $\frak
J$ and $ Ra(\frak J) $ is the radical of $\frak J$. $\frak J$ is
irreducible and not simple, then $ Ra(\frak J)\neq \{0\}.$ Since
$\frak S$ is semi-simple, then $\frak S=\bigoplus^{m}_{i=1}\frak
S_i$, where ${\frak S}_i$ is a simple ideal of $\frak S ,$ for all
$ i\in \{1,...,m\}.$ The fact that ${\frak S}_1$ is a simple ideal
of $\frak J,$ implies that $I=\frak S_2\oplus ...\oplus \frak
S_n\oplus  Ra({\frak J})$ is a maximal ideal of $\frak J$. By the
last theorem, $\frak J$ is a double extension of $(W={\cal I}/
{\cal I^{\perp}},T= \tilde B)$ by $\frak S_1$ by means of the
representation $\pi:\frak S_1\longrightarrow \mbox{End}_s ({\cal
W},T)$ defined by: $\pi(x)(s(i)):=s(xi),\,\,\,\forall x\in{\frak
S}_1,\,i\in\cal I,$ where $s:{\cal I}\longrightarrow {\cal
I}/{\cal I}^{\perp}$ is the canonical surjection. Recall that
$\tilde B$ is defined by: $\tilde B(s(i),s(j)):=B(i,j),\,\,
\forall i,j\in {\cal I}.$ $\Box$
\begin{coro}\label{dx1}
Let $(\frak{J},B) $ be an irreducible pseudo-euclidean Jordan algebra which is not simple and such that
$\frak{J}\neq\{0\}$. If $Ann(\frak{J})=\{0\}$, then $\frak{J}$ is a double extension of a pseudo-euclidean Jordan
algebra $({\cal W},T)$ by a simple Jordan algebra.
\end{coro}
\dem If $Ann(\frak{J})=\{0\}$, then $Ann(\frak{J})^\perp
=\frak{J}$ (ie. $\frak{J}^2=\frak{J}$). Consequently, $\frak J$ is
not nilpotent. Hence, by the  Corollary above, $\frak{J}$ is a
double extension of a pseudo-euclidean Jordan algebra $({\cal
W},T)$ by a simple Jordan algebra. $\Box$
%%%%%%%%%%%%%%%%%%%%%%%%%%%%%%%%%%%%%%%%%%%%%%%%%%%\`u\`u
%%%%%%%%%%%%%%%%%%%%%%%%%%%%%%%%%%%%%%%%%%%%%%%%%%%%%
\begin{theor}\label{th3}
Let $(\frak J,B) $ be a pseudo-euclidean Jordan algebra which not
the one dimensional Jordan algebra with zero product. If $Ann(\frak
J)\neq \{0\}$ and if there exist $b\in Ann(\frak J)\setminus\{0\}$
such that $B(b,b)=\{0\}$, then $\frak J $ is a generalized double
extension of a pseudo-euclidean Jordan algebra $({\cal W},T)$ by
one dimensional Jordan algebra with zero product.
\end{theor}
\dem Let $b\in Ann({\frak J})\backslash \{0\}$ such that $
B(b,b)=0$. Consider the ideal ${\cal I}:=\mathbb Kb$ of $\frak J$.
There exists $a\in {\frak J} $  such that $B(a,b)=1,\,\,B(a,a)=0$
and $ {\frak J} ={\cal I^{\perp}}\oplus \mathbb Ka.$
 Denote ${\cal W}:= ( \mathbb Ka\oplus \mathbb Kb)^\perp,$ then ${\cal I^{\perp}}=\mathbb Kb\oplus {\cal W}$. It follows
  that,
 ${\frak J} =\mathbb Ka\oplus {\cal W}\oplus\mathbb Kb$. \\
Let $x,y\in  {\cal W} ,\,\, xy=\beta(x,y)+\alpha(x,y)b,$ where
$\beta(x,y)\in {\cal W}$ and $\alpha(x,y)\in \mathbb K.$ It is
easy to see that $\cal W$ endowed with the bilinear form $
\beta:\cal W\times \cal W\longrightarrow \cal \cal W;\,\,\, $
$(x,y)\longmapsto \beta(x,y), $ is a Jordan algebra. Moreover 
$B_{\cal W}=B_{\vert_{\cal W\times\cal W}}$ is an associatif
scalar product on $\cal W$ and $\alpha: \cal W\times \cal
W\longrightarrow \mathbb K$ is a symmetric bilinear form such that

\begin{eqnarray*}
\alpha\Bigl(\beta(x,y),\beta(x,y)\Bigl)=\alpha\Bigl(x,\beta(y,\beta(x,x))\Bigl)
\,\,\,\forall x,y\in {\cal W}.
\end{eqnarray*}
Now, if $ x\in {\cal W},$ then $ax= D(x)+\varphi(x)b$ where
$D(x)\in \cal W$ and $\varphi(x)\in \mathbb K$ or $D:\cal
W\longrightarrow \cal W;$ $x\longmapsto D(x)$ is an endomorphism
of ${\cal W}$ and $\varphi:{\cal W}\longrightarrow \mathbb K$ is
an element of the dual $\cal W^*$ of ${\cal W}.$ Since $B_{\cal
W}$ is nondegenerate, then there exist $ w_0\in \cal W$ such that
$\varphi= B_{\cal W}(w_0,.).$ Moreover, there exist $x_0\in {\cal
W},k,h\in \mathbb K$ such that $a^2=kb+x_0+ha$.\\ Let $x,y\in
{\cal W},\,\,B(xy,a)=B(x,ya)$ because $B$ is associative.
Consequently, $\alpha(x,y)=B(x,D(y)).$ On the other hand,
$B(x,ya)=B(xa,y).$ Which proves that  $B(x,D(y))=B(D(x),y).$ Then 
$D\in \mbox{End}_s({\cal W},B_{{\cal W}})$ and
$xy=\beta(x,y)+B(D(x),y),\,\, \forall x,y\in {\cal W}$. Now,
$B(a^2,b)=B(a,ab)=0$ which implies that, $hB(a,b)=0 $ (i.e. $
h=0).$ Hence, $a^2=w_0+kb.$\\ Let $x\in {\cal
W},B(ax,a)=B(x,a^2)$. Thus  $\varphi(x)=B(x,x_0).$ We conclude
that,
$w_0=x_0.$ Consequently, $ax=D(x)+B(x,x_0)b.$\\
It is easy to show that
$(D,x_0)\in \mbox{End}_s({\cal W},B_{{\cal W}})\times \cal W$ is an admissible pair of the pseudo-euclidean Jordan algebra $
({\cal W} ,B_{\cal W}).$ Consequently, $(\frak J ,B)$ is a generalized double extension of  $ ({\cal W} ,T:=B_{\cal W})$
by the one dimensional Jordan algebra with zero product  $\mathbb Ka$ by means of the pair $(D,x_0)\in \mbox{End}_s({\cal W},T)\times
{\cal W}. $  $\Box$

\begin{coro} Let $(\frak J,B) $ be an irreducible nilpotent pseudo-euclidean Jordan algebra. If $ \frak J $ is not the one
 dimensional Jordan algebra with zero product, then
$\frak J $ is the generalized double extension of a nilpotent
pseudo-euclidean Jordan algebra
 $({\cal W},T)$ by the one dimensional Jordan algebra with zero product.
\end{coro}
\dem If $\frak J$ is nilpotent, then $\frak J ^2\neq \frak J $.
 Consequently $Ann(\frak J )\neq\{0\}.$ By the last Theorem $\frak J $ is
a generalized double extension of a pseudo-euclidean Jordan
algebra $({\cal W},T)$ by the one dimensional Jordan algebra to
null product. In the proof of the same Theorem,  $\cal W =\cal
I^{\perp}/\cal I$ where $ {\cal I}=\mathbb Kb\subset Ann(\frak J
)$. Since $\frak J $ is nilpotent, then ${\cal I}^{\perp}$ is
nilpotent. Consequently, $\cal W $ is nilpotent. $\Box$

Let $\cal U$ be the set constituted with $\{0\},$ the one
dimensional Jordan algebra with zero product and all simple Jordan
algebras.
\begin{theor}\label{calas}
Let $(\frak J,B) $ be a pseudo-euclidean Jordan algebra. If $\frak
J \notin \cal U$, then $\frak J $ is obtained from    elements
${\frak J}_1,...,{\frak J}_n$ of $\cal
 U$, by a finite number of orthogonal direct sums of pseudo-euclidean Jordan algebras or/and double extensions by a simple Jordan algebra or/and generalized double extension by the one dimensional Jordan algebra with zero product.
\end{theor}
\dem
We proceed by induction on $dim{\frak J}$.
If $\mbox{dim}{\frak J}=0 $ or $1,$ then $\frak J \in \cal U.$\\
Assume that $\mbox{dim}{\frak J}=2.$ If $\frak J $ is not irreducible, then $\frak J ={\cal{I}}_1 \oplus {\cal{I}}_2$ where ${\cal{I}}_1,{\cal{I}}_2$ are two nondegenerate ideals of $\frak J $ which satisfies $B({\cal{I}}_1,{\cal{I}}_2)=\{0\}$ and $\mbox{dim}{\cal{I}}_1=\mbox{dim}{\cal{I}}_2=1.$ Therefore ${\cal{I}}_1$
  (resp.${\cal{I}}_2$) is either a one dimensional simple Jordan algebra or the one dimensional algebra with zero product. Now, suppose that $\frak J$ is irreducible. If $\frak J $ is neither simple nor nilpotent, then $\frak J ={\cal S} \oplus Ra(\frak J )$ where
${\cal S}$ is a one dimensional semi-simple Jordan subalgebra and
 $\mbox{dim}Ra(\frak J )=1.$ In this case, $\frak J ={\cal S}
 \oplus {\cal S}^*$ is the double extension of $\{0\}$ by ${\cal S}$ by means of  the representation $0$. Now if $\frak J $ nilpotent, then $\frak J $ is the generalized double extension of $\{0\}$ by the one dimensional Jordan algebra with zero product. We conclude that if $\mbox{dim}\frak J =2$, the theorem is satisfied.\\
Now suppose that the theorem is satisfied for $\mbox{dim}({\frak
J})<n\in\mathbb N$. We shall prove it  in the case
where $\mbox{dim}{\frak J}=n$.  If $\frak J \notin \cal U$ and 
$\frak J $ is irreducible, then $\frak J $ is either a generalized double
extension of a pseudo-euclidean Jordan algebra $({\cal W},T)$ by
the one dimensional Jordan algebra with zero product or a double
extension of a pseudo-euclidean Jordan algebra $({\cal W},T)$ by
a simple Jordan algebra. Since $\mbox{dim}{\cal W}<
\mbox{dim}{\frak J},$ then $({\cal W},T)$ satisfies the theorem, so
$\frak J $ satisfies the theorem. Now, if $\frak J \notin \cal U$  and $\frak J $ is not
irreducible, then ${\frak J} ={\cal I}_1\oplus ...\oplus {\cal
I}_m$ where ${\cal I}_i\not= \{0\}, 1\leq i\leq m$ are   nondegenerate irreducible
ideals of $\frak J $ such that  $B({\cal I}_i,{\cal
I}_j)=\{0\}\,\,\, \forall i\neq j\in\{1,...m\}.$  Let $i\in\{1,...m\},$ the fact that 
$\mbox{dim}{\cal I}_i<\mbox{dim}{\frak J}$  implies that  $({\cal I}_i,B/_{{\cal I}_i\times {\cal
I}_i})$ satisfies the theorem. Hence, $({\frak J},B)$ satisfies
the theorem. $\Box$

Let $\cal E$ be the set constituted by $\{0\}$ and the one
dimensional Jordan algebra with zero product.
\begin{theor}
Let $(\frak J ,B) $ be a nilpotent pseudo-euclidean Jordan algebra. If $\frak J \notin \cal E,$ then $\frak J $ is obtained from   elements
${\frak J}_1,...,{\frak J}_n$ of $\cal E$ by a finite number of orthogonal direct sums of pseudo-euclidean Jordan algebras or/and generalized double extension by the one dimensional Jordan algebra with zero product.
\end{theor}
\dem
We proceed by induction on $dim{\frak J}$.
If $\mbox{dim}{\frak J}=0 $ or $1,$ then $\frak J \in \cal E.$\\
If $\mbox{dim}{\frak J}=2$, then $\frak J$ is either the orthogonal
direct sum of two one dimensional Jordan algebras with zero product
or $\frak J $
is the generalized double extension of $\{0\}$ by the one dimensional Jordan algebra with zero product.\\
Now, suppose that the theorem is satisfied for $\mbox{dim}{\frak J}<n$. We shall prove it  in the case where $\mbox{dim}{\frak J}=n$.  Suppose that $\frak J \notin \cal E$ and  $\frak J $ is irreducible. since $\frak J $ is nilpotent, then
$\frak J $ is the generalized double extension of a nilpotent
pseudo-euclidean Jordan algebra $({\cal W},T)$ by the one
dimensional algebra with zero product. Since $\mbox{dim}({\cal W})<
\mbox{dim}{\frak J},$ then $\cal W$ satisfies the theorem. Then 
$\frak J $ satisfies the theorem. If $\frak J \notin \cal E$ and $\frak J $ is not
irreducible, then ${\frak J} ={\cal I}_1\oplus ...\oplus {\cal
I}_m$ where ${\cal I}_i, 1\leq i\leq m,$ are   nilpotent nondegenerate irreducible
ideals of
 $\frak J $ which satisfy
$B({\cal I}_i,{\cal I}_j)=\{0\},\,\,\, \forall i\n ot=j\in\{1,...m\}.$ Since $\mbox{dim}{\cal
I}_i<\mbox{dim}{\frak J},\,\,\forall i\in\{1,...m\}$, then, any
$({\cal I}_i,B/_{{\cal I}_{i}\times {\cal I}_i})$ satisfies the
theorem. Hence, $(\frak J,B)$ satisfies the theorem.  $\Box$

\section{Nilpotent pseudo-euclidean Jordan algebras with dimension less than or
equal to $5$}\label{5}

In this section, we shall use   inductive description obtained in section $4$ to construct the nilpotentpseudo-euclidean Jordan algebras with dimensions $ n\leq 5.$ In particular, we shall prove that all  nilpotent pseudo-euclidean Jordan algebras with dimensions $n\leq4$   are
associative, but  there are  nilpotent pseudo-euclidean Jordan
algebras with dimension equal to $ 5$ which are nonassociative.

First case: $n=1$. The product of nilpotent one dimensional Jordan algebra vanish
(see $\cite {scha}$), (i.e. $ \frak{J}_{1,1}=\mathbb{K}a $ where
$a^2=0$). The associatif scalar product of $\frak{J}_{1,1}$ is given by $B_{1}(a,a)=1$.

Second case:  $n=2$. The generalized double extension of null algebra is the nilpotent pseudo-euclidean Jordan algebra $\frak{J}_{2}=Vect\{a_1,b_1\} $ where $b_1\in Ann(\frak{J}_{2})\mbox{
and } a_1a_1=\lambda b_1\mbox{ where } \lambda\in \mathbb K.$

If $\lambda=0$, we obtain $\frak{J}_{2,0}=Vect\{a_1,b_1\}$ where
$\frak{J}_{2,0}^2=\{0\}.$ If $\lambda\neq 0$, we obtain
$\frak{J}_{2,\lambda}$ such that $b_1\in Ann(\frak{J}_{2})\mbox{
and } a_1a_1=\lambda b_1$. Now, let $\lambda\in \mathbb K^{*},$
then $\varphi_2:\frak{J}_{2,\lambda}\longrightarrow\frak{J}_{2,1}$
defined by $$\varphi_2(a_1)=a_1 \mbox{ and }
\varphi_2(b_1)=\frac{1}{\lambda_{1}}b_1$$  is an algebras
isomorphism. Thus, for all $\lambda\neq 0$ ,
$\frak{J}_{2,\lambda}$ is isomorph to $\frak{J}_{2,1}.$
 The associatif scalar product $B_{2}$ on $\frak{J}_{2,0}$ and on $\frak{J}_{2,1}$ is given by $B_{2}(a_1,b_1)=1,$ $B_{2}(a_1,a_1)=B_{2}(b_1,b_1)=0$.

Third case: $n=3$.
 Let us consider the algebra
$\frak{J}_{1,1}=\mathbb{K}a $, where $aa=0$. The set of admissible pairs of this algebra is given by:
$$\{(D,x_0)\in\mbox{End}_s(\frak{J}_{1,1})\times\frak{J}_{1,1};\,\,\,D=0,\,\,x_0=\alpha
a,\,\, \mbox{where } \alpha\in\mathbb K\}.$$ If $\alpha=0$, the
nilpotent pseudo-euclidean Jordan algebra with dimension equal to
$3$  obtained by generalized double extension of $\frak{J}_{1,1}$
by means of the pair $(D,x_0)=(0,0)$ is
$\frak{J}_{3,0,k}=\mathbb{K}a_{2}\oplus \mathbb{K}a \oplus
\mathbb{K}b_{2}$ where $\{a,b_{2}\}\in Ann({\frak J}_{3,0,k})$ and
$a_{2}a_{2}=kb_{2},\,\,\, k\in\mathbb{K}$. If $k=0$, then ${\frak
J}_{3,0,0}=\mathbb{K}a_{2}\oplus \mathbb{K}a \oplus
\mathbb{K}b_{2}$ where ${\frak J}_{3,0,0}^2=\{0\}$. Let
$k\in{\mathbb K}^{*},$ then the map $\varphi_3:{\frak
J}_{3,0,k}\longrightarrow{\frak J}_{3,0,1}$ defined by
$$\varphi(a_2)=a_2,\,\,\,\varphi(a)=a\mbox{ et }
\varphi(b_2)=\frac{1}{k}b_2,$$ is a isomorphism of Jordan algebras.
Thus for all $k\in{\mathbb{K}}^{*} $, ${\frak J}_{3,0,k}$ is
isomorph to ${\frak J}_{3,0,1}$. Now if $\alpha\neq 0,$ then the
nilpotent pseudo-euclidean Jordan algebra with dimension equal to
$3$ obtained by generalized double extension of $\frak{J}_{1,1}$
by means of $(D,x_0)=(0,\alpha a)$ is given by
$\frak{J}_{3,\alpha,k}=\mathbb{K}a_{2}\oplus \mathbb{K}a \oplus
\mathbb{K}b_{2}$ where $b_{2}\in Ann({\frak J}_{3,\alpha,k}),$
$aa=0$, $a_{2}a=\alpha b_{2}$ and $a_{2}a_{2}=\alpha
a+kb_{2},\,\,\, k\in\mathbb{K}$. Let $k\in
{\mathbb{K}},\,\,\,\alpha\in {\mathbb{K}}^{*}$. Let us consider
the map
$\varphi:\frak{J}_{3,\alpha,k}\longrightarrow\frak{J}_{3,1,0}$,
defined by, $$\varphi(a_2)=\alpha a_2+ \frac{k}{2}a,
\,\,\,\varphi(b_2)=\alpha b_2,\,\,\, \varphi(a)=\alpha a.$$
$\varphi$ is an isomorphism of Jordan algebras. Thus, $\forall \alpha\in
{\mathbb{K}}^{*},\forall k\in {\mathbb{K}},\,\,\,
\frak{J}_{3,\alpha_1,k}$ is isomorph to $\frak{J}_{3,1,0}.$ The
associatif scalar product on $\frak{J}_{3,0,0}$, ${\frak
J}_{3,0,1}$ and on ${\frak J}_{3,1,0}$ is given by
\begin{eqnarray*}
&&B(a_{2},b_{2})=1 \hskip 1.5 cm B(a_{2},a_{2})=B(b_{2},b_{2})=0\\
&&B(a,a)=1 \hskip 1.5 cm B(a_{2},a)=B(b_{2},a)=0
\end{eqnarray*}
 Fourth case: $n=4$. Recall that in the second case, we have proved that the two dimension nilpotent
 pseudo-euclidean Jordan algebras are
  $\frak{J}_{2,0}=Vect\{a_1,b_1\}$
where$\frak{J}_{2,0}^2=\{0\}$ and $\frak{J}_{2,1}=Vect\{a_1,b_1\}$
where $b_1\in A(\frak{J}_{2})\mbox{ and } a_1a_1= b_1.$ Moreover,
the bilinear form defined on $\frak{J}_{2,0}$ (resp.
$\frak{J}_{2,1}$) by $B_{2}(a_1,b_1)=1,$
$B_{2}(a_1,a_1)=B_{2}(b_1,b_1)=0$ is an associatif scalar product on $\frak{J}_{2,0}$ (resp. $\frak{J}_{2,1}$).\\

 Let us consider ${\frak J}_{2,0}=\mathbb
Ka_1\oplus\mathbb Kb_1$, where ${\frak J}_{2,0}^2=\{0\}$. Denote
by, $A$ and $B$ the subsets of $End({\frak J}_{2,0})$ defined by:
\begin{eqnarray*}
&&A=\{D\in End({\frak J}_{2,0}); \,\,D(a)=\alpha b \mbox{ and } D(b)=0,\mbox{ where }\alpha\in\mathbb K\},\\
&&B=\{D\in End({\frak J}_{2,0}); \,\,D(a)=0 \mbox{ and } D(b)=\alpha
a,\mbox{ where }\alpha\in\mathbb K\}.
\end{eqnarray*}
The pair $(D,x_0)\in End({\frak J}_{2,0})\times{\frak J}_{2,0}$ is
an admissible pair of ${\frak J}_{2,0}$, if and only if
$(D,x_0)\in (A\cup B)\times{\frak J}_{2,0}$. Let $ D\in A$ (resp.
$\in B$) defined by $D(a_1)=\alpha b_1$ and $D(b_1)=0$ (resp.
$D(a_1)=0$ and $D(b_1)=\alpha a_1)$ and $x_0=\eta a_1+\varepsilon
b_1\in{\frak J}_{2,0}$, where $\alpha,\eta,\varepsilon\in\mathbb
K$. Let $k\in\mathbb K$. Then, the product in the nilpotent
pseudo-euclidean Jordan algebra ${\frak J}_{4,0}=\mathbb
Ka_3\oplus\mathbb Ka_1\oplus\mathbb Kb_1\oplus\mathbb Kb_3$
obtained by the generalized double extension of ${\frak J}_{2,0}$
by the one dimensional Jordan algebra with zero product by means of
$(D,x_0)\in A\times{\frak J}_{2,0}$ (resp. $\in B\times{\frak
J}_{2,0}$) is given by:
$$\begin{array}{llll}
b_3\in Ann({\frak J}_{4,0}),& a_3a_3=\eta a_1+\varepsilon b_1+kb_3,& a_3b_1=\eta b_3,\\
a_3a_1=\alpha b_1+\varepsilon b_3,& a_1a_1=b_1+\alpha b_3,&
a_1b_1=b_1b_1=0.
\end{array}$$
resp.
$$\begin{array}{llll}
b_3\in Ann({\frak J}_{4,0}),& a_3a_3=\eta a_1+\varepsilon b_1+kb_3,& a_3b_1=\alpha a_1+\eta b_3,\\
a_3a_1=\varepsilon b_3,& a_1a_1=b_1,\hskip 0.7cm a_1b_1=0,
&b_1b_1=\alpha b_3.
\end{array}$$
Now let us consider the Jordan algebra $\frak{J}_{2,1}$ obtained  in
the second case. $(D,x_0)\in \mbox{End}(\frak{J}_{2,1})\times
\frak{J}_{2,1}$ is admissible, if and only if, there exists
$\beta,\varepsilon \in \mathbb K,$ such that
$$D(a_1)=\beta b_1,\,\,\,D(b_1)=0,\mbox{ and } x_0=\varepsilon b_1.$$ Let $(D,x_0)$ be an admissible pair of $\frak{J}_{2,1}$. Then, the product on the nilpotent pseudo-euclidean Jordan algebra ${\frak J}_{4,1}=\mathbb Ka_3\oplus\mathbb Ka_1\oplus\mathbb Kb_1\oplus\mathbb Kb_3$ obtained by the generalized double extension of ${\frak J}_{2,1}$ by the one dimensional Jordan algebra with zero product by means of $(D,x_0)$ is given by:
 $$\begin{array}{llll}
b_1,b_3\in Ann({\frak J}_{3}), &a_3a_3=\varepsilon b_1+kb_3,&
a_3a_1=\beta b_1+\varepsilon b_3,&a_1a_1=b_1+\beta b_3
\end{array}$$
An associatif scalar product on ${\frak J}_{4,0}$ and ${\frak
J}_{4,1}$ is given by:
 $$\begin{array}{llll}
B_4(a_1,a_1)=B_4(b_1,b_1)=0,&B_4(a_1,b_1)=B_4(a_3,b_3)=1,&
B_4(a_3,a_3)=B_4(b_3,b_3)=0,\\B_4(a_3,a_1)=B_4(a_3,b_1)=0,&
B_4(b_3,a_1)=B_4(b_3,b_1)=0
\end{array}$$
 Fifth case: $n=5$.  We shall construct with the same process all nilpotent pseudo-euclidean Jordan algebras with dimension equal to $5$. In the first time, let us consider the nilpotent pseudo-euclidean Jordan algebra
 ${\frak
J}_{3,0,0}$ and let $(D,x_0)\in End({\frak J}_{3,0,0})\times
{\frak J}_{3,0,1}$, which satisfies the conditions $(C_1),\dots,(C_7)$. These conditions are equivalent to $D^3(x)=0,\,\,\,\forall x\in{\frak J}_{3,0,0}$ and
$D^2(x_0)=0$.

In the second time, let us consider  ${\frak
J}_{3,0,1}=\mathbb{K}a_{2}\oplus \mathbb{K}a\oplus
\mathbb{K}b_{2}$. The product on  ${\frak J}_{3,0,1}$ is defined
by $a_2a_2=b_2$ and $a,b_2\in A({\frak J}_{3,0,1})$. $(D,x_0)\in
End({\frak J}_{3,1,0})\times {\frak J}_{3,0,1}$ is admissible, if
and only if there exists,
$\alpha,\beta,\gamma,\eta_1,\eta_2,\eta_3$ in $\mathbb K,$ such
that
$$D(a_2)=\alpha a+\beta b_2,\,\,\,D(a)=\gamma b_2,\,\,\,D(b_2)=0\mbox{ and } x_0=\eta_1a_2+\eta_2 a+\eta_3 b_2.\,\mbox{ where }\eta_1\in\{0,\alpha^2\}. $$

Let $(D,x_0)$ an admissible pair of ${\frak J}_{3,0,1}$. The product on the Jordan algebra ${\frak J}_{5,0,1}=\mathbb
Ka_4\oplus\mathbb Ka_2\oplus\mathbb Ka\oplus\mathbb
Kb_2\oplus\mathbb Kb_4$, obtained by the generalized double extension
of ${\frak J}_{3,0,1}$ by the one dimensional Jordan algebra with zero product is given by: $b_4\in A({\frak
J}_3),$
$$\begin{array}{llll}
a_4a_4=\eta_1 a_2+\eta_2 a+\eta_3b_2+kb_4,&\hskip 0.4cm
aa=b_2a_2=b_2a=b_2b_2=0,&\hskip 0.4cm b_2a_4=\eta_1b_4\\
a_2a=\alpha_2b_4\hskip 1.cm aa_4=\alpha_2b_2,&\hskip 0.5cm
 a_2a_4=\alpha_2 a+\alpha_3 b_2+\eta_3b_4,&\hskip 0.5cma_2a_2=b_2+\alpha_3b_4.\\
\end{array}$$
Now let us consider ${\frak J}_{3,1,0}=\mathbb{K}a_{2}\oplus
\mathbb{K}a\oplus \mathbb{K}b_{2}$. The product on ${\frak
J}_{3,0,1}$ is defined by: $$b_2\in A({\frak
J}_{3,1,0}),\,\,\,a_2a_2=a,\,\,\,a_2a=b_2.$$ $(D,x_0)\in
End({\frak J}_{3,1,0})\times {\frak J}_{3,1,0}$ is admissible, if
and only if, there exists $\alpha,\beta,\gamma,\varepsilon,\eta$
in $\mathbb K,$ such that
$$D(a_2)=\alpha a+\beta b_2,\,\,\,D(a)=\gamma b_2,\,\,\,D(b_2)=0\mbox{ and } x_0=\varepsilon a+\eta b_2. $$
Let $(D,x_0)$ be an admissible pair of ${\frak J}_{3,1,0}$. Then,
the product on the Jordan algebra ${\frak J}_{5,1,0}=\mathbb
Ka_4\oplus\mathbb Ka_2\oplus\mathbb Ka\oplus\mathbb
Kb_2\oplus\mathbb Kb_4$, obtained by the generalized double extension
of ${\frak J}_{3,1,0}$  by the one dimensional Jordan algebra with zero product is given by:
$$\begin{array}{llll}
b_2,b_4\in Ann({\frak J}_{5,0,1}), & \hskip 0.7cm a_4a_4=\varepsilon
a+\eta b_2 +kb_4, & \hskip 0.7cm
a_2a_2=b_2+\beta b_4,\\
a_2a=\gamma b_4,\hskip 0.3cm aa=0,& \hskip 0.7cm a_4a_2=\alpha
a+\beta b_2+\eta b_4, & \hskip 0.7cm a_4a=\gamma b_2,&
\end{array}$$

Now let us consider the pseudo-euclidean Jordan algebra ${\frak
J}={\frak J}_{5,1,0}$ where $k=\alpha_3=\eta_2=\eta_3=0$ and
$\alpha_2=0.$ The product on ${\frak J}$ is given by:
$$b_2,b_4\in Ann({\frak J}),\hskip 0.5cm a_4a_4=aa=0,\hskip 0.5cm a_2a_2=a,\hskip 0.5cm a_2a=b_2+b_4,\hskip 0.5cm a_4a_2=a,\hskip 0.5cm a_4a=b_2.$$
The fact that $a_2(a_2a_4)-(a_2a_2)a_4=b_4\neq0$   implies that ${\frak J}$ is not associative.
\section{Symplectic Forms and Jordan bialgebras}\label{big}

In this section we  study the relation between solutions of 
Yang Baxter equation (EYB) in the case of Jordan algebras (see
$\cite{zhel1}$) and     symplectic structures   on pseudo euclidean Jordan
algebras. We study also the connection between symplectic pseudo euclidean Jordan
algebras and some symplectic quadratic Lie algebras.

\begin{defi} \label{def2}
 Let $\omega$ be a skew-symmetric
non-degenerate bilinear form on $\frak J$ which satisfies
$\omega(xy,z)+\omega(yz,x)+\omega(zx,y)=0,\,\,\,\forall
x,y,z\in{\frak J}$. $\omega$ is called a symplectic form and the
pair $({\frak J},\omega)$ is called symplectic Jordan algebra and
if moreover $B$ is an associatif scalar product on $\frak J$, then
$({\frak J},B,\omega)$ is called a symplectic pseudo-euclidean
Jordan algebra.
\end{defi}
\begin{defi}
Let ${\cal V}$ be a vector space and $\Delta:{\cal
V}\longrightarrow{\cal V}\otimes{\cal V}$ be a map. Then the pair
$({\cal V},\Delta)$ is called a coalgebra and $\Delta$ is called a
comultiplication.
\end{defi}
For $v\in{\cal V},$ one writes: $\Delta(v)=\sum_vv_{(1)}\otimes
v_{(2)}$. The dual space ${\cal V}^*$ of ${\cal V}$ endowed by the
following product:
$$<fg,v>=\sum_{v}<f,v_{(1)}><g,v_{(2)}>,\,\,\,\forall f,g\in{\cal V}^*,\,\,v\in{\cal V}.$$ is an algebra called the dual of the coalgebra $({\cal V},\Delta)$. Recall that if $f,g\in {\cal V}^*,$ we can consider $f\otimes g \in ({\cal V}\otimes {\cal V})^*$ defined by:
 
$$<f\otimes g,\sum_ia_i\otimes b_i>=\sum_i<f,a_i><g,b_i>,\,\,\,  \,a_i,b_i\in{\cal V},$$

which proves that ${\cal V}^*\otimes {\cal V}^* \subseteq ({\cal V}\otimes {\cal V})^*.$
 
The space ${\cal V}$ endowed with the following actions is a ${\cal V}^*-$bimodule:
\begin{eqnarray*}
&&f.v=\sum_{v}v_{(1)}<f,v_{(2)}>\mbox{ and
}v.f=\sum_{v}<f,v_{(1)}>v_{(2)},\\&&\mbox{ where } f\in {\cal
V}^*,\,\,\,v\in{\cal V}\mbox{ and
}\Delta(v)=\sum_{v}v_{(1)}\otimes v_{(2)}.
\end{eqnarray*}
It is clear that for all $f,g\in{\cal V}^*$ and $v\in{\cal V}$ the
equalities $<fg,v>=<f,g.v>=<g,v.f>,$ hold for all $f,g\in{\cal
V}^*$ and $v\in{\cal V}$. Now, let ${\cal V}$ be an algebra,
$\Delta$ be a comultiplication of ${\cal V}$ and ${\cal V}^*$ the
dual algebra of the coalgebra $({\cal V},\Delta)$.  ${\cal V}^*$ with the following actions is a ${\cal V}-$bimodule:

$$<f\ast v,w>=<f,vw>=<w\ast f,v> ,\,\,\,\forall f\in{\cal V}^*,\,\,v,w\in{\cal V}.$$
Consider the space $D({\cal V})={\cal V}\oplus{\cal V}^*$ on which we define the following product:
$$(v+f)(w+g)=(vw+f.w+v.g)+(fg+f\ast w+v\ast g),\,\,\,\forall f,g\in{\cal V}^*,\,\,v,w\in{\cal V}.$$
Then, $D({\cal V})$ endowed with the product above is an algebra
in which ${\cal V}$ and ${\cal V}^*$ are two subalgebras. The
algebra $D({\cal V})$ is called the (Drinfeld-)double of the
bialgebra $({\cal V},m,\Delta)$ where $m:{\cal V}\times{\cal V}
\longrightarrow{\cal V}$ defined by $m(v,w)=vw,\,\,\,\forall
v,w\in{\cal V}$ (i.e. $m$ is the multiplication in ${\cal V}).$

\begin{defi}
Let ${\frak J}$ be a Jordan algebra not necessarily unital
and $\Delta$ be a comultiplication of ${\frak J}$. $({\frak
J},\Delta)$ is a Jordan bialgebra  if $D({\frak J})$ is a Jordan
algebra.
\end{defi}

Let ${\frak J}$ be a Jordan algebra and consider the linear map
 \begin{eqnarray*}
\tau:{\frak J}\otimes{\frak J}&\longrightarrow&{\frak J}\otimes{\frak J}\\
\sum_{i}x_i\otimes y_i&\longmapsto&\sum_{i}y_i\otimes x_i.
\end{eqnarray*}
Let $r=\sum_ia_i\otimes b_i\in{\frak J}\otimes{\frak J}$ such that
$\tau(r)=-r$ (i.e. $r$ is skew-symmetric). We define on  ${\frak
J}$ the comultiplication $\Delta_r$ by
$\Delta_r(x)=\sum_{i=1}^na_ix\otimes b_i-a_i\otimes
xb_i,\,\,\,\forall x\in{\frak J}$. Let
$$C_{\frak J}(r)=r_{12}r_{13}-r_{12}r_{23}+r_{13}r_{23},\hskip 0.3cm where$$
\begin{eqnarray*}
&&r_{12}r_{13}=\sum_{1\leq i,j\leq n}a_ia_j\otimes b_i\otimes
b_j,\\&&r_{12}r_{23}=\sum_{1\leq i,j\leq n}a_i\otimes b_ia_j\otimes
b_j=-\sum_{1\leq i,j\leq n}b_j\otimes a_ia_j\otimes b_i\\
&&and\,\,\,\, r_{13}r_{23}=\sum_{1\leq i,j\leq n}a_i\otimes
a_j\otimes b_ib_j=\sum_{1\leq i,j\leq n}b_i\otimes b_j\otimes
a_ia_j
\end{eqnarray*}
because $r$ is antisymmetric. Thus,
$$C_{\frak J}(r)=\sum_{1\leq i,j\leq n}a_ia_j\otimes b_i\otimes b_j+b_i\otimes b_j\otimes a_ia_j
+b_j\otimes a_ia_j\otimes b_i.$$ $C_{\frak J}(r)$ is an element of
${\frak J}\otimes{\frak J}\otimes{\frak J}$ and it is well defined
even if ${\frak J}$ is not unital. The Jordan Yang Baxter equation
of ${\frak J}$ is $C_{\frak J}(r)=0$ (see $\cite{zhel})$. In this
case, $r$ is said an antisymmetric   solution of The Jordan Yang Baxter equation or $r$ is an
antisymmetric $r-$matrix. In $\cite{zhel}, $ it is proven that if
$r\in{\frak J}\otimes{\frak J}$ is an antisymmetric $r-$matrix, then the pair (${\frak J},\Delta_r)$ is a Jordan
bialgebra.

Let $r$ be an antisymmetric $r-$matrix. One poses the linear map
$R:{\frak J}^*\longrightarrow{\frak J}$ defined by
$R(f)=\sum_{i}f(a_i)b_i,\,\,\,\forall f\in{\frak J}^*$. The
equality $C_{\frak J}(r)=0$ is equivalent to the following
identity: $$<f,R(h)R(l)>+<h,R(l)R(f)>+<l,R(f)R(h)>=0.$$ In fact,
\begin{eqnarray*}
&&(f\otimes h\otimes l)(r_{12}r_{13})=\sum_{i,j}f(a_ia_j) h(b_i)
l(b_j)=f\Bigl(h(b_i)a_il(b_j)a_j\Bigl)=<f,R(h)R(l)>,
\\&&(f\otimes h\otimes l)(r_{12}r_{23})=\sum_{i,j}f(a_i)h(b_ia_j)l(b_j)=-<h,R(l)R(f)>,
\\&&(f\otimes h\otimes l)(r_{13}r_{23})=\sum_{i,j}f(a_i)h(a_j)l(b_ib_j)=<l,R(f)R(h)>.
\end{eqnarray*}
Moreover, since $r$ is antisymmetric, then $R$ is
antisymmetric. In fact, let $f,h\in {\frak J}^*$.
$$<f,R(h)>=-f(\sum_{i=1}^nh(b_i)a_i)=-h(\sum_{i=1}^nf(a_i)b_i)=-<h,R(f)>=-<R(f),h>.$$
In this case we say that $R$ satisfies the Yang Baxter equation.
If $R$ is bijective, we say that $r$ is a nondegenerate
$r-$matrix.
\begin{prop}
Let $({\frak J},B)$ be a pseudo-euclidean Jordan algebra and
$r=\sum_{i=1}^na_i\otimes b_i$ be an antisymmetric $r-$matrix of
${\frak J}$. Let $\phi:{\frak J}\longrightarrow{\frak J}^*$ be a
linear isomorphism defined by: $\phi(x)=B(x,.),\,\,\,\,\forall
x\in{\frak J}$ and $R:{\frak J}^*\longrightarrow{\frak J}$ the
linear map defined by: $R(f)=\sum_{i=1}^nf(a_i)b_i,\,\,\,\forall
f\in{\frak J}^*$. Then, ${\cal U}=R\circ \phi$ is
$B-$antisymmetric. Further, ${\frak J}$ endowed with the product
$\star$ defined by $x\star y={\cal U}(x)\,\,y+x\,\,{\cal
U}(y),\,\,\,\forall x,y\in{\frak J}$ is a Jordan algebra for which
${\cal U}:({\frak J},\star)\longrightarrow{\frak J}$ is a Jordan
algebra isomorphism.
\end{prop}
\dem Let $x,y,z\in{\frak J}.$ Consider  $f:=\phi(x)=B(x,.)$,
$h:=\phi(y)=B(y,.)$ and $l:=\phi(z)=B(z,.)$.
\begin{enumerate}
\item $\,\,\,x\star y=x\,\,{\cal U}(y)+y\,\,{\cal U}(x)=x(R\circ\phi)(y)+y(R\circ\phi)(x)=xR(h)+yR(f)$\\

$\hskip 1.3cm=\phi^{-1}(f\circ  R_{R(h)})+\phi^{-1}(h\circ  R_{R(f)})=\phi^{-1}(\phi(x)\phi(y)).$\\

Hence, $\phi(x\star y)=\phi(x)\phi(y)$. Since ${\frak J}^*$ is a
Jordan algebra and $\phi$ is bijective, then $( {\frak J},\star)$
is isomorph to
 ${\frak J}^*$.

\item Since $R$ is antisymmetric, then $<\phi(x),R(\phi(y))>=-<\phi(y),R(\phi(x))>$. Hence, $B(x,{\cal U}(y))=-B(y,{\cal U}(x))$. Thus ${\cal U}$ is $B-$antisymmetric.

\item $C_{\frak J}(r)=0$. Thus, $<f,R(h)R(l)>+<h,R(l)R(f)>+<l,R(f)R(h)>=0.$ Thus, $B(x,{\cal U}(y)\,{\cal U}(z))+B(y,{\cal U}(z)\,{\cal U}(x))+B(z,{\cal U}(x)\,{\cal U}(y))=0$.
It follows that, $B(-{\cal U}(x\,{\cal U}(y))-{\cal U}(y\,{\cal
U}(x))+{\cal U}(x){\cal U}(y),z)=0.\mbox{ Hence, } {\cal
U}(x\,\,{\cal U}(y))+{\cal U}(y\,\,{\cal U}(x))={\cal U}(x){\cal
U}(y).$ Consequently, ${\cal U}(x\star y)={\cal U}(x){\cal U}(y)$.
Thus ${\cal U}$ is a Jordan algebra  morphism.$\Box$
\end{enumerate}

\begin{rema}
The converse of the last proposition is true. If ${\cal U}$ is a
$B-$antisymmetric endomorphism of ${\frak J}$ such that the
product $\star$ given by: $x\star y={\cal U}(x)\,y+x\,{\cal
U}(y),\,\,\,\forall x,y\in{\frak J},$ define a Jordan structure on
${\frak J}$, then the linear map $R:{\frak
J}^*\longrightarrow{\frak J}$ defined by $R={\cal
U}\circ\phi^{-1}$, where $\phi(x)=B(x,.),\,\,\,\forall x\in {\frak
J}$, is a solution of the Yang Baxter equation.
\end{rema}
\begin{prop}
Let $({\frak J},B)$ be a pseudo-euclidean Jordan algebra and
$r=\sum_{i=1}^na_i\otimes b_i$ be an antisymmetric $r-$matrix.
Let ${\cal U}:{\frak J}\longrightarrow{\frak J}$ be the linear map
defined by ${\cal U}=R\circ\phi$. Then, $Im({\cal U})=\{{\cal
U}(x);\,\,\,x\in{\frak J}\}$ is a Jordan subalgebra of $\frak J$.
Further, the bilinear form $\omega:Im({\cal U})\times Im({\cal
U})\longrightarrow\mathbb K$ defined by: $$\omega({\cal
U}(x),{\cal U}(y))=B({\cal U}(x),y),\,\,\,\forall x,y\in{\frak
J},$$ is a symplectic form on $Im({\cal U})$.
\end{prop}
\dem Let $x\in ker{\cal U}$ and $y\in Im({\cal U})$. There exists
$z$ in ${\frak J}$ such that $y={\cal U}(z)$, then  $B(x,y)=B(x,{\cal U}(z))=-B({\cal U}(x),z)=0.$ Thus, $Im({\cal U})\subset(ker {\cal U})^\perp$. Since
$dim\Bigl( Im({\cal U})\Bigl)= dim\Bigl(( ker{\cal U})^\perp
\Bigl)$, then $Im({\cal U})\subset(ker {\cal U})^\perp.$

  Let $x={\cal U}(x_0),\,\,y={\cal U}(y_0)\in Im({\cal U})$ where $x_0,y_0\in{\frak J}$ and $z\in ker{\cal U}$.
$$B(xy,z)=B({\cal U}(x_0){\cal U}(y_0),z)=-B({\cal U}(z){\cal U}(x_0),y_0)-B({\cal U}(y_0){\cal U}(z),x_0)=0.$$
Thus, $Im({\cal U})$ is a subalgebra of ${\frak J}$. Now, let $x_0,y_0,z_0\in{\frak J}$ and $x={\cal
U}(x_0),\,y={\cal U}(y_0) $ et $ z={\cal U}(z_0)\in Im({\cal U})$.
\begin{eqnarray*}
\omega(x,y)=\omega({\cal U}(x_0),{\cal U}(y_0))=B({\cal
U}(x_0),y_0)=-B(x_0,{\cal U}(y_0))&&\\=-\omega({\cal U}(y_0),{\cal
U}(x_0))=-\omega(y,x).&&
\end{eqnarray*}
Thus $\omega$ is antisymmetric. Further,
$$\sum_{cyc}\omega(xy,z)=\sum_{cyc}\omega({\cal U}(x\star y),{\cal U}(z_0))=
\sum_{cyc}B({\cal U}(x_0){\cal U}(y_0),z_0)=0,$$ because $C_{\frak J}(r)=0$. Now, let $x\in \frak J$. If $\omega({\cal U}(x),{\cal U}(y))=0$ for all $y$ in ${\frak J}$, then $B({\cal U}(x),y)=0$, so  ${\cal U}(x)=0$. Consequently, $\omega_{Im({\cal U}\times Im({\cal U}} $ is nondegenerate.
 It follows that, $(Im({\cal U}),\omega_{Im({\cal U}\times Im({\cal U}})$ is a symplectic Jordan algebra.$\Box$

 \begin{coro}
 Let $({\frak J},B)$ be a pseudo-euclidean Jordan algebra and $r=\sum_{i=1}^na_i\otimes b_i\in{\frak J}\otimes{\frak J}$ such that $\tau(r)=-r$, $C_{\frak J}(r)=0$ and $r$ nondegenerate. Then ${\frak J}$ endowed with the bilinear form $\omega_{\cal U}:{\frak J}\times{\frak J}\longrightarrow\mathbb K$ defined by $$\omega_{\cal U}(x,y):=B({\cal U}^{-1}(x),y),\,\,\,\forall x,y\in {\frak J},$$ is a symplectic Jordan algebra.
\end{coro}
\dem The fact that $r$ is nondegenerate implies that ${\cal U}$ is
invertible. It follows that, $Im({\cal U})={\frak J}$. The last
proposition gives the result.  $\Box$
\begin{rema}
Let ${\cal U}$ be an invertible endomorphism of ${\frak J}$. Then,
${\cal U}$  satisfies:
$${\cal U}\Bigl({\cal U}(x)y+x\,{\cal U}(y)\Bigl)={\cal U}(x){\cal U}(y),\,\,\,\forall x,y\in{\frak J}.$$
if and only if $ {\cal U}^{-1}$ is a derivation of ${\frak
J}$.
\end{rema}
\begin{prop}
Let $({\frak J},B)$ be a pseudo-euclidean Jordan algebra. ${\frak
J}$ have a symplectic form $\omega,$ if and only if there exists
a $B$-antisymmetric invertible derivation $D$ of ${\frak J}$ such
that $\omega(x,y)=B(D(x),y),\,\,\,\forall x,y\in {\frak J}$.
\end{prop}
\dem Since $B$ and $\omega$ are two nondegenerate bilinear form of
${\frak J},$ then there exists an invertible endomorphism $\delta$
of ${\frak J}$ such that $\omega(x,y)=B(\delta(x),y),\,\,\,\forall
x,y\in {\frak J}.$ Further, since $\omega$ is symplectic, then
$$B(\delta(xy),z)=B(\delta(y),xz)+B_1(\delta(x),yz)
=B(\delta(y)x,z)+B(\delta(x)y,z),\,\,\,\forall x,y,z\in {\frak J}.$$
The fact that $B$ is nondegenerate implies that
$\delta(xy)=\delta(x)y+x\delta(y),\,\,\,\forall x,y,z\in {\frak J}.$
Hence $\delta$ is an invertible derivation of ${\frak J}.$
Conversely, if $\delta$ is a $B-$antisymmetric invertible derivation
 of ${\frak J},$ then it is clear that the bilinear form $\omega:{\frak J}\times{\frak
J}\longrightarrow\mathbb K$ defined by: $\omega(x,y)=
B(\delta(x),y),\,\,\,\forall x,y\in {\frak J}$ is a symplectic form of ${\frak J}$.  $\Box$

\begin{coro}\label{u}
Let $({\frak J},B,\omega)$ be a symplectic pseudo-euclidean Jordan
algebra. Then, there exists a $B-$antisymmetric endomorphism
${\cal U}:{\frak J}\longrightarrow{\frak J}$ satisfying ${\cal
U}({\cal U}(x)y+x\,{\cal U}(y))={\cal U}(x){\cal
U}(y),\,\,\,\forall x,y\in{\frak J}$ and such that
$\omega=\omega_{\cal U}$.
\end{coro}
\dem  Consider the $B-$antisymmetric invertible derivation $D$ of
${\frak J}$ defined by $\omega(x,y)=B(D(x),y),\,\,\,\forall
x,y\in{\frak J}$. Then, $ {\cal U}=D^{-1}$ is $B-$ antisymmetric
and satisfies $${\cal U}({\cal U}(x)y+x\,{\cal U}(y))={\cal
U}(x){\cal U}(y),\,\,\,\forall x,y\in{\frak J}.\Box$$

 The following Theorem gives a generalization of the Corollary $\ref{u}$.
\begin{theor}
Let $({\frak J},\omega)$ be a symplectic Jordan algebra. Then,
there exists a nondegenerate $r-$matrix which satisfies
$\tau(r)=-r$, $C_{\frak J}(r)=0$ and such that
$\omega=\omega_{\cal U}$ where ${\cal U}$ is the endomorphism
associate to $r$.
\end{theor}
\dem Let $\{e_1,...,e_n\}$ be a base of ${\frak J}$. One poses
$\alpha_{ij}=\omega(a_i,a_j)$. Since the form $\omega$ is
nondegenerate, then the matrix $(\alpha_{ij})_{1\leq i,j\leq n}$
is invertible. Let $(\beta_{ij})_{1\leq i,j\leq n}$ be the inverse
of the matrix $(\alpha_{ij})_{1\leq i,j\leq n}$ and
$r=\sum_{ij}\beta_{ij}(a_i\otimes a_j)$. It is clear that
$(\beta_{ij})_{1\leq i,j\leq n}$ is an antisymmetric matrix. Thus
$\tau(r)=-r$. Now we shall prove that $C_{\frak J}(r)=0$. Let
$\phi:{\frak J}\longrightarrow{\frak J}^*$ be the linear
isomorphism given by $\phi(x)=\omega(x,.)$ and $R=\phi^{-1}$.
Then, $R$ satisfies the Yang Baxter equation. In fact, Let
$f,h,l\in {\frak J}^*$. Then, there exists $x,y,z\in{\frak J}$
such that $f=\phi(x),\,\,h=\phi(y)$ and $l=\phi(z)$. Thus
$x=R(f),\,\,y=R(h)$ and $z=R(l).$ It is easy to see that
$$
<f,R(h)R(l)>+<h,R(l)R(f)>+<l,R(f)R(h)>=\omega(x,yz)+\omega(y,xz)+\omega(z,xy)=0.
$$
Hence, $R $ satisfies the $YBE$. It follows that, $C_{\frak
J}(r)=0.$ $\Box$\\\\

Recall that to
any  unital pseudo-euclidean Jordan algebra one can apply the
Tits Kantor Koecher   construction $(TKK$ construction) to obtain  a
quadratic Lie algebra $(\cite{koe1},\,\cite{koe2})$.  In the following, we are going to  make a slight modification to this
construction in order to obtain  a quadratic Lie
algebras (${\frak L}ie({\frak J}),B_{\cal L})$ starting from
the pseudo-euclidean Jordan algebras $({\frak J},B)$ which are not
necessarily unital.
 
Let $(\frak J,B)$ be a pseudo-euclidean Jordan algebra which is not
necessarily unital. $\overline{\frak J}$ is a copy of $\frak J$. One pose $L({\frak J}):=vect\{L_x,\,\,\,x\in{\frak J}\}$ and ${\cal H}({\frak J}):=L({\frak J}^2)
\oplus\left[ L({\frak J}),L({\frak J})\right] .$
On ${\cal H}({\frak J})$ we define the bilinear form $\Gamma$ by:
$$\Gamma( R_a+D_1, R_b+D_2)=B(a,b)+\Omega(D_1,D_2),\,\,\,
\forall a,b\in {\frak J}^2,\,\,D_1,D_2\in\left[ L({\frak J}),L({\frak J})\right] .$$
\begin{eqnarray*}
\mbox{where } \Omega \mbox{ is given by: }\Omega:\left[ L({\frak J}),L({\frak J})\right] \times\left[ L({\frak J}),L({\frak J})\right]&\longrightarrow&\mathbb K\\
(D_1 ,D_2=\sum_{i=1}^n\left[  R_{c_i}, R_{d_i}\right])&\longmapsto&\sum_{i=1}^nB(D_1(c_i),d_i).
\end{eqnarray*}

The bilinear form $\Gamma$ is well defined because the orthogonal of $\Bigl(Ann({\frak J})\Bigl)$ with respect $B$ is  ${\frak J}^2$.

The vector space ${\frak L}ie({\frak J})={\frak J}\oplus{\cal H}({\frak J})\oplus\overline{{\frak J}}$ with   the following bracket:

$$ \left[ T_1,T_2\right]=\left[ T_1,T_2\right]_{\cal H}, 
\left[T_1,a_2\right]=T_1(a_2), 
\left[T_1,\overline{b_2}\right]=\overline{-{T_1}(b_2)},$$
$$\left[a_1,\overline{b_2}\right]= R_{a_1{b_2}}, 
\left[a_1,a_2\right]=\left[\overline{b_1},\overline{b_2}\right]=0,\,\,\forall T_i\in {\cal H}({\frak J}),\,a_i,b_i\in{\frak J},\,i\in\{1,2\},$$

is a Lie algebra. Moreover $({\cal L}ie({\frak J}),B_{\cal L})$ is a quadratic Lie algebra, where the invariant sacalar product  $B_{\frak L}$ is defined by:
\begin{eqnarray*}
& B_{\cal L}:{\cal L}ie({\frak J})={\frak J}\oplus {\cal H}({\frak J})\oplus\overline{{\frak J}}\times{\cal L}ie({\frak J})={\frak J}\oplus {\cal H}({\frak J})\oplus\overline{{\frak J}}\longrightarrow\mathbb K\hskip6.5cm&\\
&(x_1=a_1+T_1+\overline{b_1},x_2=a_2+T_2+\overline{b_2})
\longmapsto\Gamma(T_1,T_2)+2B(a_1,b_2)+2B(a_2,b_1).&
\end{eqnarray*}

Before starting the description of symplectic pseudo-euclidean
Jordan algebras, it is natural to ask the following question: If
$({\frak J},B)$ is a pseudo-euclidean Jordan algebra and $({\frak
L}ie({\frak J}),B_{\cal L}) $ the Lie algebra constructed by the
$TKK$ construction from ${\frak J}$. Suppose that $({\frak J},B)$
has a symplectic structure. Then, does $({\frak L}ie({\frak
J}),B_{\cal L}) $ admit a symplectic structure?\\

We give the answer of this question in the following proposition.

  \begin{prop}
Let $({\frak J},B,\omega)$ be a symplectic pseudo-euclidean
Jordan algebra. Let  $D$ be the derivation of ${\frak J}$
defined by: $\omega(x,y)=B(D(x),y),\,\,\,\forall x,y\in {\frak
J}$. Define the linear map $D_{\cal L}:{\frak L}ie({\frak
J})\longrightarrow{\frak L}ie({\frak J})$ by
$$D_{\cal L}(a)=D(a),\,\,\,D_{\cal L}(\overline{a})=\overline{D(a)},\,\,D_{\cal L}(R_a)=R_{D(a)},
\mbox{ and } D_{\cal L}(\left[ R_a,R_b\right] )=\left[
R_{D(a)},b\right] +\left[ R_a,R_{D(b)}\right],$$ $\,\,\,\forall
a,b\in{\frak J}.$  Then, the bilinear form $\omega_{\cal L}$
defined on ${\frak L}ie({\frak J})$, by: $\omega_{\cal
L}(x,y)=B_{\cal L}(D_{\cal L}(x),y),\,\,\,\forall x,y\in {\frak
J}$, is a symplectic form if and only if $D_{\cal L}$ satisfies
the following condition:
\begin{eqnarray}
&& \hskip0.5cmD_{\cal L}(\left[ R({\frak J}),R({\frak J})\right])
=\left[ R({\frak J}),R({\frak J})\right].\hskip1cm\label{d1}
\end{eqnarray}
\end{prop}
\dem Since $D$ is  a derivation of ${\frak J},$ then 

$$D((x,y,z))= (D(x),y,z) + (x,D(y),z) + (x,y,D(z)),\,\,\, \forall x,y,z \in{\frak J}.$$
Consequently $D_{\cal L}$ is a derivation of ${\frak L}ie({\frak J}).$ Moreover the fact that $D$ is $B-$antisymmetric implies that $D_{\cal L}$ is $B_{\cal L}-$antisymmetric. Next, it is clear that  $D_{\cal L}$ is invertible if and only if $D_{\cal L}(\left[ R({\frak J}),R({\frak J})\right])=\left[ R({\frak J}),R({\frak J})\right].$ $\Box$

\begin{defi}
Let $({\cal G}={\cal G}_{\bar 0}\oplus{\cal G}_{\bar 1},B)$ be a quadratic $\mathbb Z_2-$graded Lie algebra. We say that
$({\cal G}={\cal G}_{\bar 0}\oplus{\cal G}_{\bar 1},B,\omega)$ is a symplectic quadratic $\mathbb Z_2-$graded Lie algebra, if $\omega:{\cal G}\times{\cal
G}\longrightarrow\mathbb K$ is a nondegenerate antisymmetric bilinear form which satisfies:
$$\omega(\left[ x,y\right] ,z)+\omega(\left[ y,z\right] ,x)+\omega(\left[ z,x\right] ,y),\,\,\,\forall x,y,z\in{\cal G}$$ and
such that $\omega({\cal G}_{\bar 0},{\cal G}_{\bar 1})=\{0\}$.
\end{defi}
 \begin{rema}\label{sympz2}
Let $({\frak J},B,\omega)$ be a symplectic pseudo-euclidean Jordan
algebra and $D$ be the antisymmetric invertible derivation of $
{\frak J}$ such that $\omega(x,y)=B(D(x),y),\,\,\,\forall
x,y\in{\frak J}$. If $D_{\frak L}$ satisfies the condition
$(\ref{d1})$, then $({\frak L}ie({\frak J})={\frak L}ie({\frak
J})_{\bar 0}\oplus{\frak L}ie({\frak J})_{\bar 1},B_{\frak
L},\omega_{\frak L}), $ where $\omega_{\frak L}(x,y)=B_{\frak
L}(D_{\frak L}(x),y),\,\,\,\forall x,y\in{\frak L}ie({\frak J}),$
is a symplectic quadratic $\mathbb Z_2-$graded Lie algebra.
\end{rema}
We are going to  give an example of nonassociative Jordan algebra which admits a derivation satisfying the condition $(\ref{d1})$.\\

\textbf{Example.} \,\,\, Let $A={\mathbb K}[X]$ $({\mathbb K}=
{\mathbb R}\mbox{ or }{\mathbb C})$ and ${\cal O}=XA$ the ideal of
$A$ generated by $\{X\}$. Let us consider ${\tilde{\cal O}}:={\cal
O}/X^n{\cal O},\,\,\, n\in {\mathbb N}^*$. ${\tilde{\cal O}}$ is
an associative algebra generated by $\{\overline
X,\overline{X^2},...,\overline{X^n}\}$. Let $\frak J$ be an $r$
dimension Jordan algebra $(r\in\mathbb N^*)$. The vector space
$\tilde{\frak J}={\frak J}\otimes{\tilde{\cal O}}$, endowed with
the following commutatif product:
$$(x\otimes \overline{P})(y\otimes \overline{Q})=(xy)\otimes(\overline{PQ}),\,\,\,\forall x,y\in\frak J
,\,\,\forall P,Q\in \cal O,$$ is a Jordan algebra. Now, let
${\tilde{\frak{J}}}\oplus(\tilde{\frak{J}})^*$ be the trivial
T$^*-$extension of ${\tilde{\frak{J}}}$ and  $I:=\{1,...,r\}\times\{1,...,n\}.$ Let
$(e_i)_{1\leq i\leq r}$ be a base of ${\frak{J}}$. Then,
${\Bigl(e_{ij}:=e_i\otimes \overline{X^j}\Bigl)}_{(i,j)\in I}$ is
a base of $\tilde {\frak J}$. Let $\Bigl(e_{ij}^*\Bigl)_{(i,j)\in
I}$ be the dual base of ${\frak{J}}$ associate to
${\Bigl(e_{ij}\Bigl)}_{(i,j)\in I}.$ For all $j\in\{1,...,n\},$
denote by ${\frak{J}}_{j}={\frak{J}}\otimes \overline{X^j}$. It is
easy to see that if $a\in{\frak{J}}_{i}$ and $b\in{\frak{J}}_{j}$,
then $ab\in{\frak{J}}_{i+j}$ if $i+j\leq n$ and $ab=0$ if $i+j>
n$. Let us consider  the endomorphism $D$ of $\tilde{\frak J}$
defined by: $D(x\otimes \overline{X^i})=ix\otimes
\overline{X^i},\,\,\,\forall x\in{\frak J},\,\forall
i\in\{1,...,n\}.$ Then, $D$ is a derivation of $\tilde{\frak J}$.
It is easy to verify that the map $\bar
D:{\tilde{\frak{J}}}\oplus(\tilde{\frak{J}})^*\longrightarrow{\tilde{\frak{J}}}\oplus(\tilde{\frak{J}})^*$
defined by:
$$\bar{D}(a+f)=D(a)-f\circ D,\,\,\,\forall a\in\tilde {\frak
J},\,f\in\tilde {\frak J}^*$$ is a derivation of
${\tilde{\frak{J}}}\oplus(\tilde{\frak{J}})^*$.
 On the other hand, let $i,k,t\in\{1,...,r\}$ and $j,l\in \{1,...,n\}.$ We have: $$\left[ R_{e_{ij}},R_{e^*_{kj}}\right](e_{tl})=e^*_{kj}\circ\Bigl(R_{e_{tl}}R_{e_{ij}}-R_{e_{ij}e_{tl}}\Bigl).$$
Since
$$e^*_{kj}\circ\Bigl(R_{e_{tl}}R_{e_{ij}}-R_{e_{ij}e_{tl}}\Bigl)(e_{pq})
=e^*_{kj}\Bigl(e_{tl}(e_{ij}e_{pq})\Bigl)-e^*_{kj}\Bigl((e_{tl}e_{ij})e_{pq}\Bigl)=0$$
because
\begin{eqnarray*}
e_{tl}(e_{ij}e_{pq}),\,\,( e_{tl}e_{ij})e_{pq}\in{\frak J}_{j+l+q} \mbox{ if } j+l+q\leq n\\
\mbox{ and } (e_{tl}(e_{ij}e_{pq}=(e_{tl}e_{ij})e_{pq}=0 \mbox{ if
}j+l+q> n.
\end{eqnarray*}
 Hence, if we consider $I_1:=\{(i,j,k,l),(i,j),(k,l)\in I,\,\,\,where\,\,\,j\neq l\}$, then
$$\left[ R(\bar{\frak{J}},R(\bar{\frak{J}}\right] =Vect\Biggl(\bigcup_{(i,j),(k,l)\in I}\Bigl(\Bigl\{\left[
R_{e_{ij}},R_{e_{kl}}\right],\,\,\left[
R_{e^*_{ij}},R_{e^*_{kl}}\right]\Bigl\}\Bigl)\Biggl)\bigcup\Biggl
(\bigcup_{(i,j,k,l)\in I_1}\Bigl\{\left[
R_{e_{ij}},R_{e^*_{kl}}\right]\Bigl\}\Biggl).$$
 By the definition of $\bar{D}$, we have:
$\bar{D}(e_{ij})=ie_{ij}$ and
$\bar{D}(e^*_{ij})(e_{kl})=-e^*_{ij}(D(e_{kl}))=-le^*_{ij}(e_{kl})$.
Hence $\bar{D}(e^*_{ij})(e_{kl})=0$ if $(i,j)\neq (k,l)$ and
$\bar{D}(e^*_{ij})(e_{ij})=-j$. It follows that,
$\bar{D}(e^*_{ij})=-je^*_{ij}$.
$$\left[ \bar{D},\left[ R_{e_{ij}},R_{e_{kl}}\right] \right] =(j+l)\left[ R_{e_{ij}},R_{e_{kl}}\right].$$
$$\left[ \bar{D},\left[ R_{e^*_{ij}},R_{e^*_{kl}}\right] \right] =-(j+l)\left[ R_{e^*_{ij}},R_{e^*_{kl}}\right].$$
It follows that, $$\left[ \bar{D},\left[
R(\bar{\frak{J}},R(\bar{\frak{J}}\right]\right] =\left[
R(\bar{\frak{J}},R(\bar{\frak{J}}\right].$$ Hence, $ \bar{D}$
satisfies the condition $(\ref{d1})$.
\section{Descriptions of symplectic pseudo-euclidean Jordan algebras}\label{dajps}

In this section, $\mathbb K$ is    algebraically closed.

\begin{theor}\label{syor}
Let $({\frak J}_1,B_1,\omega_1)$ be a symplectic pseudo-euclidean
Jordan algebras, $a_0\in{\frak J}_1$, $\lambda\in\mathbb K$ and
$\delta$ the invertible derivation of ${\frak J}_1$ such that
$$\omega_1(x,y)=B_1(\delta(x),y),\,\,\,\forall x,y\in {\frak J}_1.$$
Let $(D,x_0)\in\mbox{End}_s({\frak J}_1)\times{\frak J}_1$ be an admissible pair such that
\begin{eqnarray*}
D(a_0)=\lambda x_0+{{1}\over{2}}\delta(x_0)\hskip 1.5cm
R_{a_0}=\delta   D-D\delta+\lambda D,\,\,\,\forall x\in{\frak J}_1.
\end{eqnarray*}
Let $({\frak J}=\mathbb Ka\oplus{\frak J}_1\oplus\mathbb Kb,B)$ be
the generalized double extension  of $({\frak J}_1,B_1)$ by the
one dimensional Jordan algebra with zero product by means of
$(D,x_0)$. Then, the endomorphism $\Delta$ of ${\frak J}$ defined
by:
$$\Delta(b)=\lambda b,\hskip 0.7cm \Delta(x)=\delta(x)-B_1(a_0,x)b,\hskip
0.7cm\Delta(a)=a_0-\lambda a,\,\,\,\forall x\in{\frak J}_1,$$ is a
$B-$antisymmetric invertible derivation of ${\frak J}$. Thus the
antisymmetric bilinear form $\omega$ of ${\frak J}$ defined by:
$\omega(x,y)=B(\Delta(x),y),\,\,\,\forall x,y\in {\frak J},$ is a symplectic form on ${\frak J}$.\\
The algebra $({\frak J},B,\omega)$ obtained above is called the symplectic pseudo-euclidean
double extension of $({\frak J}_1,B_1,\omega_1)$ (by means of $(D,x_0,a_0,\lambda)\in\mbox{End}({\frak
J}_1)\times{\frak J}_1\times{\frak J}_1\times\mathbb K).$
\end{theor}
\dem
 The product on $\frak J$ is denoted by $\star$. Let $x,y\in{\frak J}_1,$ then, we have
$$\Delta(x\star a)-\Delta(x)\star a-x \star\Delta(a)=\delta D(x)-D\delta(x)-a_0x+\lambda D(x)+B_1(2\lambda x_0-a_0+\delta(x_0)-D(a_0),x)b=0.$$
$$
\Delta(x\star y)-\Delta(x)\star
y-x\star\Delta(y)=\delta(xy)-\delta(x)y-x\delta(y)
+B_1\Bigl((\lambda D-R_{a_0}-D\delta+\delta D)(x),y\Bigl)b=0.$$
$\Delta (a \star a)-2a\star D(a)=\delta(x_0)-2D(a_0)+2\lambda x_0+3(k\lambda-B_1(a_0,x_0))b=0.$\\

Hence, $\Delta$ is a derivation of ${\frak J}$. Moreover, since
$\delta$ is invertible, then $\Delta$ is invertible too. Thus, the bilinear
form $\omega$ defined on ${\frak J}$ by:
$\omega(x,y)=B(\Delta(x),y),\,\,\,\forall x,y\in {\frak J},$ is a
symplectic form on ${\frak J}$. $\Box$

\begin{lema}\label{nil1}
Let ${\frak J}$ be a Jordan algebra. If ${\frak J}$ admits  an
invertible derivation, then ${\frak J}$ is nilpotent.
\end{lema}
\dem Let $ Ra({\frak J})$ be the nilpotent radical of ${\frak J}$ and let 
$\cal S$ be a semi-simple subalgebra of ${\frak J}$ such that
${\frak J}=Ra({\frak J})\oplus\cal S$. The map
\begin{eqnarray*}
\delta&:&{\cal S}\longrightarrow{\cal S}\\
&&x\longmapsto\delta(x):=p_{s}\circ D(x),
\end{eqnarray*}
where $p_{s}:{\frak J}\longrightarrow {\cal S}$ is defined by $
p_{s}(s+r)=s,\,\,\,\forall (s,r)\in {\cal S}\times Ra({\frak J})$,
is a derivation of $\cal S$. In fact, Let $x,y\in{\cal S}.$ One
poses $D(x)=x_1+x_2,$ $D(y)=y_1+y_2$ where $x_1,y_1\in {\cal S}$
and $x_2,y_2\in  Ra({\frak J})$.
\begin{eqnarray*}
\delta(xy)=p_s(D(xy))=p_s(xD(y)+D(x)y)=xy_2+x_1y=x\delta(y)+\delta(x)y.
\end{eqnarray*}
Thus $\delta$ is a derivation of $\cal S$. Now, let $x\in{\cal S}$
such that  $\delta(x)=0$. Thus, $D(x)=r\in Ra({\frak J}).$ Hence,
$D^{-1}(r)=x\in\cal S$, which is impossible because $D^{-1}(
Ra({\frak J}))= Ra({\frak J})$. It follows that, $\delta$ is
injective. On the other hand, since $\cal S$ is semi-simple, then
it admits a unital element $e$. Thus $\delta(e)=0$. It follows
that $e=0$. Which show that ${\cal S}=\{0\}$. Hence ${\frak J}=
Ra({\frak J})$. Thus ${\frak J}$ is nilpotent. $\Box$

\begin{theor}\label{dbl}
Let $({\frak J},B,\omega)$ be a symplectic pseudo-euclidean Jordan
algebra such that ${\frak J}\neq\{0\}$. Then $({\frak
J},B,\omega)$ is a symplectic pseudo-euclidean double extension of
a symplectic pseudo-euclidean Jordan algebra $({\cal W},T,\Omega).$
\end{theor}
\dem Let $\Delta$ be the invertible $B-$antisymmetric derivation
of ${\frak J}$ defined by, $\omega(x,y)=B(\Delta(x),y) $ for all
$x,y\in{\frak J}$. Since ${\frak J}$ is nilpotent (see the Lemma
$\ref{nil1}$), then $Ann({\frak J})\neq\{0\}.$ Let $x\in
Ann({\frak J})\backslash \{0\}$, then
$$\Delta(x)y=\Delta(xy)-x\Delta(y)=0,\,\,\,\,\forall y\in {\frak J}.$$ Hence  $\Delta(x)\in Ann({\frak J})$.
It follows that, $\Delta({Ann({\frak J})})=Ann({\frak J})$. Let
$\lambda$ be an eigenvalue of $\Delta_{\vert_{Ann({\frak J})}}.$
Let  $b \in {Ann({\frak J})}$  such that $\Delta(b)= \lambda b$. Then 
$0=\omega(b,b)=\lambda B(b,b)$. Which implies that $B(b,b)=0$. One
pose ${\cal I}=\mathbb K b$ and ${\cal I}^\perp$ the
$B-$orthogonal of $\cal I$. Let $a\in {\frak J}$ such that
$B(a,a)=0$, $B(a,b)=1$ and ${\frak J}= \mathbb Ka\oplus{\cal
I}^\perp$. By the Theorem $\ref{th3}, $ the pseudo-euclidean
Jordan algebra $({\frak J},B)$ is a generalized  double extension
of the pseudo-euclidean Jordan algebra (${\cal W}=(\mathbb
Ka\oplus\mathbb Kb)^\perp,T:=B_{\vert_{{\cal W}\times{\cal W}}})$
by the one dimensional Jordan algebra with zero product by means of
the pair $(D,x_0)\in End_s({\cal W})\times{\cal W}$ defined by:
$D=p\circ {R_{a}}_{\vert_{\cal W}}$ and $x_0=p\circ R_a(a)$, where
$p:{\cal W}\oplus\mathbb Kb \longrightarrow{\cal W}$  is the
projection defined by $p(x+\alpha b)=x,\,\,\,\forall x\in{\cal
W},\,\alpha\in\mathbb K$. The product in ${\frak J}=\mathbb
Ka\oplus{\cal W}\oplus\mathbb Kb$ is denoted by $\star$ and given
by:
$$\left\{\begin{array}{l}
b\star{\frak J}=\{0\},\\
a\star a=w_0=x_0+kb,\\
x\star y:=xy+B(D(x),y)b,\hskip 1cm\forall x,y\in {\cal W}.\\
a\star x=D(x)+B(x_0,x)b\\
 \end{array}\right.$$
where $xy$ is the multiplication of $x$ by $y$ in ${\cal W}$. Let
$x\in{\cal W}$. We have $B(\Delta(x),b)=-B(x,\Delta(b))=-\lambda
B(x,b)=0.$ Thus $\Delta({\cal W})\subset {\cal I}^\perp={\cal
W}\oplus\mathbb Kb.$ We obtain:
$$
\Delta(b)=\lambda b, \hskip 0.5cm \Delta(x)=\delta(x)+\psi(x)b,
\,\,\,\forall x\in{\cal W},\hskip 0.5cm \Delta(a)=\alpha
a+a_0+\beta b,
$$
where $\delta$ is an endomorphism of ${\cal W}$, $\psi$ is a
linear form of ${\cal W}$, $\alpha,\beta\in\mathbb K$ and $a_0\in
{\cal W}$. The fact that $\Delta$ is $B-$antisymmetric and
$B(a,a)=0$ implies, that $\delta$ is $T-$antisymmetric, $\beta=0$,
$\alpha=-\lambda$ and $\psi=-B(a_0,.)$. Thus,
$$\Delta(x)=\delta(x)-B(a_0,x),\,\,\,\forall x\in{\cal W}\hskip 0.5cm \Delta(a)=\alpha a+a_0.$$
Further, the fact that $\Delta$ is an invertible derivation of
${\frak J}$, implies that $\delta$ is an invertible derivation of
${\cal W}$ and
$$
\delta(x_0)=2D(a_0)-2\lambda x_0,\hskip 0.5cm
k\lambda=B(a_0,x_0),\hskip 0.5cm {R_{a_0}}_{\vert_{{\cal
W}}}=\delta D-D\delta+\lambda D.
$$
Since $\delta$ is an invertible $T-$antisymmetric derivation of
${\cal W}$, then the bilinear form $\Omega:{\cal W}\times{\cal
W}\longrightarrow\mathbb K$ defined by:
$\Omega(x,y):=T(\delta(x),y),\,\,\,\forall x,y\in{\cal W}$ is a
symplectic form on $ {\cal W}$. $\Box$
\begin{coro}
Let $({\frak J},B,\omega)$ be a symplectic pseudo-euclidean Jordan
algebra. Then $({\frak J},B,\omega)$ is obtained from the
algebras $\{0\}$ by a finite sequence of symplectic pseudo-euclidean double extensions of
symplectic pseudo-euclidean Jordan algebras.

\end{coro}

\section{ Jordan-Manin Algebras}\label{mani}
\subsection{Double extensions of Jordan-Manin algebras}
\begin{defi}
A pseudo-euclidean Jordan algebra $({\frak J},B)$ is said
Jordan-Manin algebra if there exists two totaly isotropic
subalgebras ${\cal U}$, $\cal V$ of ${\frak J}$ such that ${\frak
J}={\cal U}\oplus{\cal V}.$ In this case the triple $({\frak
J},{\cal U},{\cal V})$ is said a Jordan-Manin triple.
\end{defi}
\begin{theor}\label{thman}
Let $({\frak J}={\cal U}\oplus{\cal V},B)$ be a Jordan-Manin
algebra. Let $(D,x_0)\in\mbox{End}_s({\frak J})\times{\frak J}$ be
an admissible pair of ${\frak J}$ such that $D({\cal
V})\subset{\cal V}$ and $x_0\in{\cal V}$. The Jordan algebra
obtained by generalized double extension of ${\frak J}$ by means
of $(D,x_0,0)$ is a Jordan-Manin algebra $(\tilde{\frak J}={\cal
U}'\oplus{\cal V}',\tilde B),$ where ${\cal U}'={\cal
U}\oplus\mathbb Kb$ is the direct product of ${\cal U}$ and
$\mathbb Kb$ and ${\cal V}'={\cal V}\oplus\mathbb Ka$ is the
generalized semi-direct product of ${\cal V}$ by $ \mathbb Ka$ by
means of $(D_{\vert_{\cal V}},x_0)$.
\end{theor}
\dem The associatif scalar product on $\tilde{\frak J}$ is given
by
\begin{eqnarray*}
\tilde B_{\vert_{{{\frak J}\times{\frak J}}}}=B,\,\, \tilde
B(a,b)=1,\,\,\tilde B(a,{\frak J})=\tilde B(b,{\frak
J})=\{0\}\mbox{ et }\tilde B(a,a)=\tilde B(b,b)=0 .
\end{eqnarray*}
Let ${\cal U}'={\cal U}\oplus\mathbb Kb$ and ${\cal V}'= {\cal
V}\oplus\mathbb Ka$. It is easy to see that ${\cal U}'$ and ${\cal
V}'$ are two totaly isotropic subalgebra of $\tilde{\frak J}$.
Further, it is clear that $\tilde{\frak J}={\cal U}'\oplus{\cal
V}'$. Thus $\tilde{\frak J}$ is a Jordan-Manin algebra. $\Box$
\begin{defi}
The Jordan-Manin algebra $(\tilde{{\frak J}}={\cal U}'\oplus{\cal
V}',\tilde B)$ is said the double extension of the Jordan-Manin
algebra $({\frak J}={\cal U}\oplus{\cal V},B)$ by the one
dimensional algebra with zero product by means of $(D,x_0)$.
\end{defi}
\begin{theor}\label{man}
Let $({\frak J}={\cal U}\oplus{\cal V},B)$ be a Jordan-Manin
algebra. If ${\cal U}\cap Ann( {\frak J})\neq\{0\}$ or ${\cal
V}\cap Ann({\frak J})\neq\{0\}$, then ${\frak J}$ is the double
extension of the Jordan-Manin algebra $({\cal W}={\cal
U}'\oplus{\cal V}',T)$ by the one dimensional Jordan algebra to
null product.
\end{theor}
\dem Suppose that ${\cal U}\cap Ann( {\frak J})\neq\{0\}$. Let
$b\in{\cal U}\cap Ann( {\frak J})\setminus\{0\}$. There exists
$a\in {\cal V}$ such that $B(a,a)=0$, $B(a,b)=1 $ and ${\frak
J}=(\mathbb Kb)^\perp\oplus\mathbb Ka$. By the proof of the
Theorem $\ref{th3}$, ${\frak J}$ is the generalized double
extension of ${\cal W}=(\mathbb Ka\oplus\mathbb Kb)^\perp$ by the
one dimensional Jordan algebra with zero product by means of the
pair $(D,x_0)\in End_s({\cal W})\times{\cal W}$ defined by:
$D=p\circ {R_a}_{\vert_{{\cal W}}}$ and $x_0=p\circ R_a(a)$, where
$p:{\cal W}\oplus\mathbb Kb\longrightarrow{\cal W}$ is the
projection defined by $p(x+\alpha b)=x,\,\,\,\forall x\in{\cal
W},\,\alpha\in\mathbb K$. The product in ${\frak J}=\mathbb
Ka\oplus{\cal W}\oplus\mathbb Kb$ is given by:
$$\left\{\begin{array}{ll}
b\star \tilde{\frak{J}}=\{0\},\\
 a\star a=x_{0}+kb,\\
x\star y=\beta(x,y)+B(D(x),y)b,\hskip 0.4cm &\,\,\forall x,y \in{\cal W}.\\
a \star x=D(x)+B(x_{0},x)b,
\end{array}\right.$$
where $\beta$ is a Jordan product on ${\cal W}$ and $k\in\mathbb
K$. We have ${\cal I}=\mathbb Kb\subset {\cal U},$ thus ${\cal
U}^\perp\subset {\cal I}^\perp$. Since ${\cal U}={\cal U}^\perp$,
then ${\cal U}\subset{\cal I}^\perp$. It follows that, ${\cal
I}^\perp={\cal U}\oplus ({\cal V}\cap{\cal I}^\perp).$ Further,
${\cal V}\cap{\cal I}^\perp\subset {\cal W}$ because if $ x\in
{\cal V}\cap{\cal I}^\perp$, then $B(a,x)=B(b,x)=0$. Hence, ${\cal
W}=({\cal U}\cap{\cal W})\oplus({\cal V}\cap{\cal I}^\perp).$ One
poses ${\cal U}'={\cal U}\cap{\cal W}$ and ${\cal V}'={\cal
V}\cap{\cal I}^\perp$. It is clear that ${\cal V}'$ is an ideal of
${\cal V}$. Thus, for all $v\in{\cal V}'$, $a\star v\in{\cal V}'$.
Hence, $D(v)+B(v,x_0)b\in{\cal V}'\,\,\,\forall v\in {\cal V}'$.
It follows that, $B(v,x_0)=0$ and $D(v)\in{\cal V}',\,\,\,\forall
v\in{\cal V}'. $ Hence, $x_0\in({\cal V}')^\perp\cap{\cal W}={\cal
V}'$ and $D({\cal V}')\subset {\cal V}'$. Further, ${\cal U}'$ and
${\cal V}'$ are two subalgebra of ${\cal W}.$ In fact, let
$x,y\in{\cal U}'$, then $\alpha(x,y)=x\star y-B(D(x),y)b\in{\cal
U}$ because $b\in{\cal U}$. Since $\alpha(x,y)\in {\cal W}$, then
$\alpha(x,y)\in {\cal U}'.$ The fact that ${\cal V}'$ is stable
under $D$ and that ${\cal V}'\subset {\cal V}$ is totally
isotropic implies that ${\cal V}'$ is a subalgebra of ${\cal W}.$
Further , it is clear that ${\cal U}'$ and ${\cal V}'$ are totally
isotropic. It follows that, $(D,x_0)$ is an admissible pair of
${\cal W}={\cal U}'\oplus{\cal V}'$ such that $x_0\in {\cal V}'$
and $D({\cal V}')\subset {\cal V}'$. Hence ${\frak J}={\cal
U}\oplus{\cal V}$ is a double extension of the Jordan-Manin
algebra ${\cal W}={\cal U}'\oplus{\cal V}'$ by means of $(D,x_0)$.
For the other case ${\cal V}\cap Ann( {\frak J})\neq\{0\}$, the
proof is similar.$\Box$

\begin{lema}\label{inter}
Let $({\frak J}={\cal U}\oplus{\cal V},B)$ be a Jordan-Manin
algebra. If ${\frak J}$ is nilpotent, then ${\cal U}\cap Ann(
{\frak J})\neq\{0\}$ and ${\cal V}\cap Ann( {\frak J})\neq\{0\}$.
\end{lema}
\dem Suppose that ${\cal U}\cap Ann( {\frak J})=\{0\}$. One poses
$L_0={\cal V}Ann( {\cal U})=\{xy;\,\,x\in{\cal V}\mbox{ and } y\in
Ann( {\cal U})\}$. It is easy to check that $L_0\not= \{0\}$  and
include in ${\cal U}$. Consider ${\cal I}_0$ the ideal of ${\cal
U}$ generated by $L_0$. Then ${\cal I}_0$ is a non-zero ideal of
the nilpotent algebra ${\cal U}$. Thus, ${\cal I}_0\cap Ann({\cal
U})\neq\{0\}.$ Now, let $L_1={\cal V}\Bigl({\cal I}_0\cap
Ann({\cal U})\Bigl)$. It is clear that $L_1\not=\{0\}$   and
include in ${\cal U}$. Let ${\cal I}_1$ be the ideal of ${\cal U}$
generated by $L_1$. By the same argument, we show that ${\cal
I}_1\cap Ann({\cal U})\neq\{0\}$. We repeat this process and we
obtain the subsets  sequence $(L_{n}){_{n\in\mathbb N}}$ of ${\cal
U}$ and the ideals sequence $({\cal I}_{n}){_{n\in\mathbb N}}$ of
$\cal U$ defined by $L_0={\cal V}Ann( {\cal U})$, $L_n={\cal
V}({\cal I}_{n-1}\cap Ann({\cal U}))$, and ${\cal I}_n$ is an
ideal generated by $L_n\,\,\,\forall n\in\mathbb N$. It is easy to
see that for all $n\in\mathbb N,$ ${\cal I}_n$ is a non-zero ideal
of ${\cal U}$ and that for all
 $n\in\mathbb N$ ${\cal I}_n\subset C^{n+1}({{\frak J}}).$ which contradict the fact that ${\frak J}$ is nilpotent.
For the other case ${\cal V}\cap Ann( {\frak J})\neq\{0\}$ the
proof is similar.$\Box$

The following Corollary is a consequence of the Theorem $\ref{man}$ and the Lemma $\ref{inter}$.
\begin{coro}
Let $({\frak J}={\cal U}\oplus{\cal V},B)$ be a Jordan-Manin algebra. If ${\frak J}$ is nilpotent, then ${\frak J}$ is the double extension of a Jordan-Manin algebra $({\cal W},T)$
by the one dimensional algebra with zero product.
\end{coro}
\begin{coro}
Let $({\frak J}={\cal U}\oplus{\cal V},B)$ be a nilpotent
Jordan-Manin algebra. Then, ${\frak J}$ is obtained from the two
dimension Jordan-Manin algebra by a finite sequence of Manin
double extensions of nilpotent Jordan-Manin algebras.
\end{coro}
 \subsection{Symplectic Jordan-Manin algebras}\label{ajmsp}
\begin{defi}
Jordan-Manin $({\frak J}={\cal U}\oplus{\cal V},B)$ is said
symplectic Jordan-Manin algebra if there exists a symplectic
structure $\omega$ on ${\frak J}$ satisfying $\omega({\cal
U},{\cal U})=\omega({\cal V},{\cal V})=\{0\}.$
\end{defi}
The study of symplectic Jordan-Manin algebra is intersting because every pseudo-euclidean symplectic Jordan algebras is a symplectic Jordan-Manin algebra (see. Proposition $\ref{spe}$). This result, will be proved in the Proposition $\ref{spe}$. The proof of the Proposition $\ref{spe}$ is based on  the following Lemma (Proposition$2.$ $(ii).$ p. 8
of $\cite{bou})$.
\begin{lema}$\cite{bou})$\label{bou}
Let ${\cal E},\,\,{\cal E}',\,\,{\cal F}$ be a vector spaces and
$S$ be a set. Let $r:S\longrightarrow End({\cal
E}),\,\,r':S\longrightarrow End({\cal E}')$ and
$q:S\longrightarrow End({\cal F})$. For all $\lambda:
S\longrightarrow \mathbb K$, one pose:
$${\cal E}(\lambda)=\{x\in{\cal E};\,\,\forall s\in S,\,\,(r(s)-\lambda(s))^n(x)=0, \,\,\mbox{for } n\in\mathbb N^*\}.$$
Let $\sigma:{\cal E}\times{\cal E}'\longrightarrow{\cal F}$ be a bilinear map such that:
$$q(s)\sigma(x,x')=\sigma(r(s)x,x')+\sigma(x,r'(s)x'), \,\,\,\forall s\in S,\,\, x\in{\cal E},\,\,x'\in{\cal E}'. $$
Then, For all maps $\lambda,\,\mu: S\longrightarrow \mathbb K$, we
have  $\sigma\Bigl({\cal E}(\lambda),{\cal E}'(\mu)\Bigl)\subset
{\cal F}({\lambda+\mu}).$
\end{lema}
\begin{prop}\label{spe}
Let $({\frak J},B,\omega)$ be a symplectic pseudo-euclidean Jordan
algebra over an algebraically closed fields $\mathbb K$. Then,
there exists two subalgebras $\cal U$ and $\cal V$ of ${\frak J}$
such that $({\frak J}={\cal U}\oplus{\cal V},B,\omega)$ is a
symplectic Jordan-Manin algebra.
\end{prop}

\dem Let $D$ be the derivation defined by
$\omega(x,y)=B(D(x),y)$ for all $x,y\in{\frak J}$. Let $Sp(D)$
the set of the eigenvalues of $D$. One pose
\begin{eqnarray*}
Sp^+=\{\lambda\in
Sp(D);\,\,\,\mbox{Re}(\lambda)>0\}\cup\{\lambda\in
Sp(D);\,\,\,\mbox{Re}(\lambda)=0 \mbox{ and }
Im(\lambda)>0\}&&\\=(Sp(D)\cup i\mathbb R_+) \hskip 0.8cm\mbox{
and } \hskip 0.8cmSp^- =\{\lambda\in Sp(D);\,\,\,-\lambda\in
Sp^+\}.\hskip 0.8cm&&
\end{eqnarray*}
It is clear that $Sp(D)=Sp^+\cup Sp^-$ and $Sp^+\cap
Sp^-=\emptyset.$ Now, one poses $${\cal U}=\sum_{\lambda\in
Sp^+}{\frak J}(\lambda),\,\,\,\,\,\,{\cal V}=\sum_{\lambda\in
Sp^-}{\frak J}({\lambda}),$$ where ${\frak J}({\lambda})=\{x\in
{\frak J};\,\,(D-\lambda id_{\frak J})^{dim({\frak J})}(x)=0\}.$
${\cal U}$ and ${\cal V}$ are two subalgebras of ${\frak J}.$ In
fact, if we consider in the Proposition $\ref{bou}$ ${\cal
E}={\cal E}'={\cal F}=S={\frak J}$ and we define
\begin{eqnarray*}
\sigma: {\frak J}\times{\frak J}\longrightarrow{\frak J} ,\hskip
0.5 cm \hskip 0.5 cm r:{\frak
J}\longrightarrow\mbox{End}({\frak J})\hskip 0.5 cm and\hskip 0.5 cm r'=q=r.\\
(x,y)\longmapsto xy,\hskip 1.4 cmx\longmapsto D \hskip 3.9cm
\end{eqnarray*}
 Since $D$ is a
derivation, then
$$q(x)\sigma(y,z)=\sigma(r(x)y,z)+\sigma(y,r'(x)z),\,\,\,\forall
x,y,z\in{\frak J}.$$ Thus, by  the Proposition $\ref{bou},$ we
have $\sigma\Bigl({\frak J}({\lambda}),{\frak
J}({\mu})\Bigl)\subset {\frak J}({\lambda+\mu}),\,\,\,\forall
\lambda,\,\,\mu\in\mathbb K$. Thus ${\frak J}({\lambda}){\frak
J}({\mu})\subset  {\frak J}({\lambda+\mu}),\,\,\,\forall
\lambda,\,\,\mu\in\mathbb K$. Thus, ${\cal U}$ and ${\cal V}$ are
two subalgebras of ${\frak J}.$ Now, we use the same proposition,
to prove that ${\cal U}$ and ${\cal V}$ are totally isotropic. One
pose ${\cal E}={\cal E}'=S={\frak J},$  ${\cal F}=\mathbb K$ and
define
\begin{eqnarray*}
r=r':{\frak J}\longrightarrow End({\frak J}) \hskip 0.5 cm
et\hskip 0.5 cm q:{\frak J}\longrightarrow End(\mathbb
K)\\x\longmapsto D,\hskip 2.6 cm x\longmapsto 0,\hskip 0.9 cm
\end{eqnarray*}
Since $D$ is $B-$antisymmetric, then
$$q(x)B(y,z)=B(r(x)y,z)+B(y,r(x)z)=0,\,\,\,\forall x,y,z\in{\frak
J}.$$ It follows that, by the Proposition $\ref{bou}$, we have
$B({\frak J}(\lambda),{\frak J}(\mu))\subset \mathbb
K(\lambda+\mu).$
$$\mbox{ However, }\hskip 0.3 cm \mathbb K(\lambda):=\{\alpha\in\mathbb K;\,\,\,(q(x)-\lambda)^n(\alpha)=0;\,\,
\mbox{ where } n\in\mathbb N^*\}=\{\alpha\in\mathbb K
;\alpha\lambda=0\}.\hskip 3 cm$$ It follows that, for all
$\lambda\neq 0$, we have $\mathbb K(\lambda)=\{0\}.$ Consequently,
since the sum of two elements of $Sp^+$ (resp. $Sp^-$) is not
zero, then ${\cal U}$ and ${\cal V}$ are totally isotropic (i.e.
$B({\cal U},{\cal U})={\cal U}({\cal V},{\cal V})=0$). It follows
that, ${\cal U}$ and ${\cal V}$ are two $D-$stable totally
isotropic subalgebras of ${\frak J}$ wich satisfies ${\frak
J}={\cal U}\oplus{\cal V}$. Hence, $({\frak J}={\cal U}\oplus{\cal
U},B,\omega)$ is a symplectic Jordan-Manin algebra.
 $\Box$
\begin{theor}\label{232}
Let $({\frak J}_1={\cal U}_1\oplus{\cal V}_1,B_1,\omega_1)$ be a symplectic Jordan-Manin algebra. Let $a_0\in{\cal V}_1$,
$\lambda\in\mathbb K$ and $\delta$ be the invertible derivation of
${\frak J}_1$ defined by
$$\omega_1(x,y)=B_1(\delta(x),y),\,\,\,\forall x,y\in {\frak J}_1.$$
Let $(D,x_0)\in\mbox{End}_s({\frak J}_1)\times{\frak J}_1$ be an admissible pair such that:
\begin{eqnarray*}
D({\cal V})\subset {\cal V},\hskip1cm x_0\in{\cal V}\hskip1cm
D(a_0)=\lambda x_0+{{1}\over{2}}\delta(x_0)\hskip 1cm
R_{a_0}=\delta D-D\delta+\lambda D.
\end{eqnarray*}
Then, the symplectic pseudo-euclidean double extension
$({\frak J},B,\omega)$ of $({\frak J}_1={\cal U}_1\oplus{\cal
V}_1,B_1,\omega_1)$ by means of $(D,x_0,a_0,\lambda)$ is a symplectic Jordan-Manin algebra.
\end{theor}
\dem By the Theorem $\ref{syor}$, the Jordan algebra $({\frak
J},B,\omega)$ is symplectic and pseudo-euclidean, where
$\omega(x,y)=B(\Delta(x),y),\,\,\,\,\forall x,y\in{\frak J,}$ and
$\Delta$ is the $B-$antisymmetric invertible derivation of
$({\frak J},B)$ defined by:
$$\Delta(b)=\lambda b,\hskip 0.7cm \Delta(x)=\delta(x)-B_1(a_0,x)b,\hskip
0.7cm\Delta(a)=a_0-\lambda a,\,\,\,\forall x\in{\frak J}_1.$$
Further, by the Theorem $\ref{thman}$, The algebra $({\frak
J}={\cal U}\oplus{\cal V},B)$ where ${\cal U}={\cal
U}_1\oplus\mathbb Kb$, ${\cal V}={\cal V}_1\oplus\mathbb K a$ est
is a Jordan-Manin algebra. It remains to be shown that
$\Delta({\cal U})={\cal U}$ and $\Delta({\cal V})={\cal V}$. In
fact, let $u=u_1+\alpha b\in{\cal U}$.
$$\Delta(u)=\Delta(u_1)+\alpha\Delta(b)=\delta(u_1)+(\alpha\lambda b-B_1(a_0,u_1))b\in{\cal U}.$$
Let $v=v_1+\beta a\in{\cal V}$.
$$\Delta(v) =\Delta(v_1)+\beta\Delta(a)=\delta(v_1)-\beta a_0-\beta\lambda a\in{\cal V}.\Box$$
\begin{defi}
The algebra $({\frak J}={\cal U}\oplus{\cal V},B,\omega)$
constructed above is said the symplectic-Manin double extension of
the symplectic Jordan-Manin algebra $({\frak J}_1={\cal
U}_1\oplus{\cal V}_1,B_1,\omega_1)$ by means of
$(D,x_0,a_0,\lambda)$.
\end{defi}
\begin{theor}
 Let $({\frak J}={\cal U}\oplus{\cal V},B,\omega)$ be a symplectic Jordan-Manin algebra.
  Let $\Delta$ be the $B-$antisymmetric invertible derivation defined by $\omega(x,y)=B(\Delta(x),y),\,\,\,
  \forall x,y\in{\frak J}$. If $\Delta$ admits  an eigenvector $c\in Ann({\frak J})\cap{\cal U}+Ann({\frak J})
  \cap{\cal V}$, then $({\frak J}={\cal U}\oplus{\cal V},B,\omega)$ is the symplectic-Manin double
  extension of the symplectic Jordan-Manin algebra $({\cal W}={\cal U}'\oplus{\cal
V}',T,\Omega)$ by means of
  $(D,x_0,a_0,\lambda)\in End({\cal W})\times{\cal V}'\times{\cal V}'\times\mathbb K$ satisfying $D({\cal V}')
  \subset {\cal V}'$.
\end{theor}
\dem Let $c\in Ann({\frak J})\cap{\cal U}+Ann({\frak J})\cap{\cal
V}$ such that $c\neq 0$ and $\Delta(c)=\lambda c$ where
$\lambda\in\mathbb K$. One poses $c=b+b' $ where $b\in Ann({\frak
J})\cap{\cal U}$ and $b'\in Ann({\frak J})\cap{\cal V}$. Since $
\cal U$ and $\cal V$ are stable by $\Delta$, then
$\Delta(b)=\lambda b$ and $\Delta(b')=\lambda b'$. Since $c\neq
0$, then $b\neq 0$ or $b'\neq 0$. Suppose that $b\neq 0$ (if
$b'\neq 0$ the same demonstration is remade). Since $b \in
Ann({\frak J})\cap\,{\cal U}\setminus\{0\}$, then there exists
$a\in{\cal V}$ such that $B(a,a)=0,\, \,\,B(a,b)=1$ and ${\frak
J}=\mathbb Ka\oplus(\mathbb Kb)^\perp$. By the Theorem $\ref{dbl}$
$({\frak J},B,\omega)$ is the symplectic pseudo-euclidean double
extension of $({\cal W}=(\mathbb Ka\oplus\mathbb
Kb)^\perp,T,\delta)$ by means of $(D,x_0,a_0,\lambda)\in End({\cal
W})\times{\cal W}\times{\cal W}\times\mathbb K$, where
$a_0=p'\circ\Delta(a)$ and $(D,x_0)$  an admissible pair defined
by: $D=p\circ {R_a}_{\vert_{{\cal W}}}$ and $x_0=p\circ R_a(a)$,
where $p$ (resp. $p'$) is the projection of ${\cal W}\oplus\mathbb
Kb$ (resp. ${\cal W}\oplus\mathbb Ka)$ in ${\cal W}$ defined by
$p(x+\alpha b)=x$ (resp. $p'(x+\alpha a)=x),\,\,\,\forall
x\in{\cal W},\,\alpha\in\mathbb K$.

Recall that, $\delta=p\circ \Delta_{\vert_{{\cal W}}}$. On the
other hand, by the Theorem $\ref{man}$, $({\frak J}={\cal
U}\oplus{\cal V},B)$ is the Manin double extension of the
Jordan-Manin algebra $({\cal W}=(\mathbb Ku\oplus\mathbb
Kw)^\perp,T)$ by the one dimensional algebra with zero product by
means of the admissible pair $(D,x_0)$ where$D({\cal
V}')\subset{\cal V}'$ and $x_0\in{\cal V}'$. Recall that, ${\cal
W}={\cal U}'\oplus{\cal V}'$ where ${\cal U}'={\cal U}\cap{\cal
W}$ and ${\cal V}'={\cal V}\cap(\mathbb Ku)^\perp$. Since ${\cal
U}$ and ${\cal V}$ are stable by $\Delta$, then ${\cal U}'$ et
${\cal V}'$ are stable by $\delta$ too. Moreover,
$\Delta(a)=a_0-\lambda a$. Thus, $a_0-\lambda a\in{\cal V}$
because $a\in{\cal V}$. Hence, $a_0\in{\cal V}$. On the other
hand,
$$B(\Delta(a),b)=-B(a,\Delta(b))=-\lambda B(a,b)=B(-\lambda a,b).$$
Thus, $B(a_0,b)=0$. It follows that $a_0\in(\mathbb Kb)^\perp$.
Hence, $a_0\in{\cal V}'$. We conclude by the Theorem $\ref{232}$
that, $({\frak J}={\cal U}\oplus{\cal V},B,\omega)$ is the
symplectic-Manin double extension  of a symplectic Jordan-Manin
algebra $({\cal W}={\cal U}'\oplus{\cal V}',T,\Omega)$ by means of
$(D,x_0,a_0,\lambda)$. $\Box$
\section {Characterizations of semisimple Jordan algebras}\label{nou}

In this section we are going to give some new charcterizations of semisimple Jordan algebras among the pseudo-euclidean Jordan algebras.
 
\subsection{Characterization via the  representations of Jordan algebra}\label{rep2}

Recall that, by Theorem of Albert $(\cite{scha}$, Theorem 4.5),  a Jordan algebra ${\frak J}$ is   semisimple if and only if  its Albert form   ${\frak A}:{\frak J}\times{\frak J}\rightarrow\mathbb K$ defined by: ${\frak A}(x,y)=tr(R_{xy}),\,\,\,\forall x,y\in{\frak J}$ is non-degenerate. Now, we are going to give an answer to the following question: Let  ${\frak{J}}$ be a Jordan algebra and  $\pi: {\frak{J}}\rightarrow \mbox{End}({\cal V})$ be  a representation of finite dimension of  ${\frak{J}}.$ Consider  the bilinear form  $B_{\pi}:  {\frak{J}}\times {\frak{J}}\rightarrow \mathbb{K}$ defined by: $B_{\pi}(x,y)= tr\Bigl(\pi(xy)\Bigl), \forall x,y \in{\frak{J}}$ (we say that  $B_{\pi}$ is the biliner form of ${\frak{J}}$ associate to representation $\pi).$ What about the structure of Jordan algebra  ${\frak{J}}$ such that $B_{\pi}$ is non-degenerate?
 
\begin{prop}\label{syin}
Let   ${\frak{J}}$ be a Jordan  algebra and  $\pi:\frak{J}\longrightarrow End(\cal{V})$ be a finite-dimensional   representation of $\frak{J}$. Then, the  bilinear form  $B:\frak{J}\times \frak{J} \longrightarrow \mathbb K $ defined  by: $$B_{\pi}(x,y)=tr\Bigl(\pi(xy)\Bigl)\,\,\,, \forall x,y\in \frak{J},$$ 
is symmetric and associative.
\end{prop}
\dem
Let $x,y,z\in \frak{J}$. $B_{\pi}(x,y)=tr\Bigl(\pi(xy)\Bigl)=tr\Bigl(\pi(yx)\Bigl)=B(y,x),$ then  $B_{\pi}$ is  symmetric. Now,
 
$$B_{\pi}(xy,z)-B_{\pi}(x,yz)= tr\Bigl(\pi\Bigl((xy)z-x(yz)\Bigl)\Bigl)=tr\Bigl(\pi\Bigl((x,y,z)\Bigl)\Bigl).$$

By Corollary $\ref{lie},$ we have 
$\pi\Bigl((x,y,z)\Bigl)=[\pi(y),[\pi(x),\pi(z)]].$ It follows that 
$tr\Bigl(\pi\Bigl((x,y,z)\Bigl)\Bigl)=tr([\pi(y),[\pi(x),\pi(z)]])=0,$
so  $B_{\pi}(xy,z)-B_{\pi}(x,yz)=0.$ Consequently, $B$ is associative. $\Box$\\
 
 The following Lemma is the first step to study the non-degeneracy of $B_{\pi}.$  

\begin{lema}
 Let  $ \frak{J}$ be a Jordan Algebra. Then, $\frak{J}$ is nilpotent if and only if  $\pi(x)$ is  nilpotent, for all finite-dimensional representation    $\pi$  of  $\frak{J}$  and for all   $x\in \frak{J}$.
\end{lema}
\dem
Let  $ \frak{J}$ be a Jordan Algebra and let  $\pi :\frak{J} \longrightarrow End(\cal{V}) $   be a finite-dimensional representation  of $\frak{J}$. Then the vector space $\tilde{\frak{J}}=\frak{J}\oplus \cal{V}$ with the following product 
\begin{eqnarray*}
(x+v)(y+w)=xy+\pi(x)w+\pi(y)v, \forall x,y \in \frak{J}, v,w \in \cal{V},
\end{eqnarray*} 
is a Jordan Algebra. It is clear that  $\cal{V}$ is an ideal of $\tilde{\frak{J}} $ such that $vw=0,\,\,\, \forall v,w\in \cal{V}$, then  $\cal{V}$ is an ideal nilpotent of  $\tilde{\frak{J}}$. Consequently, $\cal{V}$ is contained in the radical  $Ra(\tilde{\frak{J}})$ of $\tilde{\frak{J}}.$

Suppose that  $\tilde{\frak{J}}$ is not nilpotent,  then, $\tilde{\frak{J}}=S\oplus Ra(\tilde{\frak{J}})$ where $ S\neq \{0\}$  is a  semi-simple subalgebra of $\tilde{\frak{J}}$. Consequently,  the Jordan  algebras $\tilde{\frak{J}}/\cal{V} $  and   $ S\oplus Ra(\tilde{\frak{J}}/ \cal{V})$ are isomorphic. 
Let $\phi:\tilde{\frak{J}}\rightarrow \tilde{\frak{J}}/\cal{V}$ be the  canonical surjection. Since $S\cap {\cal V}=\{0\},$ then  $\phi:S\longrightarrow \phi(S) $ is an isomorphism of Jordan Algebras. It follows that $\phi(S)$ is semi-simple. 
Moreover ${\tilde\frak{J}}/\cal{V}$ is nilpotent because   $\frak{J} $ is nilpotent, then   $\phi(S)$ is also nilpotent, so  $\phi(S)=\{0\}. $ Therefore, $S\subset \cal{V},$ which contradicts  $S \neq \{0\}$. We conclude that  $\tilde{\frak{J}} $ is nilpotent. consequently, $\pi(x)$ is  nilpotent for all   $x\in \frak{J}$.  

Conversely,  we have, in particular,  $R_x$ is nilpotent for all $x\in \frak{J}$ where $R$ is the adjoint representation of $\frak{J}$, so $\frak{J}$ est nilpotente  ( $\cite{scha},$  Theorem $4.3)$.$\Box$

 \begin{theor}\label{theor} Let $\frak{J}$ be a Jordan algebra. Then the following assertions are equivalent:
\begin{enumerate}
\item [(i)] $\frak{J}$ is semi-simple; 

\item [(ii)] There exists a finite-dimensional representation   $\pi:\frak{J}\times \frak{J} \rightarrow End(\cal{V}) $ of $\frak{J}$ such that the bilinear form $B: \frak{J}\times \frak{J} \rightarrow {\Bbb K}$  defined by $B(x,y)=tr(\pi(xy)), \forall x,y \in \frak{J},$ is non-degenerate.
\end{enumerate}
\end{theor} 
 
\dem
 
If  $\frak{J}$ is  semi-simple, then  by Theorem of Albert ($\cite{scha}, Theorem 4.5)$, the bilinear form  $B: \frak{J}\times \frak{J} \rightarrow {\Bbb K}$ defined by : $B(x,y)=tr(R_{xy}),\forall x,y \in \frak{J},$  is non-degenerate.    

Conversely, suppose that there exists a representation   $\pi:\frak{J}\times \frak{J} \longrightarrow End(\cal{V}) $ such that the bilinear form $B$ de $\frak{J}$ defined by:  $B(x,y)=tr(\pi(xy)),\forall x,y \in \frak{J},$  is non-degenerate. Let  $x\in {\frak{J}} \mbox{ and } r \in Ra(\frak{J})$, where $Ra(\frak{J})$ is the radical de $\frak J$.  The previous Lemma implie that $\pi(rx)$ is  nilpotent, hence $B(r,x)=tr(\pi(rx))= 0.$ The fact that  $B$ is non-degenerate implies that $r=0.$ Which proves that $Ra(\frak{J})=\{0\}.$ We conclude that  $\frak J$ is  semi-simple. $\Box$
 
%%%%%%%%%%%%%%%%%%%%%%%%%%%%%%%%%%%%%%%%%%%%%%%%%%%%%%%%%%%%%%%%%%%%%%%%%%%%
%%%%%%%%%%%%%%%%%%%%%%%%%%%%%%%%%%%%%%%%%%%%%%%%%%%%%%%%%%%%%%%%%%%%%%%%%%%%%%%%%%%%
%%%%%%%%%%%%%%%%%%%%%%%%%%%%%%%%%%%%%%%%%%%%%%%%%%%%%%%%%%%%%%%%%%%%%%%%%%%%%%%
\subsection{Opertors of  Casimir type of  pseudo-euclidean Jordan algebra}\label{cass}

Now, we are going to give a characterization of semi-simple Jordan algebra by using an operator called opertor of  Casimir type.  
Let  $({\frak J},B)$ be a  pseudo-euclidean Jordan algebra of dimension $n$.  
We consider $\beta=\{e_1,...,e_n\},\, \beta'=\{e'_1,...,e'_n\}$, two basis of  ${\frak J}$ such that  $B(e_i,e'_j)=\delta_{ij},\,\, \forall i,j\in \{1,...,n\}$ where  $\delta_{ij}$ is  Kronecker's symbol. Let us consider  $c=\sum_{i=1}^{n}e_{i}e'_{i},$ then  the operator $R_c:{\frak J}\rightarrow {\frak J}$ is defined by: $R_c(x)=\sum_{i=1}^{n} x(e_{i}e'_{i}), \forall x \in  {\frak J}.$ If   $\Gamma=\{f_1,...,f_n\},\,\,\, \Gamma'=\{f'_1,...,f'_n\}$   two other basis of ${\frak J}$ such that  \\$B(f_i,f'_j)=\delta_{ij}\,\, \forall i,j\in\{1,...,n\},$ in the following Lemma, we   prove in particular that    $R_c=R_{c'}$, where $c'=\sum_{i=1}^{n}f_{i}f'_{i}.$

\begin{lema}\label{cas}
Let  $({\frak J},B)$ be a pseudo-euclidean Jordan algebra.
\begin{enumerate}
\item[(1)]For all  $x,y\in {\frak J}$, we have ${\frak A}(x,y)=B(R_c(x),y)$ (where ${\frak A}$ is the Albert form of ${\frak J}$);
\item[(2)]  $R_c=R_{c'}$;
\item[(3)]  For all  $x\in {\frak J},$ on a $R_c\circ R_x=R_x\circ R_c$.
\end{enumerate}
\end{lema}
\dem
(1) Let $x,y \in{\frak J}$, $B(R_c(x),y)= \sum_{i=1}^{n}B(R_{e_{i}e'_{i}}(x),y)
=\sum_{i=1}^{n}B({e_{i}e'_{i}},xy)= \sum_{i=1}^{n}B((xy){e_{i},e'_{i}})=trR_{xy}={\frak A}(x,y).$
 
(2) By $(1)$, it is clear that  $B(R_{c'}(x),y)={\frak A}(x,y)=B(R_c(x),y),\,\,\forall x,y\in{\frak J}$. The fact that  $B$ is non-degenerate  implies that $R_c=R_{c'}$. 
 
(3) Let $x,y,z\in {\frak J}$,   $B(R_cR_x(y),z)-B(R_xR_c(y),z)=B(R_c(yx),z)-B(R_c(y),xz)= {\frak A}(yx,z)-{\frak A}(y,xz)= 0.$ The fact that $B$ is  non-degenerate implies that $R_cR_x(y)=R_xR_c(y).$ $\Box$
  
\begin{defi}
  $R_c$ is called the operateur of Casimir type of the pseudo-euclidean Jordan algebra $({\frak J},B)$.
\end{defi}

\begin{prop}\label{cas1}
 ${\frak J}$ is  semi-simple if and only if  $R_c$ is  invertible.
\end{prop}
\dem Let  $x,y\in {\frak J}$, by Lemma \ref{cas} we have 
 ${\frak A}(x,y)=B(R_c(x),y)$. Consequently,  the Albert form $\frak A$ is non-degenerate if and only if  the operator  $R_c$ est invertible. Then we conclude, by Theorem 4.5 of \cite{scha},  that ${\frak J}$ is  semi-simple if and only if  $R_c$ is  invertible. $\Box$
  
\begin{prop}\label{cas2}
Let  $({\frak J},B)$ be an  $B-$irreducible pseudo-euclidean Jordan algebra. Then  the operator of Casimir Type $R_c$ of  $({\frak J},B)$ is either  nilpotent or  invertible.
\end{prop}
\dem
It is known that   there exists  $n\in \mathbb N$ such that   ${\frak J}=ker(R_c)^n\oplus Im(R_c)^n$. It is easy to see that $B(ker(R_c)^n,Im(R_c)^n)=\{0\},$ so   ${(ker(R_c)^n)}^{\perp}=Im(R_c)^n.$ Moreover, by Lemma $\ref{cas}$,  $ ker(R_c)^n$  is an ideal of    ${\frak J}$. Consequently,  $ker(R_c)^n=\{0\}$ or $ker(R_c)^n={\frak J}$. We conclude that $R_c$  is either  nilpotent or  invertible. $\Box$
 
\begin{prop}
Let  $({\frak J},B)$ be a pseudo-euclidean Jordan algebra and  $R_c$ be the operator of Casimir Type of  $({\frak J},B)$. Then  $R_c$ is  nilpotent if and only if  ${\frak J}$ is a Jordan algebra without nonzero semisimple ideal.
\end{prop}
\dem   Suppose that  ${\frak J}$ contain a semisimple $\cal I$. By Lemma $\ref {lem1},$ $\cal I$ and ${\cal I}^{\perp}$ are  non-degenerate ideals of ${\frak J}$   and ${\frak J}={\cal I}\oplus{\cal I}^{\perp}$. It follows that $R_{c_{\vert_{{\cal I}}}}$ is the operateur of Casimir type   of $({\cal I}, B_{\vert_{{\cal I}\times{\cal I}}})$ where $R_c$ is the operator of Casimir type of $({\frak J},B)$. If $ {\cal I}\not= \{0\},$ then $R_{c_{\vert_{{\cal I}}}}$ is invertible, so $R_c$ is not nilpotent.

Conversely, suppose that ${\frak J}$ is without nonzero semisimple ideal then ${\frak J}$ is an orthogonal direct sum of irreducible non-degenerate ideals ${\frak J}_{k},\, k\in \{1,\dots,p\}.$ Let us set $R_k$ the operator of Casimir type of the pseudo-euclidean Jordan algebra $({\frak J}_{k},B_{\vert_{{\frak J}_{k}\times{\frak J}_{k}}}),$ where $k\in \{1,\dots,p\}.$ It is clear that $R_{c_{\vert_{{\frak J}_{k}}}}= R_k, $ for all $k\in \{1,\dots,p\}.$ Since, by Propositions $\ref{cas1}$ and $\ref{cas2}$,  $R_k$ is nilpotent, for all $k\in \{1,\dots,p\},$ then $R_{c}$ is nilpotent. $\Box$

\begin{coro}
Let  $({\frak J},B)$ be a pseudo-euclidean Jordan Algebra and $R_{c}$ be its operator of Casimir type. Let  $S$ be the largest semisimple ideal of ${\frak J}.$  Then ${\frak J}=S\oplus S^\perp$ is the    Fitting decomposition of ${\frak J}$ relative to $R_{c}$.
\end{coro}
\dem By Lemma $\ref {lem1}$ $S$ and $S^{\perp}$ are  non-degenerate ideals of ${\frak J}$   and ${\frak J}={S}\oplus{S}^{\perp}$. It follows that the operator of Casimir type of $(S,B_{\vert_{S\times S}})$  (resp. $(S^\perp,B_{\vert_{S^\perp\times S^\perp}})$) is $T:={(R_c)}_{\vert_{S}}\mbox{ (resp. } 
L:={(R_c)}_{\vert_{S^\perp}}).$ Moreover,  $T$ is invertible because $S$ is semisimple and $L$ is nilpotent because ${S^\perp}$ is a Jordan algebra without nonzero semisimple ideal. $\Box$

%%%%%%%%%%%%%%%%%%%%%%%%%%%%%%%%%%%%%%%
%%%%%%%%%%%%%%%%%%%%%%%%%%%%%%%%%%%%%%%%%%%%

\subsection{Characterization by means of  the index of Jordan algebra}
In this subsection, ${\mathbb K}= {\mathbb R}\, or \,{\mathbb C}.$
Let $({\frak J},B)$ be a   pseudo-euclidean Jordan algebra.  We denote by  $\frak F({\frak J})$ the vector space  of all associative symmetric bilinear forms on ${\frak J}$ and by  $\frak B({\frak J})$ the subspace of $\frak F({\frak J})$ spanned by the set of the associative scalar products  on ${\frak J}$.
\begin{lema}
 $\frak B({\frak J})=\frak F({\frak J})$.
\end{lema}
\dem
Let  $B$ be  an associative scalar product on  ${\frak J}$ and  $\varphi\in \frak F({\frak J})$. Let  $T$ be a basis of ${\frak J}$, $M(B)$ and  $M(\varphi)$ the associated matrices of $B$  and $\varphi$ in  $T$. Since  $B$ is non-degenerate, then there exists  $f$, a linear map of ${\frak J}$ into itself, such that  $\varphi(x,y)=B(f(x),y),\,\forall x,y\in {\frak J}.$ Let  $M(f)$ be the matrix of   $f$ in $T$, so  $M(\varphi)={^tM}(f)M(B).$ Consider the polynomial  $P(X)=det(M(\varphi)-XM(B))\in 
\mathbb C[X]$, then  $P(X)= det(M(B))det(M(f)-XI_n)$  where  $n$ is the dimension of  ${\frak J}.$ Consequently,  $P(X)$ is a non-zero polynomial; so there exists    $\lambda\in {\mathbb K}$ such that $P(\lambda)\neq 0$. Therefore   $\varphi-\lambda B$ is  non-degenerate, it follows that $\varphi=(\varphi-\lambda B)+\lambda B\in \frak B({\frak J}).$ We conclude that  $\frak F({\frak J})=\frak B({\frak J})$. $\Box$
  
\begin{defi}
Let  $({\frak J},B)$ be a pseudo-euclidean Jordan algebra. The dimension of  
${\frak B}({\frak J})$  is called the index of ${\frak J}$ and will be denoted by $ind({\frak J})$.
\end{defi}
%%The followin proposition is well known $(\cite{far}).$
\begin{prop}\label{sind}
If  ${\frak J}$ is a simple Jordan algebra over ${\mathbb C}$, then  $ind({\frak J})=1$.
\end{prop}
\dem
Since  ${\frak J}$ is simple then,  by Theorem of Albert (\cite{scha}, Theorem 4.5),   the Albert form   ${\frak A}$ of ${\frak J}$ is an associative scalar product on ${\frak J}$. Now, let $B$ be an associative scalar product on  ${\frak J}$. Then there exists   an endomorphism  $D$ of ${\frak J}$ such that $B(x,y)={\frak 
A}(D(x),y),\,\forall x,y\in{\frak J}.$ If  $\lambda\in {\mathbb K}$, the bilinear form  $B'$ of ${\frak J}$ defined by:  $B'(x,y)={\frak A}(D(x)-\lambda x,y),\,\forall x,y\in {\frak J},$ is  
associative. If  $\lambda$ is  an eigenvalue of  $D$ and  $v$ is an eigenvector for $\lambda$ , then  $B'(v,y)=0,\,\forall y\in {\frak J}$ and  
$rad(B'):=\{x\in {\frak J};B'(x,y)=0,\,\,\forall y\in {\frak J}\}\neq \{0\}.$ Moreover, $rad(B')$ is an ideal of  ${\frak J}$, it follows that $rad(B')={\frak J}$. Therefore $D-\lambda id_{{\frak J}}=0$, so $B(x,y)=\lambda{\frak A}(x,y),\,\forall x,y\in {\frak J}.$ $\Box$
 
\begin{prop}\label{11}
Let $({\frak J},B)$ be a pseudo-euclidean Jordan algebra. If  $ind({\frak J})=1$, then  ${\frak J}$ is either a simple Jordan algeba   or ${\frak J}$ is the one-dimensional algebra with zero product.
\end{prop}
 
\dem  If  $ind({\frak J})=1$, then  ${\frak J}$ is irreducible. Assume that  ${\frak J}$ is neither simple  nor the one-dimensional algebra with zero product. Then, by Corollary $\ref{dx1}$ and Theorem $\ref{th3}$, ${\frak J}$ is either a double extension of a pseudo-euclidean Jordan algebra by a simple Jordan algebra or a generalized double extension  of a pseudo-euclidean Jordan algebra by the one-dimensional algebra with zero product.  
If ${\frak J}^2={\frak J}$ (i.e. $Ann({\frak J})=\{0\}$), then  ${\frak J}$ is a double extension of a pseudo-euclidean Jordan algebra $({\cal W},T)$ by a simple Jordan algebra   $\cal S$. Moreover,  by Theorem $\ref{th1}$, if $\gamma$ is an associative symmetric bilinear form on  ${\cal S}$, then the bilinear form $\tilde \gamma$ on ${\cal S}\oplus{\cal W}\oplus{\cal S}^*$ defined by: ${\tilde \gamma}
(x+y+f,x'+y'+f')= \gamma(x,x')+B(y,y')+f(x')+f'(x),\, \forall (x,y,f),(x',y',f') \in {\cal S}\times {\cal W}\times {\cal S}^*,$ 
 is an associative scalar product on ${\frak J}.$ Let us consider $\gamma_1={\frak A}$ the Albert form of ${\cal S}$ andthe bilinear form $\gamma_2=0$  on ${\cal S}$ . It is clear that  $\tilde{\gamma_1}$ et $\tilde{\gamma_2}$ are two linearly independant elements of    ${\frak B}({\frak J})$, so $ind({\frak J})\geq2$ which contradicts the fact that $ind({\frak J})= 1$. Now, if  $Ann({\frak J})\neq\{0\}$, then there exists   $b\in Ann({\frak J})$ such that $b\neq 0$. Since ${\frak J}$ is irreducible and  ${\frak J}$ is not the one-dimensional algebra with zero product, then $B(b,b)= 0.$ By Theorem    $\ref{th3}$, ${\frak J}= {\mathbb K}a\oplus{\cal W}\oplus\mathbb Kb$, where $a\in{\frak J}$ such that  $B(a,a)=0$, $B(a,b)=1$ and ${\cal W}^{\perp}={\mathbb K}a\oplus {\mathbb K}b$ .   Let us consider the associative symmetric bilinear form $T$ on ${\frak J}$ defined by: $T(a,a)=1,\,$ and $\,T(x,y)=0,\,\forall (x,y)\in{\frak J}\times{\frak J}\setminus{\mathbb K}a\times{\mathbb K}a.$ It is easy to see that $B$ et $T$ are two linearly independant elements of    ${\frak B}({\frak J}).$ therefore, $ind(\frak{J})\geq 2$, which contradicts the hypothesis $ind({\frak J})= 1.$
 
We conclude that if $ind({\frak J})=1$, then ${\frak J}$ is either a simple Jordan algebra or  ${\frak J}$ is the one-dimensional algebra with zero product. $\Box$\\

For any real Jordan algebra ${\frak J}$, we will denote by ${\frak J}^{\mathbb C}:= {\frak J}{\otimes_{\mathbb R}}{\mathbb C}$ its complexification which is a Jordan algebra over ${\mathbb C}$. It is clear if $({\frak J},B)$ is a real pseudo-euclidean Jordan algebra, then $({\frak J}^{\mathbb C},B^{\mathbb C})$ is a complex pseudo-euclidean Jordan algebra where $B^{\mathbb C}$ is the bilnear form on defined ${\frak J}^{\mathbb C}$ defined by: $B^{\mathbb C}(x\otimes a,y\otimes b):= ab B(x,y), \forall x,y \in {\frak J}, \forall a,b \in 
{\mathbb C}.$

\begin{coro}\label{am1}
Let  $({\frak J},B)$ be a pseudo-euclidean Jordan algebra   over $\mathbb C$. Then the two following assertions are equivalent:
\begin{enumerate}
\item $ind({\frak J})=1$
\item ${\frak J}$ is either a simple Jordan algebra   or ${\frak J}$ is the one-dimensional algebra with zero product. 
\end{enumerate}
\end{coro}
\begin{coro}
Let  $({\frak J},B)$ be a pseudo-euclidean Jordan algebra   over $\mathbb R$ such that ${\frak J}^2\neq\{0\}$. Then the  following assertions are equivalent:  
\begin{enumerate}
\item[(i)] $ind({\frak J})=1$,
\item[(ii)] $ind({\frak J}^{\mathbb C})=1$ (i.e. dim$_{\mathbb C}B({\frak J}^{\mathbb C})=1$),
\item[(iii)]${\frak J}^{\mathbb C}$ is simple.
\end{enumerate}
\end{coro}
\dem By Corollary $\ref{am1}$, the assertions  $(ii)$ and $(iii)$ are equivalent. In the following we are going to prove that  $(i)$ and  $(ii)$ are equivalent. Assume that   dim$_{\mathbb C}B({\frak J}^{\mathbb C})=1$. Let  $T_1,$ $T_2$ be two associative scalar products on   $\frak J$. For $k\in\{1,2\}$, we consider the symmetric bilinear form   $\tilde{T_k}$  on  ${\frak J}^{\mathbb C}$ defined by:
$$\tilde T_k(x,y)=T_k(x,y),\hskip 0,5 cm \tilde T_k(ix,iy)=-T_k(x,y),\hskip 0,5 cm \tilde T_k(ix,y)=iT_k(x,y),\,\,\forall x,y\in{\frak J}.$$

It is easy to verify that  $\tilde{T_k} $ is an associative scalar product on   ${\frak J}^{\mathbb C}$. Since dim$_{\mathbb C}B({\frak J}^{\mathbb C})=1,$ then there exists   $\lambda:=\alpha+i\beta\in {\mathbb C}$ such that   $\tilde{T_1}=\lambda \tilde{T_2}$. In particular,  for all  $x,y\in{\frak J}$  we have
$\tilde{T_1}(x,y)=\lambda\tilde {T_2}(x,y), $ so $T_1(x,y)=(\alpha+i\beta)T_2(x,y)$. It follows that  $\beta=0,$ consequently $T_1=\alpha T_2$. We conclude that
dim$_{\mathbb R}B(\frak J)=1.$ 

Conversely, assume that dim$_{\mathbb 
R}B(\frak J)=1.$ By Proposition $\ref{11}$, 
$\frak J$ is  simple. Suppose that dim$_{\mathbb C}B({\frak J}^{\mathbb C})\neq 1.$ Then, by Theorem $\ref{am1}$, ${\frak J}^{\mathbb C}$ is not 
simple. Therefore $\frak J$ admits a complexe structure. This means that there exists a simple complex Jordan algebra  $\frak S$ such that  ${\frak J}={\frak S}_{\mathbb  R}$ is the underlying real algebra $\frak S.$ By Theorem 8.5.2 of $\cite{far}$, $\frak S$ admits an euclidean real form  ${\cal E}$, so $ {\frak J}= {\cal E}\oplus i{\cal E} $. Let $\frak A$ be the Albert form of ${\cal E},$ we consider the two bilinear forms $B_1$ and  $B_2$ on $\frak J$ defined by:
 
\begin{eqnarray*}
&&B_1(x_1+ix_2,y_1+iy_2)={\frak A}(x_1,y_1)-{\frak A}(x_2,y_2),\\
&&B_1(x_1+ix_2,y_1+iy_2)={\frak A}(x_1,y_2)+{\frak A}(x_2,y_1),\,\,\,\,\forall x_1,x_2,y_1,y_2\in {\cal E}.
\end{eqnarray*} 
An easy calculation proves that  $B_1$ and  $B_2$ are linearly independent associative scalar product on ${\frak J}$ which contradicts that $ind({\frak J})=1.$ We conclude that  dim$_{\mathbb C}B({\frak J}^{\mathbb C})=1.$ $\Box$
 
%%%%%%%%%%%%%%%%%%%%%%%%%%%%%%%%%%%%%%%%%%
%%%%%%%%%%%%%%%%%%%%%%%%%%%%%%%%%%%%%%%%%%%%%%%%%%%%%
%%%%%%%%%%%%%%%%%%%%%%%%%%%%%%%%%%%%%%%%%%%%
 
\begin{defi}
Let  ${\frak J}$ be a  Jordan algebra.   ${\frak J}$ is called  reductive if ${\frak J}={\frak S}\oplus Ann({\frak J})$ where $\frak S$ is a  semisimple ideal of  ${\frak J}$.
\end{defi}
\begin{rema}
If ${\frak J}=\sigma\oplus Ann({\frak J})$ where $\frak S$ is a  semisimple ideal of  ${\frak J}, $
then $\frak S$ is the greatest semisimple ideal of $\frak J$ and ${\frak S}=\bigoplus^{n}_{i=1}{\frak S}_i$ where $\{ {\frak S}_i,\,1\leq i\leq n\}$ is the set of all simple ideals of $\frak J.$
 \end{rema}

\begin{coro}
Let  $({\frak J},B)$ be a pseudo-euclidean Jordan algebra and   ${\frak J}_1,...,{\frak J}_r$ ($r\in\mathbb N$) be $B-$irreducible non-degenerate ideals  of  ${\frak J}$ such that ${\frak J}=\bigoplus^{r}_{i=1}{\frak J}_i$ and  $B({\frak J}_k,{\frak J}_l)= \{0\},$ for all $k,l \in \{1,\dots,r\}.$  Then the following assertions are equivalent:
\begin{enumerate}
\item  $ind({\frak J})=r.$
\item ${\frak J}$ is reductive and $dim Ann({\frak J})\leq1$.
\end{enumerate}
\end{coro}
\dem
Suppose that  $ind({\frak J})=r$. Let $i\in \{1,...r\}$, we consider  $\{T_{i1},...,T_{in_i}\} $ a  basis of  ${\frak B}({\frak J}_i)$, where $n_i=ind({\frak J}_i)$. Now, for  $i\in\{1,...,r\}$ and $j\in\{1,...,n_i\}$, we consider the bilinear form $\tilde{T_{ij}}:{\frak J}\times{\frak J}\rightarrow {\mathbb C}$ defined by: $\tilde{T_{ij}}_{{\vert}_{{\frak J}_i\times{\frak J}_i}}=T_{ij},\,\mbox{and} \, \tilde{T_{ij}}(x,y)=0,\,\forall (x,y)\in{\frak J}\times{\frak J}\setminus{\frak J}_i\times{\frak J}_i.$  It is clear that the elements of 
$\bigcup_{1\leq i\leq r}\{\tilde{T_{ij}},\,\,1\leq j\leq n_i\}$ are linearly independant, so   $ind({\frak J})\geq\sum_{1\leq i\leq r}n_i=\sum_{1\leq i\leq r}ind({\frak J}_i).$ Since $ind({\frak J})=r$, then  $ind({\frak J}_i)=1,\,\forall i\in\{1,...,r\}$ (i.e. ${\frak B}({\frak J}_i)=Vect\{{T_{i1}}\})$. Consequently, by Corollary $\ref{am1}$,for all  $i$ in  $\{1,...,r\}$, ${\frak J}_i$ is either simple or  the one-dimensional algebra with zero product. Let us suppose that there exist $i\neq j\in\{1,...,r\}$ such that  ${\frak J}_i={\mathbb C}a$ and  ${\frak J}_j={\mathbb C}b$ with $a^2=b^2=0$ and  consider the symmetric bilinear form defined on  ${\frak J}$ by: $T(a,b)=1\mbox{ and } T(x,y)=0,\,\,\,\forall (x,y)\in{\frak J}\times{\frak J}\setminus({\mathbb C}a\times{\mathbb C}b).$ It is clear that  $T$ is an element of  ${\frak F}({\frak J})$. Therefore, there exist $\alpha_1,...,\alpha_r\in {\mathbb C}$ such that $T=\sum_{i=1}^{r}\tilde{T_{i1}}$. Which contradicts the fact that $T(a,b)\neq 0$,  so  $\frak J$ is  reductive with $dim Ann({\frak J})\leq 1$.

Conversely, Assume that ${\frak J}$ is reductive and  $dim Ann({\frak J})\leq1$. Then, without lost of generality, we can suppose that for all  $i\in\{2,...,r\}$ 
${\frak J}_i$ is  simple and  ${\frak J}_1$ is either simple or  the one-dimensional algebra with zero product. It follows that for all  $i\in\{1,...,r\}, $ we have $ind({\frak J}_i)=1$ and  ${\frak B}({\frak J}_i)=\mathbb CB_i$, where $B_i:=B_{{\vert}_{{\frak J}_i\times{\frak J}_i}}.$  Now, If   $T$ is  an associative symmetric bilinear form on ${\frak J},$ then $T_i:={T_{\vert}}_{{\frak  J}_i\times{\frak J}_i}$ is  an associative symmetric bilinear form on  ${\frak J}_i$. Consequently, there exists ${(\alpha_i)}_{1\leq i\leq r}$ in $\mathbb C$ such that  
$T_i=\alpha_iB_i,\,\forall i\in\{1,...,r\}$. Moreover, since $T_i({\frak J}_i,{\frak J}_k)=0,\,\forall i\neq k\in\{1,...,r\},$ then  $T=\sum^{r}_{i=1}\alpha_i{\tilde B}_i$ where $\tilde B_i$ is the bilinear form on ${\frak J}$ defined by: ${\tilde {B_i}}_{\vert_{{\frak 
J}_i}\times{{\frak J}_i}}=B_i \mbox{ and } {\tilde B_i}(x,y)=0,\,\forall x,y\in 
{\frak J}\times{\frak J}\setminus{{\frak J}_i\times{\frak J}_i}.$ Moreover, it is easy to see that the elements of  $\{\tilde{(B_i)},\,{1\leq i\leq r}\}$ are linearly independent. We conclude that $ind({\frak J})=r.$ $\Box$

In the  following result, we give a characterization of semisimple Jordan algebra by using the notion of index

\begin{coro}
Let  $({\frak J},B)$ be a pseudo-euclidean Jordan algebra and  ${\frak J}_1,...{\frak J}_r$ ($r\in\mathbb N$) be $B-$irreducible non-degenerate ideals  of  ${\frak J}$ such that ${\frak J}=\bigoplus^{r}_{i=1}{\frak J}_i$ and  $B({\frak J}_k,{\frak J}_l)= \{0\},$ for all $k,l \in \{1,\dots,r\}.$  Then the two following assertions are equivalent: 
\begin{enumerate}
\item${\frak J}$ is semisimple,
\item ${\frak J}^2={\frak J}$ and $ind({\frak J})=r.$
\end{enumerate}
\end{coro}

{\bf Acknowledgments} \\

 We thank W. Bertram, M. Bordemann and A. Elduque for very interesting discussions. Moreover, we are very grateful to A. Elduque for his remarks, which improved the readability of this article.


\begin{thebibliography}{99}
\bibitem{aub} A. Aubert, Structures affines et pseudo-m\'etriques invariantes \`a gauche sur des groupes de Lie, Th\`ese, Universit\'e Montpellier II, 1996.
\bibitem{ben2} I. Bajo, S. Benayadi, A. Medina. $``$Symplectic structures on quadratic Lie algebras$"$, Journal of Algebra 316 (1) (2007) 174--188.
 
\bibitem{bert} W. Bertram, The geometry of Jordan and Lie
structures, Lecture Notes in Math., Vol. 1754,  Springer-Verlag,
Berlin, 2000.
\bibitem{Bord}  M. Bordemann, ``Nondegenerate invariant bilinear forms on nonassociative algebras,   Acta Math. Univ. Com  LXVI (2) (1997) 
151--201.
\bibitem{bou} N. Bourbaki, Groupe et alg\`ebres de Lie. Ch. 7 et 8, Hermann, Paris, (1975).
\bibitem{med2} A. Diatta, A. Medina, $``$Classical Yang-Baxter equation and left invariant affine geometry on Lie groups$"$, Manuscr. Math 114  (4) (2004) 477--486.
\bibitem{far} J.Faraut and  A. Kor\'anyi,  Analysis on Symmetric Cones, Clarendon Press, Oxford,  (1994).
\bibitem{Jacob} N. Jacobson,  Structure and representations of Jordan algebras, American Mathematical Society Colloquium Publications, Vol.XXXIX.  American Mathematical Society, Providence, R.I.(1968).
\bibitem{koe1}M. Koecher, Imbedding of Jordan algebras into Lie algebras. I.  Amer. J. Math  89  (1967) 787--816.
\bibitem{koe2}M. Koecher, Imbedding of Jordan algebras into Lie algebras. II. Amer. J. Math  90  (1968) 476--510.
\bibitem{mac} K. McCrimmon, A taste of Jordan algebras. Universitext. Springer-Verlag, New York, (2004).
\bibitem{med}  A. Medina, Ph. Revoy, Alg\`ebre de Lie et produit scalaire invariant, Ann. Scient. Ec. Norm. Sup., 4\`eme s\'erie, 18  (1985) 553--561.
\bibitem{nei} O'Neill, Semi-Riemannian geometry with applications to relativity, Academic Press, New York, (1983).
\bibitem{scha}   R.D. Schafer, An introduction to nonassociative algebras, Academic Press, New York, (1966).
\bibitem{Tits} J. Tits, Une classe d'alg\`ebres de Lie en relation avec les alg\`ebres de Jordan, Nederl. Akad. Wetensch. Proc. Ser.A 65 = Indag. Math 24  (1962)  530--535.
\bibitem{zhel} V.N Zhelyabin. $``$Jordan D-Bialgebras And Symplectic Forms On Jordan Algebras$"$, Siberian Adv. Math   10   (2) (2000) 142--150.
\bibitem{zhel1} V.N Zhelyabin. $``$On a class of Jordan $D$-bialgebras$"$.
St. Petersbg. Math. J  11 (4) (2000) 589--609 ; translation from
Algebra Anal  11  (4) (1999) 64--94.
\end{thebibliography}
\end{document}